\documentclass[11pt]{amsart}
\usepackage{amsmath}
\usepackage{amssymb}
\usepackage{amsthm}
\usepackage{amscd}

\newcommand{\mpc}[1]{\marginpar{}}
\newcommand{\npc}[1]{\marginpar{}}

\numberwithin{equation}{section}

\newcounter{AbcT}

\newtheorem {Theorem}    {Theorem}[section]
\newtheorem* {Definition} {Definition}

\newtheorem {Lemma}      [Theorem]    {Lemma}
\newtheorem {Corollary}  [Theorem]    {Corollary}
\newtheorem {Proposition}[Theorem]    {Proposition}
\newtheorem {Claim}      [Theorem]    {Claim}

\newtheorem*{thma}{Theorem~\ref{Thmain}(a)}
\newtheorem*{thmb}{Theorem~\ref{Thmain}(b)}
\newtheorem* {TheoremA.2}    {Theorem A.2}
\newtheorem* {TheoremA.9}    {Theorem A.9}
\newtheorem* {Theorem6.1}    {Theorem 6.1}
\newtheorem* {Theorem6.2}    {Theorem 6.2}

\newcommand{\ex}[5]{ #1 \stackrel{#4}{\hookrightarrow} #2 \stackrel{#5}{\twoheadrightarrow} #3 }
\newcommand {\proofs}     {\proof}

\newcounter{DM@bibnum}

\newcommand  {\QED}    {\def\qedsymbol{$\square$}\qed}

\newcommand {\infdep} {{\rm infdep}}
\newcommand {\comdep} {{\rm comdep}}

\newcommand{\fe} { f }

\newcommand{\la}{\langle}
\newcommand{\ra}{\rangle}
\newcommand{\lla}{\langle\!\langle}
\newcommand{\rra}{\rangle\!\rangle}

\newcommand{\frc}{\displaystyle\frac}

\def\skv{{\vskip .12cm}}

\def\Exp{{\rm Exp}}
\def\Log{{\rm Log}}

\def\log{{\rm log\,}}
\def\deg{{\rm deg\,}}

\def\Aut{{\rm Aut\,}}
\def\Gal{{\rm Gal\,}}
\def\Ker{{\rm Ker\,}}
\def\Im{{\rm Im\,}}
\def\Hom{{\rm Hom\,}}
\def\Tr{{\mathbf {Tr}\,}}
\def\tr{{\mathbf {tr}\,}}
\def\Tred{{\rm T_{\rm red}\,}}

\def\Nred{{\rm N_{\rm red}\,}}

\def\id{{\rm{id}}}

\def\diag{{\rm diag\,}}

\def\NSL_2{{\mathcal N SL_2}}

\def\res{{\mathbf {res}}}

\def\sl{{\mathfrak{sl}}}

\def\lam{\lambda}            
\def\sig{\sigma}                
\def\phi{\varphi}

\def\calA{{\mathcal A}}

\def\calE{{\mathcal E}}
\def\calF{{\mathcal F}}
\def\calG{{\mathcal G}}

\def\calL{{\mathcal L}}

\def\cbar{\bar c}               
               
               \def\Ebar{\overline{E}}
               \def\Fbar{\overline{F}}
               
\def\hbar{\bar h}

               \def\Kbar{\overline{K}}
               \def\Lbar{\overline{L}}
               \def\Mbar{\overline{M}}

               \def\Wbar{\overline{W}}

               \def\Ghat{\widehat {\mathstrut G}}
               \def\Hhat{\widehat {\mathstrut H}}

               \def\Nhat{\widehat {\mathstrut N}}

               \def\Shat{\widehat {\mathstrut S}}

\def\xhat{\hat x}

 \def\gra{{\mathfrak a}}

 \def\grf{{\mathbf f}}
\def\grG{{\mathfrak G}} \def\grg{{\mathfrak g}}
 \def\grh{{\mathfrak h}}

  \def\bfk{{\mathbf k}}
\def\grL{{\mathfrak L}} \def\grl{{\mathfrak l}} \def\bfl{{\mathbf l}}
 \def\grm{{\mathfrak m}} \def\bfm{{\mathbf m}}
 \def\grn{{\mathfrak n}}
 \def\gro{{\mathfrak o}}
 
 \def\grq{{\mathfrak q}}

  \def\bfu{{\mathbf u}}
 
 \def\grw{{\mathbf w}}

 \def\grz{{\mathfrak z}}

\def\gref{\grg_{df+1}}

\def\dbC{{\mathbb C}}

\def\dbF{{\mathbb F}}
\def\dbG{{\mathbb G}}

\def\dbN{{\mathbb N}}

\def\dbQ{{\mathbb Q}}
\def\dbR{{\mathbb R}}

\def\dbZ{{\mathbb Z}}

\def\grhat{{\widehat{\grh}}}

\def\Fp{{\dbF_p}}

\newcommand{\lpar}{\left(}
\newcommand{\rpar}{\right)}

\def\str{\stackrel}
\def\lra{\longrightarrow}

\def\Ext{{\rm Ext}}

\begin{document}

\title[Second cohomology of the norm one group]
{On the second cohomology of the norm one group of a $p$-adic
division algebra}
\author{Mikhail Ershov}
\address{University of Virginia}
\email{ershov@virginia.edu}
\thanks{The author is grateful for the support and hospitality of the
Institute for Advanced Study where part of this research was carried out.
The author was supported by the National Science
Foundation under agreement No. DMS-0111298. 
}
\maketitle
\vskip -.4cm
\centerline{ with an appendix on }
\vskip .1cm
\centerline {\bf $\exp-\log$ correspondence for powerful $p$-central} 
\centerline{\bf $\dbZ_p$-Lie algebras and pro-$p$ groups}
\vskip .2cm
\centerline{by {\sc Mikhail Ershov} and {\sc Thomas Weigel}}
\vskip .5cm

\centerline{\it To Gopal Prasad on the occasion of his 75th
birthday}

\begin{abstract}
Let $F$ be a $p$-adic field, that is, a finite extension of $\dbQ_p$.
Let $D$ be a finite-dimensional central division algebra over $F$ and let $SL_1(D)$
be the group of elements of reduced norm $1$ in $D$.
Prasad and Raghunathan proved that $H^2(SL_1(D),\dbR/\dbZ)$ is a cyclic $p$-group whose order is bounded from below by
the number of $p$-power roots of unity in $F$, unless $D$ is a quaternion algebra over $\dbQ_2$.
In this paper we give an explicit upper bound for the order of $H^2(SL_1(D),\dbR/\dbZ)$ for $p\geq 5$,
and determine $H^2(SL_1(D),\dbR/\dbZ)$ precisely when $F$ is cyclotomic, $p\geq 19$ and the degree
of $D$ is not a power of $p$.
\end{abstract}

\section{Introduction}
Let $F$ be a nonarchimedean local field of residue characteristic $p$,
that is, a finite extension of $\dbQ_p$ or the field of Laurent series
over a finite field of characteristic $p$.
Let $G$ be the group of rational points of a
connected simply-connected simple algebraic group $\dbG$ defined over $F$.
By results of Moore~\cite{Mo2} and Prasad and Raghunathan~\cite{PR1},\cite{PR2},
the second continuous cohomology group
\footnote{Here $\dbR/\dbZ$ is endowed with its usual topology, and the action of $G$
on $\dbR/\dbZ$ is trivial}
 $H^2(G,\dbR/\dbZ)$ classifies 
topological central extensions of $G$ (see \cite[Chapter 10]{PR1} for a detailed discussion).
If $F$ has characteristic zero, finiteness of $H^2(G,\dbR/\dbZ)$ follows from a general theorem of Raghunathan~\cite{Ra}.
\footnote{A very short proof of this fact was later found by Prasad~\cite{Pr1}.}
However, the exact determination of the above group (for $F$ of either characteristic)
is a deeper problem which received a lot of attention since mid 60's
starting with a work of Moore \cite{Mo1} and
culminating in works of Prasad and Raghunathan \cite{PR1} and \cite{PR2}.
It is now known that if $\dbG$ is \bf{isotropic }\rm over $F$ then $H^2(G,\dbR/\dbZ)$ 
is isomorphic to the group of roots of unity in $F$. Using Moore's paper [Mo1],
Matsumoto~\cite{Ma} proved this result for $F$-split groups, and the case of 
$F$-quasi-split groups is due to Deodhar~\cite{De} and Deligne (unpublished). 
Almost all remaining cases were handled in \cite{PR1}, and finally the complete 
answer was obtained in \cite{PRp}.
\footnote{The argument in \cite{PRp} uses global fields, but later Prasad~\cite{Pr2}
found a purely local proof.}

If $\dbG$ is anisotropic over $F$, then by Tits' classification
$G$ is isomorphic to $SL_1(D)$ for some finite-dimensional central division algebra $D$
over $F$. Since $SL_1(D)$ is profinite, $H^2(SL_1(D),\dbR/\dbZ)$ is isomorphic to
$H^2(SL_1(D),\dbQ/\dbZ)$, where $\dbQ/\dbZ$ is endowed with discrete topology (see \cite[2.0]{PR2}).
Moreover, $H^2(SL_1(D),\dbQ/\dbZ)$ is isomorphic to $H^2(SL_1(D),(\dbQ/\dbZ)_p)$,
where $(\dbQ/\dbZ)_p$ is the $p$-primary component of $\dbQ/\dbZ$ (see \cite[2.1]{PR2}).
The main result of \cite{PR2} asserts that $H^2(SL_1(D),(\dbQ/\dbZ)_p)$
is a cyclic $p$-group whose order is bounded from below by the number of
$p$-power roots of unity in $F$, unless $D$ is the quaternion division algebra over $\dbQ_2$.
Moreover, $H^2(SL_1(D),(\dbQ/\dbZ)_p)$ is trivial if
$F$ has no $p$-power roots of unity (in particular, if $F$ has characteristic $p$)
and $D$ is not a quaternion division algebra over $\dbQ_3$.

The goal of this paper is to obtain an explicit upper bound
for the order of $H^2(SL_1(D),\dbR/\dbZ)$ when $F$ has characteristic zero and $p\geq 5$.
Part (a) of the following theorem provides such a bound in the general case,
parts (b) and (c) give stronger bounds in some special cases, and
part (d) gives a precise order for $H^2(SL_1(D),\dbR/\dbZ)$ in the case of
cyclotomic fields: 

\begin{Theorem}
Let $F$ be a finite extension of $\dbQ_p$,
let $e$ be the ramification index of $F$, and let
$p^w$ be the highest power of $p$ dividing $e$.
Let $D$ be a finite-dimensional central division algebra over $F$, and let
$p^N$ be the order of $H^2(SL_1(D),\dbR/\dbZ)\cong H^2(SL_1(D),(\dbQ/\dbZ)_p)$.
\begin{itemize}
\item[(a)] Assume that $p\geq 5$. Then $N\leq w+6$.
\item[(b)] Assume that $p\geq 4w+15$. Then $N\leq w+1$.
\item[(c)] Assume that $p\geq 19$, the degree of $D$ is not a power of $p$,
the extension $F/\dbQ_p$ is Galois, and $F$ contains $p^2\rm{th}$ primitive
root of unity. Then $N\leq w+1$.
\item[(d)] Assume that $p\geq 19$, the degree of $D$ is not a power of $p$, 
and $F$ is a cyclotomic field. Then $N= w+1$.
\end{itemize} 
\label{Thmain}
\end{Theorem}
\vskip -.1cm
Note that Theorem~\ref{Thmain}(d) immediately follows from Theorem~\ref{Thmain}(b)(c)
and \cite[Theorem 8.1]{PR2}. Indeed, if $F=\dbQ_p(\sqrt[n]{1})$ and
$p^k$ is the highest power of $p$ dividing $n$, then $e=p^{k-1}(p-1)$ and $w=k-1$. 
Thus $|H^2(SL_1(D),\dbR/\dbZ)|\leq p^{k}$ by Theorem~\ref{Thmain}(c) if $k\geq 2$
and by Theorem~\ref{Thmain}(b) if $k=1$, while \cite[Theorem 8.1]{PR2} yields 
$|H^2(SL_1(D),\dbR/\dbZ)|\geq p^k$.
\vskip .3cm

We now give a brief sketch of the proof of
Theorem~\ref{Thmain}. Let $G=SL_1(D)$. In \cite{PR2} it is shown
that $H^2(G,(\dbQ/\dbZ)_p)$ is isomorphic to $H^2(G,\dbZ/p^k\dbZ)$
for sufficiently large $k$.
Both $G$ and $\dbZ/p^k\dbZ$ are $p$-adic analytic, so
it is natural to ask if the order of $H^2(G,\dbZ/p^k\dbZ)$
can be computed using Lie algebras.
In the theory of $p$-adic analytic groups
there is a well-known exp-log correspondence between (finitely generated)
powerful torsion-free pro-$p$ groups and powerful torsion-free
$\dbZ_p$-Lie algebras. This is not enough for our purposes;
however, what we can use is a work of Weigel~\cite{We}, who extended
the above correspondence to the classes of powerful $p$-central
pro-$p$ groups and Lie algebras (see \S~3 for definitions).

Now consider the congruence subgroup $H=SL_1^{de+1}(D)$ where $d$ is the degree of $D$.
It is easy to see that $H$ is powerful and torsion-free.
Let $\res_{G|H}$ be the restriction map from
$H^2(G,\dbZ/p^k\dbZ)$ to $H^2(H,\dbZ/p^k\dbZ)$.
At the end of \S~5 we will show
that $\res_{G|H}$ has small kernel
(see Proposition~\ref{infkernel}). Next we will prove that any cohomology class lying in the image of
$\res_{G|H}$ is represented by a (central) extension
$1\to \dbZ/p^k\dbZ\to \Hhat\to H\to 1$ where $\Hhat$ is
powerful and $p$-central (see Lemma~\ref{LemmaA}). 
Applying Weigel's $\log$ functor we obtain the corresponding extension
of powerful $p$-central Lie algebras. This extension, in turn,
represents some cohomology class in $H^2(\grh,\dbZ/p^k\dbZ)$,
where $\grh$ is the $\dbZ_p$-Lie algebra of $H$.
Moreover, the obtained cohomology class is invariant under
the natural action of $G$ on $H^2(\grh,\dbZ/p^k\dbZ)$.
These results lead to an upper bound for the order of $H^2(G,\dbZ/p^k\dbZ)$ in terms
of the exponent of the $G$-invariant part of $H^2(\grh,\dbZ/p^k\dbZ)$.
Finally, in \S~7 we obtain an explicit description of $G$-invariant classes in
$H^2(\grh,\dbZ/p^k\dbZ)$, which yields the bound
given in Theorem~\ref{Thmain}(a).
\skv
The proof of Theorem~\ref{Thmain}(b) is based on similar ideas,
but is considerably more technical. Instead of Weigel's correspondence
we use Lazard's $\exp$-$\log$ correspondence between finite
groups and finite Lie rings of $p$-power order and nilpotency class
less than $p$. The reduction of Theorem~\ref{Thmain}(b) to computation
of cohomology of finite $p$-groups is based on the analysis
of the inflation map $H^2(G/G_m,\dbZ/p^k\dbZ)\to H^2(G,\dbZ/p^k\dbZ)$
for $m\in\dbN$, where $G_m=SL_1^m(D)$.
\skv
Finally, to prove Theorem~\ref{Thmain}(c) we use the following
simple fact pointed out to the author by Gopal Prasad:
If $F/F_0$ is an extension of $p$-adic fields and $D$ is a central
division algebra over $F$ whose degree is relatively prime to $[F:F_0]$,
then $D\cong D_0\otimes_{F_0} F$ for some division algebra $D_0$
over $F_0$. Furthermore, certain information about the restriction map 
$H^2(SL_1(D),\dbR/\dbZ)\to H^2(SL_1(D_0),\dbR/\dbZ)$ is provided by \cite{PR2}.
Using this idea, we reduce the proof of Theorem~\ref{Thmain}(c)
to the case of division algebras over $p$-adic fields of small
degree, where Theorem~\ref{Thmain}(b) becomes applicable.
\vskip .12cm

\noindent
\bf{Organization. }\rm In \S~2 we discuss the basic relation between
cohomology and central extensions for profinite groups and Lie algebras.
In \S~3 we describe $\exp$-$\log$
correspondence between certain classes of $p$-adic analytic pro-$p$ 
groups and $\dbZ_p$-Lie algebras. We then use this correspondence
to establish a relationship between the second cohomology of
pro-$p$ groups and Lie algebras belonging to those classes.
In \S~4 we review basic facts about division algebras
over local fields. In \S~5 we study group-theoretic
properties of central extensions of $SL_1(D)$ where $D$ is
a division algebra over a $p$-adic field. 
In \S~6 we deduce parts (a) and (b) of Theorem~\ref{Thmain}
from certain results on Lie algebra cohomology which, in turn,
are proved in \S~7. Finally, in \S~8 we prove Theorem~\ref{Thmain}(c).

\vskip .12cm
\noindent
\bf{Restrictions on $p$. }\rm While all parts of Theorem~\ref{Thmain} require at least that 
$p\geq 5$, most of our auxiliary results hold for $p=3$ and some even for $p=2$. In \S~5 and \S~7
we will impose a general restriction on $p$ and $d$ (where $d$ is the degree of the division algebra $D$
under consideration) at the beginning of the section; in other sections we will state the restrictions
separately for each result.
\vskip .12cm
\noindent
\bf{Basic notations. }\rm Throughout the paper $\dbZ$ will stand for integers,
$\dbN$ for positive integers, $\dbZ_p$ for $p$-adic integers and
$\Fp$ for a finite field of order $p$.
If $G$ is a topological group, $\gamma_n G$ will denote the closure of the 
$n^{\rm th}$ term of the lower central series of $G$,
and $G^n$ the closed subgroup of $G$ generated by $n^{\rm th}$ powers.
If $A$ and $B$ are subsets of $G$, let $[A,B]$ be the closed subgroup generated
by $\{[a,b]\colon a\in A, b\in B\}$, where $[a,b]=a^{-1}b^{-1}ab$.
If $x_1,\ldots,x_k$ are elements of some Lie ring $L$, the left-normed commutator $[x_1,\ldots,x_k]$
is defined inductively by $[x_1,\ldots,x_k]=[[x_1,\ldots,x_{k-1}],x_k]$ for $k\geq 3$.
\vskip .12cm

\noindent
\bf{Acknowledgements. }\rm
I am very grateful to Gopal Prasad for posing the problem, very 
interesting conversations, and suggesting an idea that resulted
in significant improvement of the results of this paper. I am
very thankful to Thomas Weigel for sending me his unpublished 
manuscript \cite{We}. I would also like to thank Andrei Jaikin and Thomas Weigel for
helpful discussions related to this paper and the anonymous referee for a number of
useful suggestions.

\section{Central extensions and second cohomology}

\subsection{Central extensions and cohomology for profinite groups}
In this subsection we briefly recall the connection between topological central extensions of a profinite group
and its second continuous cohomology.

Let $H$ be a profinite group and let $A$ be an abelian profinite
group, considered as a trivial $H$-module. Denote by $\Ext(H,A)$ the set
of equivalence classes of {\it topological central extensions} of $H$ by $A$,
that is, extensions of the form 
$$1\to A\str{\iota}{\lra}\Hhat\str{\phi}{\lra}H\to 1
\quad\mbox{ or, in abbreviated form, }\quad 
\ex{A}{\Hhat}{H}{\iota}{\phi},$$
where $\iota$ and $\phi$ are continuous and $\iota(A)$ is central in $\Hhat$ (note that
the group $\Hhat$ is automatically profinite). Recall that two extensions
$1\to A\str{\iota_1}{\lra}\Hhat_1\str{\phi_1}{\lra}H\to 1$ and $1\to A\str{\iota_2}{\lra}\Hhat_2\str{\phi_2}{\lra}H\to 1$
are equivalent if there exists a (continuous)\footnote{Since profinite groups are compact and Hausdorff, 
the inverse map $\alpha^{-1}$ is automatically continuous as well.} isomorphism $\alpha:\Hhat_1\to\Hhat_2$ which makes the following diagram commutative:

\begin{equation}
\begin{CD}
1@ >>> A@ >\iota_1>> \Hhat_1@ >\phi_1 >> H@ >>> 1\\
@ VVV @ VV{\rm id}V @ VV\alpha V @ VV{\rm id}V @ VVV\\
1@ >>> A@ >\iota_2>> \Hhat_2@ >\phi_2 >> H@ >>> 1\\
\end{CD}
\end{equation}

To avoid complicated notations we will not be explicitly distinguishing between extensions and their equivalence classes.
\skv

Next recall that $\Ext(H,A)$ is an abelian group with respect to the Baer sum defined as follows:
let $\calE_1=(\ex{A}{\Hhat_1}{H}{\iota_1}{\phi_1})$ and 
$\calE_2=(\ex{A}{\Hhat_2}{H}{\iota_2}{\phi_2})$
be two elements of $\Ext(H,A)$. Then
$$\calE_1+\calE_2= (\ex{A}{\Hhat/\Nhat}{H}{\iota}{\phi})$$
where $\Hhat=\{(x_1,x_2)\in\Hhat_1\times\Hhat_2 : \phi_1(x_1)=\phi_2(x_2)\}$,
$\Nhat=\{(\iota_1(a),\iota_2(a^{-1})) :  a\in A\}$, 
$\iota(a)=(\iota_1(a),1)\Nhat=(1,\iota_2(a))\Nhat$, and 
$\phi((x_1,x_2)\Nhat)=\phi_1(x_1)=\phi_2(x_2)$.

For any element $\calE=(\ex{A}{\Hhat}{H}{\iota}{\phi})$ of $\Ext(H,A)$ its inverse is
$-\calE=(\ex{A}{\Hhat}{H}{\iota'}{\phi})$ where $\iota'(a)=\iota(a)^{-1}$ for all $a\in A$.
\skv

Now let $H^2(H,A)$ be the second continuous cohomology of $H$, that is, the second cohomology based on the standard
complex where cochains are required to be continuous. Thus $H^2(H,A)=Z^2(H,A)/B^2(H,A)$ where
the group $Z^2(H,A)$ of $2$-cocycles consists of continuous functions $f:H\times H\to A$ such that
$$f(xy,z)+f(x,y)=f(y,z)+f(x,yz) \mbox { for all }x,y,z\in H$$  
and the group $B^2(H,A)$ of $2$-coboundaries consists of continuous functions $f:H\times H\to A$ of the form
$f(x,y)=u(x)+u(y)-u(xy)$ for some continuous function $u:H\to A$.

There exists a canonical isomorphism 
of abelian groups $H^2(H,A)\cong \Ext(H,A)$ defined as follows:

Given $C\in H^2(H,A)$, let $Z:H\times H\to A$ be a $2$-cocycle whose cohomology
class is equal to $C$. Let $\Hhat$ be the set of pairs $\{(h,a): h\in H, a\in A\}$
with multiplication given by 
\begin{equation}
\label{eq:cocycle}
(h_1,a_1)\cdot (h_2,a_2)=(h_1 h_2, a_1+a_2+Z(h_1,h_2))
\end{equation}
The central extension corresponding to $C$ is
$ \ex{A}{\Hhat}{H}{\iota}{\phi}$
where $\iota(a)=(1,a)$ and $\phi((h,a))=h$ for any $a\in A$ and $h\in H$.
We will denote (the equivalence class of) this extension by $\Ext(C)$.

Conversely, let 
$\calE= (1\to A\str{\iota}{\lra}\Hhat\str{\phi}{\lra}H\to 1)$
be an element of $\Ext(H,A)$. Let $\psi: H\to\Hhat$ be a continuous section of $\phi$, 
that is, a continuous map $H\to\Hhat$ such that $\phi\circ\psi=\id_H$
(such a section exists since $\Hhat$ and $H$ are profinite -- see, e.g. \cite[1.3.3]{Wil}), 
and define
$Z:H\times H\to A$ by $$Z(h_1,h_2)=\iota^{-1}(\psi(h_1 h_2)^{-1}\psi(h_1)\psi(h_2)).$$
Then $Z$ is a $2$-cocycle whose cohomology class $[Z]$ is independent of the choice of $\psi$, and $\Ext([Z])=\calE$.

\subsection{Central extensions and cohomology for Lie algebras}

On the Lie algebra side, we will discuss the connection between central extensions and cohomology in the abstract (discrete)
case since we are not aware of the corresponding reference in the topological case. 
The discrete case will be sufficient for the purposes of this paper since we will work with extensions inside the category of $\dbZ_p$-Lie algebras which are finitely generated as $\dbZ_p$-modules, and any homomorphism between such Lie algebras is automatically continuous with respect to pro-$p$ topologies (see the discussion at the beginning of \S~3).

Fix a commutative ring $R$ with $1$, let $\grh$ be an $R$-Lie algebra, and let $\gra$ be an $R$-module considered as an abelian $R$-Lie algebra and as a trivial $\grh$-module. As in the group case, denote by $\Ext(\grh,\gra)$ the set of equivalence classes of central $R$-Lie algebra extensions $0\to\gra\str{\iota}{\lra}\grhat\str{\phi}{\lra}\grh\to 0$. The addition is defined in complete analogy with the group case.

It is still true that $\Ext(\grh,\gra)$ is isomorphic to $H^2(\grh,\gra)$; however, we are not aware of a simple-to-describe
cochain complex that can be used to compute the cohomology groups. Here the cohomology groups are defined functorially as
$H^k(\grh,\gra)=\Ext^k_{\mathcal U\grh}(R,\gra)$ where $\mathcal U\grh$ is the universal enveloping algebra of $\grh$
and $R$ is considered as a trivial $\grh$-module. 

Consider now the important special case where $\grh$ is a free $R$-module. In this case, one can compute the cohomology groups $H^k(\grh,\gra)$ similarly
to the group case, using the Chevalley-Eilenberg complex. In particular, $H^2(\grh,\gra)$ is isomorphic to 
$Z^2(\grh,\gra)/B^2(\grh,\gra)$ where $Z^2(\grh,\gra)$ consists of $R$-bilinear functions $f:\grh\times \grh\to\gra$ such that
\begin{itemize}
\item[(i)] $f(x,x)=0$ for all $x\in \grh$
\item[(ii)] $f([x,y],z)+f([y,z],x)+f([z,x],y)=0$ for all $x,y,z\in\grh$.
\end{itemize}
and $B^2(\grh,\gra)\subseteq Z^2(\grh,\gra)$ is the set of all $R$-bilinear maps $f:\grh\times \grh\to\gra$
of the form $f(x,y)=u([x,y])$ for some $R$-linear function $u:\grh\to\gra$.

\vskip .2cm
The canonical isomorphism between $H^2(\grh,\gra)$ and $\Ext(\grh,\gra)$ in this special case is defined similarly to groups:

If $c\in H^2(\grh,\gra)$ and $z:\grh\times\grh\to\gra$ is a $2$-cocycle representing
$c$, we define $\Ext(c)\in \Ext(\grh,\gra)$ to be the extension 
$0\to\gra\str{\iota}{\lra}\grhat\str{\phi}{\lra}\grh\to 0$,
where $\grhat=\grh\oplus \gra$ as an $R$-module with Lie bracket
$[(h_1,a_1), (h_2,a_2)]=([h_1,h_2], z(h_1,h_2))$.

\vskip .1cm
Conversely, let $\calE= (0\to\gra\str{\iota}{\lra}\grhat\str{\phi}{\lra}\grh\to 0)$
be an element of $\Ext(\grh,\gra)$. Since $\grh$ is a free $R$-module, 
there exists an $R$-linear section $\psi:\grh\to \grhat$, and it is easy to check that
$\calE=\Ext([z])$ where $z\in Z^2(\grh,\gra)$ is defined by
\begin{equation}
\label{LA:coc}
z(h,k)=\iota^{-1}\left(\psi([h,k])-[\psi(h),\psi(k)]\right).
\end{equation}
\vskip .1cm

Going back to the general case (where $\grh$ is not necessarily free), one can still define the groups
$Z^2(\grh,\gra)$ and $B^2(\grh,\gra)$ and the map $\Ext:Z^2(\grh,\gra)/B^2(\grh,\gra)\to \Ext(\grh,\gra)$ as above;
however, this map need not be surjective. In \cite{HZ} it was shown that at least in the special case $R=\dbZ$
there are natural generalizations
$\widetilde{Z^2}(\grh,\gra)$ and $\widetilde{B^2}(\grh,\gra)$ of the spaces $Z^2(\grh,\gra)$ and $B^2(\grh,\gra)$
defined above with the property that $\widetilde{Z^2}(\grh,\gra)/\widetilde{B^2}(\grh,\gra)\cong \Ext(\grh,\gra)$.
However, we do not know if there is a natural cochain complex whose cohomology groups are isomorphic to $H^k(\grh,\gra)$
which yields the same definitions of 2-cocycles and 2-coboundaries.

\subsection{Restriction and inflation maps}

In this subsection we will define the restriction and inflation maps on the second cohomology groups for profinite groups and Lie algebras directly in terms of the associated central extensions. We will also describe the image of the inflation map in terms of the associated central extensions.

We start by recalling the usual definition of the restriction and inflation maps for profinite groups in terms of cocycles. So let
$H$ be a profinite group, $N$ a closed subgroup of $H$ and $A$ an abelian profinite group considered as a trivial $H$-module.
The restriction ${\rm res}:H^2(H,A)\to H^2(N,A)$ is induced from the restriction map (in the obvious sense) on cocycles ${\rm Res}: Z^2(H,A)\to Z^2(N,A)$.

If $N$ is normal in $H$, the inflation ${\rm inf}:H^2(H/N,A)\to H^2(H,A)$ is induced from the map ${\rm Inf}: Z^2(H/N,A)\to Z^2(H,A)$ given by
$({\rm Inf\,}f)(x,y)=f(xN,yN)$ for all $f\in Z^2(H/N,A)$ and $x,y\in H$.

Let us now rephrase these definitions in terms of central extensions. So let $c\in H^2(H,A)$, and let
$\calE=(\ex{A}{\Hhat}{H}{\iota}{\phi})\in \Ext(H,A)$ be the corresponding central extension. Then it is easy to see that
${\rm res}(\calE)\in \Ext(N,A)$, the restriction of $\calE$, is given by ${\rm res}(\calE)=(\ex{A}{\phi^{-1}(N)}{N}{\iota}{\phi})$.

Assume now that $N$ is normal in $H$, $c\in H^2(H/N,A)$ and $$\calE=(\ex{A}{E}{H}{\iota}{\phi})\in \Ext(H/N,A)$$
is the corresponding central extension. Let $\pi:H\to H/N$ be the natural projection, and let
$P$ be the pullback of the diagram $E\stackrel{\phi}{\longrightarrow}H/N\stackrel{\pi}{\longleftarrow} H$, that is,
$$P=\{(x,y)\in E\times H :\, \phi(x)=\pi(y)\}.$$
It takes a little work (see, e.g. \cite[Exercise~6.6.4, p.~285]{Wb}) to check that the inflation
${\rm inf}(\calE)\in \Ext(H,A)$ is given by ${\rm inf}(\calE)=(\ex{A}{P}{H}{\iota'}{\phi'})$ where
$\iota'(a)=(\iota(a),1)$ for all $a\in A$ and $\phi'(x,y)=y$ for all $(x,y)\in P$.
\skv

Both definitions of the restriction and inflation maps have direct analogues for $R$-Lie algebras. The definitions in terms of central extensions
(which are most suitable for our purposes) remain exactly the same. The definitions in terms of cocycles extend (with obvious changes) to the subsets of the cohomology groups corresponding to extensions which split on the level of $R$-modules; for the general case (which is still quite similar) we refer the reader to \cite{HZ}.
\skv

The next lemma provides a characterization for
the image of the inflation map between second cohomology groups
in terms of the associated central extensions.
Although this is a standard result, we are not aware
of a reference in the literature.
\begin{Lemma}
Let $H$ be a profinite group, $N$ a closed normal subgroup of $H$, and let
$A$ be an abelian profinite group (considered as a trivial $H$-module).
Fix $c\in H^2(H,A)$ with
$\Ext(c)=(1\to A\str{\iota}{\lra} \widehat{H}\str{\phi}{\lra} H\to 1)$.
Let $\inf: H^2(H/N,A)\to H^2(H,A)$ be the inflation map.
\begin{itemize}
\item[(a)] Suppose that $c=\inf(c')$ for some $c'\in H^2(H/N,A)$.
Then there exists a continuous section $\psi$ of $\phi$
(i.e. a map $\psi:H\to\widehat{H}$ with $\phi \psi=\id_{H}$)
such that $\psi(N)$ is a normal subgroup of $\widehat{H}$.

\item[(b)] Conversely, suppose that $\phi$ has a continuous section $\psi$
such that $\psi(N)$ is a normal subgroup of $\widehat{H}$. 
Let $c'\in H^2(H/N,A)$ be the cohomology class corresponding
to the extension $$\calE'=(1\to A\str{\iota'}{\lra} \widehat{H}/\psi(N)\str{\phi'}{\lra} H/N\to 1)$$
where $\iota'$ and $\phi'$ are induced by $\iota$ and $\phi$, respectively. 
Then $c=\inf(c')$. 
\label{cohbasic}
\end{itemize}
\end{Lemma}
\proofs
In both parts we let $Q=H/N$ and $\pi:H\to Q$ is the natural projection.
\skv

(a) Let $(\ex{A}{\widehat{Q}}{Q}{\iota'}{\phi'})$ be the extension corresponding to $c'$.
Then (up to equivalence) we have $\Ext(c)=(\ex{A}{\Hhat}{H}{\iota}{\phi})$
where $\Hhat=\{(x,y)\in \widehat{Q}\times H\,: \phi'(x)=\pi(y)\}$,
and $\iota$ and $\phi$ are defined by $\iota(m)=(\iota'(m),1)$ for any $m\in A$
and $\phi(x,y)=y$ for any $y\in H$.

Now choose a continuous section $\psi': Q\to\widehat{Q}$ of $\phi'$
such that $\psi'(1)=1$, and define $\psi:H\to \Hhat$ by
$\psi(y)=(\psi'(\pi(y)),y)$. Clearly, $\psi$ is a continuous section of $\phi$,
and $\psi(N)=\{(1,y):y\in N\}$ is a normal subgroup of $\Hhat$.
\vskip .15cm

(b) Let $c''=\inf(c')$. Then by definition $\Ext(c'')=(\ex{A}{P}{H}{}{})$
where $P=\{(x,y)\in \widehat{H}/\psi(N)\times H\,: \phi'(x)=\pi(y)\}$.
It is straightforward to check that the map $\alpha: \Hhat\to P$ given by
$\alpha(z)=(z\psi(N),\pi(z))$ is a continuous isomorphism 
(hence also bi-continuous since both groups are compact) 
which yields equivalence of the extensions $(\ex{A}{P}{H}{}{})$ and $(\ex{A}{\Hhat}{H}{}{})$.
Thus $c=c''$, as desired.
\QED

\subsection {Equivariant extensions and cohomology classes}\rm
Let $H$ and $A$ be profinite groups with $A$ abelian, and suppose that $G$ is another group
which acts on both $H$ and $A$ by continuous automorphisms. 
An extension $\ex{A}{\Hhat}{H}{\iota}{\phi}$ will be called \bf{$G$-equivariant }\rm
if there exists an action of $G$ on $\Hhat$ by continuous automorphisms which is compatible with
the $G$-action on $H$ and $A$, that is, 
\skv
\centerline{$\iota(a)^g=\iota(a^g)$
for any $a\in A$ and $g\in G$, and $\phi(x^g)=\phi(x)^g$ for any
$x\in \Hhat$ and $g\in G$.}
\skv
We will denote the subset of $G$-equivariant extensions by $\Ext(H,A)^G$.
\skv
An element $c\in H^2(H,A)$ will be called $G$-equivariant if
$\Ext(c)\in \Ext(H,A)^G$, and $H^2(H,A)^G$ will denote the set of 
$G$-equivariant elements. 
Note that the standard meaning of $H^2(H,A)^G$ is different from ours:
$H^2(H,A)^G$ usually denotes the set of cohomology classes which are 
invariant with respect to the canonical action of $G$ on $H^2(H,A)$.
It is not hard to show that $G$-equivariant cohomology
classes are invariant under the $G$-action of $H^2(H,A)$, but the converse
is not necessarily true (these facts will not be used in our paper).

Likewise if $\grh,\gra$ are $R$-Lie algebras, with $\gra$ abelian, and $G$ acts on
both $\grh$ and $\gra$ by $R$-Lie algebra automorphisms,
we define $\Ext(\grh,\gra)^G$ and $H^2(\grh,\gra)^G$ in the same way.
\skv

All parts of Lemma~\ref{equiv_subgp} below are straightforward and follow immediately from the definitions of addition and inversion in
$\Ext$ groups  and the definitions of restriction and inflation maps in terms of central extensions.

\begin{Lemma}
\label{equiv_subgp}
The following hold:
\begin{itemize}
\item[(a)] Let $H$ and $A$ be profinite groups, with $H$ abelian, and let $G$ be a group which acts
on both $H$ and $A$ by continuous automorphisms. Then
\begin{itemize}
\item[(i)] $\Ext(H,A)^G$ is a subgroup of $\Ext(H,A)$.
\item[(ii)] Let $N$ be a $G$-invariant closed subgroup of $H$. Then the restriction map $\Ext(H,A)\to \Ext(N,A)$
sends $\Ext(H,A)^G$ to $\Ext(N,A)^G$.
\item[(iii)] Let $N$ be a $G$-invariant closed normal subgroup of $H$. Then the inflation map $\Ext(H/N,A)\to \Ext(H,A)$
sends $\Ext(H/N,A)^G$ to $\Ext(H,A)^G$.
\end{itemize}
\item[(b)] Let $R$ be a commutative ring with $1$, let $\grh$ and $\gra$ be $R$-Lie algebras, with $\grh$ abelian, and let $G$ be a group which acts
on both $\grh$ and $\gra$ by $R$-Lie algebra automorphisms. Then
\begin{itemize}
\item[(i)] $\Ext(\grh,\gra)^G$ is a subgroup of $\Ext(\grh,\gra)$.
\item[(ii)] Let $\grn$ be a $G$-invariant Lie subalgebra of $\grh$. Then the restriction map $\Ext(\grh,\gra)\to \Ext(\grn,\gra)$
sends $\Ext(\grh,\gra)^G$ to $\Ext(\grn,\gra)^G$.
\item[(iii)] Let $\grn$ be a $G$-invariant ideal of $\grh$. Then the inflation map $\Ext(\grh/\grn,\gra)\to \Ext(\grh,\gra)$
sends $\Ext(\grh/\grn,\gra)^G$ to $\Ext(\grh,\gra)^G$.
\end{itemize}
\end{itemize}
\end{Lemma}

\section{Exp-log correspondence}

\subsection{Exp and log functors and their basic properties}
We will say that a $\dbZ_p$-Lie algebra
$L$ is of {\it finite rank} if $L$ is finitely generated as a $\dbZ_p$-module.
Any such $L$ is isomorphic as a $\dbZ_p$-module to $\dbZ_p^n\oplus P$ for some $n\in\dbZ_{\geq 0}$ and finite abelian $p$-group $P$, and thus $(L,+)$ has a structure of a finitely generated pro-$p$ group. This structure is independent of the choice of an isomorphism
$L\cong \dbZ_p^n\oplus P$, e.g., since any homomorphism between finitely generated pro-$p$ groups is continuous\footnote {This fact is not very hard to prove -- see \cite[Corollary~1.19]{DDMS}. It is even true that any homomorphism between finitely generated profinite groups is continuous, but this is a very deep theorem of Nikolov and Segal~\cite{NS1,NS2}}. Moreover,
in this way $L$ becomes a (Hausdorff) topological $\dbZ_p$-Lie algebra, that is, the addition and Lie bracket are continuous as maps from $L\times L$ to $L$ and scalar multiplication is continuous as a map from $\dbZ_p\times L$ to $L$.
 \skv
 
Let $\grL_{\dbZ_p}$ (resp. $\grG_{\dbZ_p}$) be 
the category whose objects are $\dbZ_p$-Lie algebras of finite rank 
(resp. $p$-adic analytic pro-$p$ groups) and whose morphisms are 
$\dbZ_p$-Lie algebra (resp. group) homomorphisms. Since $p$-adic analytic pro-$p$ groups are finitely generated,
by the explanation in the previous paragraph any morphism in $\grL_{\dbZ_p}$ or $\grG_{\dbZ_p}$ is automatically continuous.

If $\grL$ is a subcategory of $\grL_{\dbZ_p}$ and $\grG$ is a subcategory
of $\grG_{\dbZ_p}$, by an \bf{$\exp$-$\log$ correspondence }\rm between $\grL$
and $\grG$ we mean a pair of functors $\exp:\grL\to\grG$ and $\log:\grG\to\grL$
such that the compositions $\exp\circ\log$ and $\log\circ\exp$ are naturally
equivalent to the identity functors on $\grG$ and $\grL$, respectively.
We will describe such correspondence in the following cases
(all relevant definitions are given later in this section). In all cases below
we assume that $\grL$ and $\grG$ are full subcategories of $\grL_{\dbZ_p}$ and $\grG_{\dbZ_p}$, respectively.

1. $(\grL,\grG)=(\grL_{<p},\grG_{<p})$ where
$\grL_{<p}$ (resp. $\grG_{<p}$) is the category of finite Lie rings (resp. finite groups) of 
$p$-power order and nilpotency class $<p$.

2. $(\grL,\grG)=(\grL_{ptf},\grG_{ptf})$ where
$\grL_{ptf}$ (resp. $\grG_{ptf}$) is the category of powerful torsion-free $\dbZ_p$-Lie algebras
of finite rank (resp. finitely generated powerful torsion-free pro-$p$ groups).\footnote{Finitely generated powerful torsion-free pro-$p$ groups are precisely the uniformly powerful
(or {\it uniform}) pro-$p$ groups \cite[Theorem~4.5]{DDMS}.}

3. $(\grL,\grG)=(\grL_{ppc},\grG_{ppc})$ where
$\grL_{ppc}$ (resp. $\grG_{ppc}$) is the category of powerful $p$-central $\dbZ_p$-Lie algebras
of finite rank
(resp. finitely generated powerful $p$-central pro-$p$ groups) and $p\geq 5$.

\vskip .12cm
Cases 1 and 2 of $\exp$-$\log$ correspondence are due to Lazard.
The correspondence $\grL_{<p}\cong \grG_{<p}$ is a special case of
\cite[Theorem~4.6]{La1}. Equivalence between $\grL_{ptf}$ and $\grG_{ptf}$ 
is essentially established in Lazard's famous 1965 paper on $p$-adic analytic
groups~\cite{La}, although the notion of a powerful group was introduced more 
than 20 years later by Lubotzky and Mann~\cite{LM}. For a detailed
account of the theory of poweful groups the reader is referred to 
the excellent book on analytic pro-$p$ groups~\cite{DDMS};
we shall just state the main definitions and results.

\begin{Definition}\rm
A pro-$p$ group $G$ (resp. a $\dbZ_p$-Lie algebra $L$)
is called \it{powerful, }\rm if $[G,G]\subseteq G^q$ (resp. $[L,L]\subseteq qL$)
where $q=p$ if $p>2$ and $q=4$ if $p=2$.
\end{Definition}

The following well-known criterion of analyticity of pro-$p$ groups
was first stated in \cite{LM} and is easily deduced from results in \cite{La}.

\begin{Theorem}
\label{powpa}
A finitely generated pro-$p$ group is $p$-adic analytic if and only if
it contains a finite index powerful subgroup. Moreover,
every finitely generated powerful pro-$p$ group contains a finite index subgroup which
is powerful and torsion-free.
\end{Theorem}

The book \cite{DDMS} contains a full proof of Theorem~\ref{powpa} ``from scratch''
as well as an explicit proof of equivalence $\grG_{ptf}\cong\grL_{ptf}$ 
between the categories of finitely generated powerful torsion-free pro-$p$ groups and $\dbZ_p$-Lie algebras
of finite rank.
Lazard's counterpart of this result \cite[Chapter IV, Theorem 3.2.6]{La}
is a correspondence between the categories of ``$p$-saturable'' pro-$p$ groups and Lie algebras.
Any finitely generated torsion-free powerful pro-$p$ group is $p$-saturable; conversely, 
a $p$-saturable pro-$p$ group is $p$-adic analytic and torsion-free, but not
necessarily powerful. Thus, Lazard's $\exp$-$\log$ correspondence is more general than the one 
between $\grG_{ptf}$ and $\grL_{ptf}$; however, powerful torsion-free pro-$p$ groups are usually easier to work
with than $p$-saturable ones (for more on this see \cite{K}).

The last case of $\exp$-$\log$ correspondence used in this paper
is Weigel's generalization of the correspondence $\grL_{ptf}\cong\grG_{ptf}$
to certain classes of pro-$p$ groups and $\dbZ_p$-Lie algebras that are 
powerful but not necessarily torsion-free.

\begin{Definition}\rm
Assume that $p>2$. A pro-$p$ group $G$ (resp. a $\dbZ_p$-Lie algebra $L$)
is called \it{$p$-central, }\rm if any $g\in G$ such that $g^p=1$
(resp. $u\in L$ such that $pu=0$) lies in the center of $G$ (resp. $L$).
\end{Definition}
In \cite{We}, 
\footnote{In fact, Weigel introduced a general technique for establishing
$\exp$-$\log$ correspondence between categories of $p$-adic analytic groups
and $\dbZ_p$-Lie algebras satisfying certain conditions. This technique
is applicable to all cases of $\exp$-$\log$ correspondence discussed in this paper. }
Weigel constructed $\exp$-$\log$ correspondence $\grL_{ppc}\cong \grG_{ppc}$
between the categories of powerful $p$-central pro-$p$ groups and $\dbZ_p$-Lie algebras
for $p\geq 5$. 
Note that a torsion-free pro-$p$ group is always $p$-central, and more
generally, a central extension of a torsion-free pro-$p$ group is always $p$-central.
Thus, Weigel's correspondence is well suited for computing the second cohomology of 
powerful torsion-free pro-$p$ groups.
\skv

\noindent
\bf{Construction of the $\exp$ functor. }\rm
We shall now explain how to construct the pro-$p$ group $\exp(L)$
corresponding to a $\dbZ_p$-Lie algebra $L$ where $L\in Ob(\grL_{<p})$
or $L\in Ob(\grL_{ppc})$.  While there exist distinct ways to define $\exp(L)$
formally, they are all based on the Baker-Cambell-Hausdorff (BCH) formula.
 
Let $\calA=\dbQ\la\la x_1,x_2\ra\ra$ 
be the algebra of power series over $\dbQ$ in two non-commuting variables $x_1$ and $x_2$.
The power series $\Phi=\Log(\Exp(x_1)\cdot \Exp (x_2))$ is called
the \it{Baker-Campbell-Hausdorff series }\rm
(here $\Exp(x)=1+x+\frac{x^2}{2}+\ldots$ and $\Log(1+x)=x-\frac{x^2}{2}+\frac{x^3}{3}-\ldots$).

\begin{Theorem}
\label{CHF}
The Baker-Campbell-Hausdorff (BCH) series $\Phi$ lies
in the closed $\dbQ$-Lie subalgebra of $\calA$ generated by $x_1$ and $x_2$.
In other words,
$\Phi=\sum_{c\in S}\lam_c c$,
where $S$ is the set of all left-normed commutators
in $x_1, x_2$ and each $\lam_c\in\dbQ$.
Moreover, if $wt(c)$ denotes the weight of a commutator $c$, then
$$p^{\,\lfloor \frac{wt(c)-1}{p-1}\rfloor}\lam_c\in \dbZ_p \mbox{ for any }c\in S
\quad\mbox{ where } \lfloor t\rfloor \mbox{ is the largest integer }\leq t.$$
\end{Theorem}
\begin{Remark}
There is an explicit expression for $\Phi$
(as a linear combination of commutators),
called the Baker-Campbell-Hausdorff formula.
The last assertion of Theorem~\ref{CHF} due to
Lazard~\cite{La} is a consequence of that formula.
\end{Remark}
\skv

Now let $\grL=\grL_{<p}$ or $\grL_{ppc}$, and let $L$ be an object of $\grL$. We define the pro-$p$ group 
$\exp(L)$ as follows. As a set $\exp(L)=L$, and the group operation $\cdot$ is defined by
\vskip -.3cm
\begin{equation*}
u_1 \cdot u_2=\Phi(u_1,u_2).
\label{eqBCH}
\end{equation*}
where $\Phi(u_1,u_2)\in L$ is defined below.
Informally, one should think of $\Phi(u_1,u_2)$ as the result of ``evaluating''
the BCH series at $x_1=u_1$ and $x_2=u_2$. The formal definition
of $\Phi(u_1,u_2)$ will be different in the cases $\grL=\grL_{<p}$ and $\grL=\grL_{ppc}$.
\skv

\bf{Case 1: }\rm $\grL=\grL_{<p}$. 
Given a left-normed commutator $c$ in $x_1, x_2$, we 
define $c(u_1,u_2)$ by substituting
$u_1$ for $x_1$ and $u_2$ for $x_2$ in $c$.
Thus, if $c=[x_{i_1},x_{i_2},\ldots, x_{i_k}]$,
then $c(u_1,u_2)=[u_{i_1},u_{i_2},\ldots, u_{i_k}]$.
Since the nilpotency class of $L$ is less than $p$,
we have $c(u_1,u_2)=0$ whenever $wt(c)\geq p$. 
Thus, we set $$\Phi(u_1,u_2)=\sum\limits_{c\in S, wt(c)<p}\lam_c c(u_1,u_2)$$
(using the notations of Theorem~\ref{CHF}).
In other words, $\Phi(u_1,u_2)$ is obtained by plugging in $u_1$ and $u_2$
into the BCH series truncated after degree $p-1$. Since $\lam_c\in\dbZ_p$
whenever $wt(c)<p$, the obtained expression is well-defined.
\skv

\bf{Case 2: }\rm $\grL=\grL_{ppc}$ and $p\geq 5$. 
Once again, let $c=[x_{i_1},x_{i_2},\ldots, x_{i_k}]$ be a left-normed 
commutator of weight $k\geq 2$. In this case for any $u_1,u_2\in L$
we can define the value $\frac{1}{p^{k-2}}c(u_1,u_2)$ as follows.
If $k=2$, this is done in the obvious way. 

Assume now that $k=3$, so that $c=[x_{i_1},x_{i_2},x_{i_3}]$.
Since $L$ is powerful, there exists
$v_1\in L$ such that $[u_{i_1},u_{i_2}]=pv_1$,
and we define $\frac{1}{p}c(u_1,u_2)=[v_1,u_{i_3}]$.
The last expression is independent of the choice of $v_1$
because $L$ is $p$-central. Indeed, if $[u_{i_1},u_{i_2}]=pv'_1$
for some $v'_1\neq v_1$, then $p(v'_1-v_1)=0$, whence
$v'_1-v_1$ lies in the center of $L$.

For $k\geq 4$ we proceed inductively -- write $c=[d,x_{i_k}]$
where $d$ is a left-normed commutator of weight $k-1$,
and define $\frac{1}{p^{k-2}}c(u_1,u_2)=[\frac{1}{p}\cdot \frac{1}{p^{k-3}}d(u_1,u_2), u_{i_k}]$
where $\frac{1}{p}w$ is any element of $L$ such that $p\cdot \frac{1}{p}w=w$
(again the choice of such an element does not affect the value $\frac{1}{p^{k-2}}c(u_1,u_2)$
since $L$ is $p$-central).

Now recall from  Theorem~\ref{CHF} that
 $p^{\,\lfloor \frac{k-1}{p-1}\rfloor}\lam_c\in\dbZ_p$
where $\lam_c$ is the coefficient of $c$ in the (chosen expansion of) BCH series $\Phi$.
Since $\lfloor \frac{k-1}{p-1}\rfloor\leq k-2$ for $p\geq 3$, we can define
$\lam_c c(u_1,u_2)\in L$ by setting
$\lam_c c(u_1,u_2)=(p^{k-2}\lam_c)\cdot \frac{1}{p^{k-2}}c(u_1,u_2)$.
Moreover, the series $\sum\limits_{c}\lam_c c(u_1,u_2)$
converges in $L$, and we let $\Phi(u_1,u_2)$ be its sum.
\skv

It is now clear how to define the functor $\exp:\grL\to\grG$
where $(\grL,\grG)= (\grL_{<p},\grG_{<p})$ or $(\grL_{ppc},\grG_{ppc})$:
\begin{itemize}
\item if $L$ is an object of $\grL$, the corresponding
object of $\grG$ is the group $\exp(L)$ as defined above
\item if $L_1, L_2$ are objects of $\grL$ and $f:L_1\to L_2$
is a Lie algebra homomorphism, the corresponding group 
homomorphism $\exp(f): \exp(L_1)\to \exp(L_2)$
coincides with $f$ as a set map.
\end{itemize}

Of course, there is a number of things to check here for both $\grL=\grL_{<p}$
and $\grL=\grL_{ppc}$, but the arguments in the latter case are much more involved.
Proving that $G=\exp(L)$ is a group with respect to the operation defined by $\Phi$ turns out to be technically most difficult, and requires the assumption $p\geq 5$ (so far we only used that
$p\geq 3$). Once this is done, it is fairly easy to check that $G$ is an object
of $\grG$ and for any Lie algebra homomorphism $f$ between two objects of $\grL$, the corresponding map $\exp(f)$ is a group homomorphism.

The functor $\log:\grG\to\grL$ is defined somewhat similarly, but both the construction and its justification are more demanding. Let $\calL$ be the closed $\dbQ$-Lie subalgebra of 
 $\dbQ\la\la x_1,x_2\ra\ra$ generated by $x_1$ and $x_2$. To define $\log$,
one first expresses the addition and the Lie bracket on $\calL$ in terms of the  multiplication given by $\Phi$. The obtained expressions, denoted by $\Psi_A$ and $\Psi_B$, are no longer power series; instead they are infinite products of ``root-commutators'' which are defined as left-normed commutators, except that one is allowed
to extract $p^{\rm th}$-root a limited number of times (e.g. the root-commutator
$[[[x_1,x_2]^{1/p},x_3]^{1/p},x_4]$ would correspond to the element $\frac{1}{p^2}[x_1,x_2,x_3,x_4]$
in the Lie algebra case). Then, given a pro-$p$ group $G\in Ob(\grG)$,
one defines the addition and the Lie bracket on $G$ by evaluating the
infinite products $\Psi_A$ and $\Psi_B$ on $G$
(similarly to how we evaluated $\Phi$ on Lie algebras in $Ob(\grL)$).

\skv
We now formally state the main result on $\exp$-$\log$ correspondence (Theorem~\ref{Weigel0} below) 
established
by Lazard~\cite{La1} for the pair $(\grL_{<p},\grG_{<p})$ and by Weigel~\cite{We} 
for the pair $(\grL_{ppc},\grG_{ppc})$. Since Weigel's manuscript~\cite{We} is not easily available,
a proof of Theorem~\ref{Weigel0} in the case of Weigel's correspondence will be given in the appendix.

\begin{Theorem} 
\label{Weigel0}
Let $(\grG,\grL)=(\grG_{ppc},\grL_{ppc})$ with $p\geq 5$ or $(\grG,\grL)= (\grG_{<p},\grL_{<p})$.
There exist mutually inverse functors $\exp:\grL\to\grG$ and $\log:\grG\to\grL$
(that is, the compositions $\exp\circ\, \log$ and $\log\circ\, \exp$ are both the identity functors)
satisfying the following properties:
\begin{itemize}
\item[(1) ] For any object $L$ of $\grL$ we have $\exp(L)=L$ as a set and likewise
for any object $G$ of $\grG$ we have $\log(G)=G$ as a set.
\item[(2) ] If $\phi:L\to M$ is any morphism in $\grL$
and $\exp(\phi):\exp(L)\to\exp(M)$ is the corresponding
morphism in $\grL$, then $\exp(\phi)=\phi$ as set maps.
The analogous property holds for the $\log$ functor.
\item[(3)] For any object $L$ of $\grL$ the group operation on $G=\exp(L)$
is defined as described above (in terms of the BCH series $\Phi$).
\end{itemize}
\end{Theorem}

Below we collect several basic properties of the functors $\log$ and $\exp$ which will be used in this paper and 
easily follow from Theorem~\ref{Weigel0}.

\begin{Proposition} 
\label{Weigel}
Let $(\grG,\grL)=(\grG_{ppc},\grL_{ppc})$ with $p\geq 5$ or $(\grG,\grL)= (\grG_{<p},\grL_{<p})$,
and let $\exp:\grL\to\grG$ and $\log:\grG\to\grL$ be as in Theorem~\ref{Weigel0}.
Then $\log$ satisfies properties (a)-(d) below and $\exp$ satisfies their natural
analogues:

\begin{itemize}

\item[(a)] If $K,G,H\in Ob(\grG)$ and $1\to K\str{\iota}{\lra}G\str{\phi}{\lra}H\to 1$ is an exact
sequence, then the sequence 
$0\to \log(K)\str{\log(\iota)}{\lra}\log(G)\str{\log(\phi)}{\lra}\log(H)\to 0$
is also exact.

\item[(b)] If $G\in Ob(\grG)$ is abelian, then $\log(G)\in\grL$ is also abelian.

\item[(c)] For any $G,H\in Ob(\grG)$ we have $\log(G\times H)=\log(G)\times \log(H)$.

\item[(d)] If $H,G\in Ob(\grG)$ and $H$ is a central subgroup of $G$, then 
$\log(H)$ is a central subalgebra of $\log(G)$.
\end{itemize}
\end{Proposition}
\begin{proof} The analogues of (a)-(d) for the $\exp$ functor follow immediately
from Theorem~\ref{Weigel0}. For (b) and (d) we use the obvious fact that
$\Phi(u,v)=u+v$ whenever $u$ and $v$ commute in some $L\in Ob(\grL)$.
\skv

We could have argued similarly for the $\log$ functor had we described the latter more
explicitly, but it is actually easy to deduce
(a)-(d) for $\log$ directly from (a)-(c) for $\exp$ and Theorem~\ref{Weigel0}. 
Property (a) for $\log$ is automatic by Theorem~\ref{Weigel0}.
Below we will give an argument for (b) and (d), with (c) being similar to (d).

\skv

(b) Suppose that $G\in Ob(\grG)$ is abelian. Then $G$ has a unique structure of an abelian 
$\dbZ_p$-Lie algebra where the addition is given by the group operation in $G$. If $L$ is the obtained
$\dbZ_p$-Lie algebra, it is clear that $\exp(L)=G$, and hence $\log(G)=\log(\exp(L))=L$ is abelian
by (b) for $\exp$. 
\skv

(d) Take any $z\in H$ and $g\in G$. By assumption, $z$ and $g$ commute in $G$, and we need to show that they
also commute  in $\log(G)$. The (closed) subgroup of $G$ generated by $z$ and $g$ is abelian
and hence also lies in $\grG$. Call this subgroup $K$, and let $\iota:K\to G$ be the inclusion map. By 
Theorem~\ref{Weigel0}, $\iota$ considered as a map from $\log(K)$ to $\log(G)$ is a Lie algebra homomorphism, and by 
(b) $\log(K)$ is abelian. Since $\log(K)$ contains both $z$ and $g$, it follows that $z$ and $g$ commute in $\log(K)$ (and hence also in $\log(G)$), as desired. 
\end{proof}

\subsection{Central extensions and $\exp$-$\log$ correspondence}

For the rest of the paper $\grG$ (resp. $\grL$) will always denote either
$\grG_{ppc}$ (resp. $\grL_{ppc}$) with $p\geq 5$ or $\grG_{<p}$ (resp. $\grL_{<p}$).
Whenever we specify a choice for $\grG$ in a particular result, it will be assumed that
$\grL$ denotes the corresponding Lie algebra category.

\skv

Given $H\in Ob(\grG)$ and $A\in Ob(\grG)$, with $A$ abelian,
we define $\Ext_{\grG}(H,A)$ to be the subset of $\Ext(H,A)$ consisting of
extensions $\ex{A}{\Hhat}{H}{}{}$ such that $\Hhat\in Ob(\grG)$ as well.
Likewise, given $\grh,\gra\in Ob(\grL)$, with $\gra$ abelian, we define
$\Ext_{\grL}(\grh,\gra)$ to be the subset of $\Ext(\grh,\gra)$ consisting of
extensions $\ex{\gra}{\grhat}{\grh}{}{}$ such that $\grhat\in Ob(\grL)$ as well.

\begin{Proposition} Let $\grL$ and $\grG$ be as above. Let $H,A\in Ob(\grG)$
where $A$ is abelian, and let $\grh=\log(H)$, $\gra=\log(A)$.
\begin{itemize}
\item[(a)] There exists a natural bijection $\Log:\Ext_{\grG}(H,A)\to \Ext_{\grL}(\grh,\gra)$.

\item[(b)] Suppose that $\grG=\grG_{<p}$, or $\grG=\grG_{ppc}$ and $H$ is torsion-free. Then
\begin{itemize}
\item[(i)] $\Ext_{\grG}(H,A)$ is a subgroup of $\Ext(H,A)$
\item[(ii)] $\Ext_{\grL}(\grh,\gra)$ is a subgroup of $\Ext(\grh,\gra)$
\item[(iii)] The map $\Log:\Ext_{\grG}(H,A)\to \Ext_{\grL}(\grh,\gra)$ is an isomorphism
of abelian groups.
 \end{itemize}
 \end{itemize}
\label{extcor}
\end{Proposition}
\proofs

(a) is an immediate consequence of Proposition~\ref{Weigel}. Since the sum of extensions
is defined in the same way for groups and Lie algebras, (b)(iii) also becomes automatic once (b)(i) and (b)(ii) are established.

The proofs of (b)(i) and (b)(ii) are similar; in fact, the proof of (b)(ii) is slightly easier, so we will only prove (b)(i).

Let $\calE_1=(\ex{A}{\Hhat_1}{H}{\iota_1}{\phi_1})$ and 
$\calE_2=(\ex{A}{\Hhat_2}{H}{\iota_2}{\phi_2})$
be two elements of $\Ext_{\grG}(H,A)$. Recall that 
$\calE_1+\calE_2= (\ex{A}{\Hhat/\Nhat}{H}{\iota}{\phi})$
where $\Hhat=\{(x_1,x_2)\in\Hhat_1\times\Hhat_2 : \phi_1(x_1)=\phi_2(x_2)\}$,
$\Nhat=\{(\iota_1(a),\iota_2(a^{-1})) :  a\in A\}$, 
$\iota(a)=(\iota_1(a),1)\Nhat=(1,\iota_2(a))\Nhat$, and 
$\phi((x_1,x_2)\Nhat)=\phi_1(x_1)=\phi_2(x_2)$.
\skv
To prove (i) we need to show that $\Hhat/\Nhat\in Ob(\grG)$. If $\grG=\grG_{<p}$,
this is obvious since $\Hhat_1,\Hhat_2\in Ob(\grG_{<p})$ and $\grG_{<p}$ is closed
under subgroups, quotients and direct products. 

Now assume that $\grG=\grG_{ppc}$ and $H$ is torsion-free.
We need to show that $\Hhat/\Nhat$ is powerful and $p$-central.
The $p$-centrality condition clearly holds
since $\Hhat/\Nhat$ is a central extension of $H$. To prove
that $\Hhat/\Nhat$ is powerful it is sufficient to prove
that $\Hhat$ is powerful. We shall use the following
well-known criterion \cite[Lemma~3.4]{DDMS}.

\begin{Lemma} A finitely generated pro-$p$ group $G$ is powerful if and only if
for any $x,y\in G$ there exists $z\in G$ such that $[x,y]=z^p$.
\label{powchar}
\end{Lemma}

Now take any $x,y\in \Hhat$. Thus $x=(x_1,x_2)$ and $y=(y_1,y_2)$,
where $x_i,y_i\in \Hhat_i$ for $i=1,2$, $\phi_1(x_1)=\phi_2(x_2)$ and $\phi_1(y_1)=\phi_2(y_2)$.
Since $\Hhat_1$ and $\Hhat_2$ are powerful, there exist $z_i\in\Hhat_i,$
$i=1,2$ such that $[x_i,y_i]=z_i^p$. We have
\begin{equation}
\label{eqpairs}
[x,y]=[(x_1,x_2),(y_1,y_2)]=([x_1,y_1],[x_2,y_2])=(z_1^p,z_2^p)=(z_1,z_2)^p
\end{equation}
Since $(z_1^p,z_2^p)=[x,y]\in\Hhat$, we must have $\phi_1(z_1^p)=\phi_2(z_2^p)\in H$.
Since $H$ is powerful torsion-free, the equality $\phi_1(z_1)^p=\phi_2(z_2)^p$
implies that $\phi_1(z_1)=\phi_2(z_2)$ by \cite[Lemma~4.10]{DDMS}, whence $(z_1,z_2)\in \Hhat$.
Thus $\Hhat$ is powerful by \eqref{eqpairs} and Lemma~\ref{powchar}. 
The proof of (i) is complete.
\QED

The next result relates $G$-equivariant central extensions in the categories 
$\grG$ and $\grL$.

\begin{Proposition} 
\label{exthom1}
Let $\grL$ and $\grG$ be as above. Let $H,A\in Ob(\grG)$
where $A$ is abelian, and let $\grh=\log(H)$, $\gra=\log(A)$.
Let $G$ be a group which acts by (continuous) automorphisms on both $H$ and $A$, 
and define $\Ext_{\grG}(H,A)^G=\Ext_{\grG}(H,A)\cap \Ext(H,A)^G$.
The following hold:
\begin{itemize}
\item[(a)] $G$ acts on $\grh$ by $\dbZ_p$-Lie algebra automorphisms (where
the action is the same as the action on $H$).
Moreover, the map
$\Log$ defined in Proposition~\ref{extcor}(a) maps
$\Ext_{\grG}(H,A)^G$ onto $\Ext_{\grL}(\grh,\gra)^G$ (where by definition
$\Ext_{\grL}(\grh,\gra)^G=\Ext_{\grL}(\grh,\gra)\cap \Ext(\grh,\gra)^G$).
 
\item[(b)] Assume $\grG=\grG_{<p}$ or $\grG=\grG_{ppc}$ and $H$ is torsion-free.
Then $\Ext_{\grG}(H,A)^G$ and $\Ext_{\grL}(\grh,\gra)^G$
are subgroups of $\Ext(H,A)$ and $\Ext(\grh,\gra)$, respectively,
and $\Ext_{\grG}(H,A)^G\cong \Ext_{\grL}(\grh,\gra)^G$.
\end{itemize}
\end{Proposition}
\proofs
(a) is an immediate consequence of Theorem~\ref{Weigel0} and Proposition~\ref{Weigel}(a)(d).
\skv
(b) $\Ext_{\grG}(H,A)^G=\Ext(H,A)^G\cap \Ext_{\grG}(H,A)$ and $\Ext_{\grL}(\grh,\gra)^G=\Ext(\grh,\gra)^G\cap \Ext_{\grL}(\grh,\gra)$
are subgroups of $\Ext(H,A)$ and $\Ext(\grh,\gra)$, respectively, by 
Lemma~\ref{equiv_subgp} and Proposition~\ref{extcor}(b)(i)(ii).
Proposition~\ref{extcor}(b)(iii) and Proposition~\ref{exthom1}(a)
imply that these subgroups are isomorphic to each other.
\QED

\section{Division algebras over local fields and their norm one groups}
\label{sec:norm1}

In this paper by a local field we will always mean a nonarchimedean local field. While the results of the paper only deal with division algebras over local fields of characteristic zero, that is, finite extensions of $\dbQ_p$ for some $p$ (we will call them
{\it $p$-adic fields}), in this section we do not impose the characteristic restriction, as the basic structure of division algebras over local fields is very similar in zero and positive characteristics. We will only need very basic results about local fields which can be found, e.g., in Serre's book~\cite{Se}. Our basic reference for the structure of division algebras over local fields and their norm one groups is Riehm's paper~\cite{Ri}.
\skv

{\bf Notation.} Throughout the paper, if $D$ is a finite-dimensional division algebra over a local field (in particular, $D$ itself could be a local field), we will denote the ring of integers of $D$ by $O_D$ and the (unique) maximal ideal of $O_D$ by $\grm_D$. 
A {\it uniformizer} of $D$ is any generator of $\grm_D$.
If $K$ is a local field, its residue field $O_K/\grm_K$ will be denoted by $\Kbar$. 

\subsection{ Division algebras over local fields}
For the rest of this subsection we fix a prime $p$ and a local field $F$ of residue characteristic $p$.
Let $D$ be a finite-dimensional central division algebra over $F$,
and let $d$ be the degree of $D$. It is well known that there exists a maximal subfield
$W$ of $D$ such that the extension $W/F$ is unramified (note that $[W:F]=d$). Moreover there exist
a uniformizer $\pi$ of $D$ and a generator $\sigma$ of the Galois group
$\Gal(W/F)$ such that
\begin{equation}
\label{pi}
\pi w \pi^{-1}=\sigma(w) \mbox{ for all } w\in W.
\end{equation}
Then $\tau=\pi^d$ is a uniformizer of $F$, so
$\grm_D=\pi O_D$, $\grm_D\cap W=\grm_W=\tau O_W$ and $\grm_D\cap F=\grm_F=\tau O_F$.
Moreover we have the following direct sum decompositions:
\begin{gather}
D=W\oplus W\pi\oplus W\pi^2\oplus \dots\oplus W\pi^{d-1} \mbox{ as a left vector space over } W \\
\label{eq:directsum}
O_D=O_W\oplus O_W\pi\oplus O_W\pi^2\oplus \dots\oplus O_W\pi^{d-1} \mbox{ as a left module over } O_W
\end{gather}
As an immediate consequence of the second decomposition, the quotient
$O_D/\pi O_D$ is naturally isomorphic to $O_W/\pi^d O_W=O_W/\tau O_W=\Wbar$.
\skv
For each $n\in\dbZ$ the map $a\mapsto a\pi^n$ induces an isomorphism of abelian groups
$\Wbar \to \pi^n O_D/\pi^{n+1}O_D$. We will frequently identify $\pi^n O_D/\pi^{n+1}O_D$ with $\Wbar$
using this isomorphism.
\skv

The following formula for additive commutators in $O_D$ is an immediate consequence of \eqref{pi}:
\begin{equation}
\label{eq:Liebracket}
[a \pi^i,b \pi^j]=(a\sigma^i(b)-b\sigma^j(a))\pi^{i+j}
\mbox{ for }a,b\in W \mbox{ and }i,j\in\dbZ.
\end{equation}
\skv

\paragraph{\bf The norm one group $SL_1(D)$.}
Let $\Nred$ (resp. $\Tred$) denote the reduced norm (resp. reduced trace) map
from $D$ to $F$. One can define $\Nred(a)$ (resp. $\Tred(a)$) as the determinant (resp. trace) 
of the endomorphism of the left $W$-vector space $D$ given by $x\mapsto xa$. The restriction of $\Nred$ (resp. $\Tred$)
to $W$ coincides with the norm (resp. trace) map of the extension $W/F$.

It is straightforward to check that for $a\in W$ and $i\in\dbZ$ we have

\begin{equation}
\label{eq:trace}
\Tred(a\pi^i)=0 \,\,\,\mbox{ if }\,\,\,
d\nmid i \quad \mbox{ and }\quad 
\Tred(a\pi^i)=\tr_{W/F}(a)\tau^{\frac{i}{d}} \,\,\,\mbox{ if }\,\,\,
d\mid i.
\end{equation}

\begin{Lemma}
\label{lem:Riehm5} Let $n\in\dbZ_{\geq 0}$. Then
$\Tred(\pi^n O_D)=\tau^{\lceil \frac{n}{d}\rceil }O_F$ and
$\Nred(1+\pi^n O_D)=1+\tau^{\lceil \frac{n}{d}\rceil }O_F$
where $\lceil x\rceil$ is the smallest integer $\geq x$.
\end{Lemma}
\begin{proof} The first assertion follows immediately from 
\eqref{eq:trace}, \eqref{eq:directsum} and the fact that $\tr_{W/F}(O_W)=O_F$.
The second assertion is Lemma~5 in \cite{Ri}. 
\end{proof}

Let $G=SL_1(D)$ denote the group of elements of reduced norm one in $D$.
For $n\geq 1$ let $$G_n=SL_1^n(D)=\{g\in G : g\equiv 1\mod \pi^n O_D\}.$$
It is easy to check that each $G_n$ is a finite index pro-$p$ subgroup of $G$. 
\skv
For each $n\in\dbN$ define the map $\rho_n: 1+\pi^n O_D\to \Wbar$
by $$\rho_n(1+\pi^n a)=a\,\, {\rm mod}\,\, \pi O_D$$ (where we identify $O_D/\pi O_D$ with $\Wbar$ as before).
Note that $\rho_n$ is a group homomorphism and
$\Ker\rho_n = 1+\pi^{n+1} O_D$, so $\rho_n(G_n)\cong G_n/G_{n+1}$.
The following result appears as Lemma~6 in \cite{Ri}. 

\begin{Lemma}
\label{lem:Riehm6}
Let $n\in\dbN$. If $d\nmid n$, then $\rho_n(G_n)=\Wbar$. If
$d\mid n$, then $\rho_n(G_n)=\sl(\Wbar)$
where $\sl(\Wbar)=\{x\in \Wbar: \tr_{\Wbar/\Fbar}(x)=0\}$.
\end{Lemma}

\begin{Proposition}
\label{sl1basic} The following hold:
\begin{itemize}
\item[(a)] $G$ decomposes as a semi-direct product $G=G_1\rtimes \Delta$
where $\Delta$ is the group consisting of roots of unity in $W$
which have order prime to $p$ and norm $1$ over $F$.

\item[(b)] $[G_i,G_j]\subseteq G_{i+j}$ for any $i,j \geq 1$. In particular, 
$\gamma_k G_i\subseteq \gamma_{ki}G_1$ for all $i,k\geq 1$.
Assume now that $p\neq 2$ or $d\neq 2$.
Then
$$
[G_i,G_j]=
\left\{
\begin{array}{ll}
G_{i+j}&\mbox{ if } d\nmid (i+j) \mbox{ or } d\mid (i+j)  \mbox{ and } d\mbox{ is coprime to }i.\\
G_{i+j+1}&\mbox{ if } d\mid i\mbox{ and }d\mid j
\end{array}
\right.
$$
In particular, $G_i=\gamma_i G_1$ for all $i\geq 1$. 
\skv

\item[(c)] Assume that ${\rm char}\, F=0$, let $e$ be the ramification index of $F$, and let
$i> \frac{de}{p-1}$. Then $G_i^p=G_{i+de}$. Moreover, if $g\in G_i\backslash G_{i+1}$,
then $g^p\in G_{i+de}\backslash G_{i+de+1}$.
\end{itemize}
\end{Proposition}
\begin{Remark} Note that the above formula for $[G_i,G_j]$ does not cover all possible pairs $(i,j)$.
A complete description of $[G_i,G_j]$ in the remaining cases is given in the proof below, but
will not be needed for our purposes. 
\end{Remark}

\begin{proof}
(a) Let $\rho: O_D\to \Wbar$ be the natural projection (where we identify $O_D/\pi O_D$ with $\Wbar$ as before).
It is straightforward to check that for any $a\in O_D$
we have $\rho(N_{red}(a))=N_{\Wbar/\Fbar}(\rho(a))$
where $N_{\Wbar/\Fbar}$ is the norm of the extension $\Wbar/\Fbar$.

Thus, if $\overline\Delta=\{x\in \Wbar^*: N_{\Wbar/\Fbar}(x)=1\}$, then $\rho(G)\subseteq \overline\Delta$.
Note that $\Delta$ has trivial intersection with $\Ker\rho\cap G=G_1$ (e.g. since $G_1$ is a pro-$p$ group
and any element of $\Delta$ has finite order prime to $p$).
Thus, to prove (a) it remains to show that $\rho(\Delta)=\overline\Delta$.

Take any $x\in \overline\Delta$. By Hensel's lemma there exists a root of unity $a\in W$ 
such that $\rho(a)=x$. Moreover, the (multiplicative) order of $x$ is equal to the order of $a$
and thus is coprime to $p$. Then $N_{red}(a)$ is also a root of unity of order coprime to $p$;
on the other hand, $a\in \Ker\rho\cap O_D^*$ which is a pro-$p$ group. Thus, we must
have $N_{red}(a)=1$, so $a\in \Delta$.
\skv

(b) The inclusion $[G_i,G_j]\subseteq G_{i+j}$ can be verified by an easy direct computation. Assume now that
$(p,d)\neq (2,2)$. The second assertion of Proposition~\ref{sl1basic}(b) is a special case of \cite[Theorem~21]{Ri}; 
however, due to the technical nature of the statement of this theorem,
we will briefly explain how to obtain the result. First, \cite[Theorem~21]{Ri} asserts that
$[G_i,G_j]\supseteq G_{i+j+1}$, so the group $[G_i,G_j]$ is determined by its image in $G_{i+j}/G_{i+j+1}$.
The second assertion of Proposition~\ref{sl1basic}(b) can now be reformulated in terms of the map $\rho_{i+j}$:
we need to show that $\rho_{i+j}([G_i,G_j])=\rho_{i+j}(G_{i+j})$ if $d\nmid (i+j)$ or
$d\mid (i+j)$ and $d$ is coprime to $i$ and $\rho_{i+j}([G_i,G_j])=\{0\}$ if $d\mid i$ and $d\mid j$.

\skv
The equation \eqref{eq:Liebracket} easily implies that $\rho_{i+j}([G_i,G_j])$ is equal to the
set $S=\{\alpha \sigma^i(\beta)-\beta\sigma^j(\alpha): \alpha\in \rho_i(G_i),\beta\in \rho_j(G_j)\}$.
Here we use the fact that the natural map $Gal(W/F)\to Gal(\Wbar/\Fbar)$ is an isomorphism since
$W/F$ is unramified, and the image of $\sigma\in Gal(W/F)$ under this isomorphism is also denoted by $\sigma$.

Recall from Lemma~\ref{lem:Riehm6} that 
$\rho_k(G_k)=\Wbar$ if $d\nmid k$ and $\rho_k(G_k)=\sl(\Wbar)$ if $d\mid k$. It is clear
that if $d\mid i$ and $d\mid j$, then $S=\{0\}$ (since $\sigma^d=id$). If $d\mid (i+j)$ but $d\nmid i$,
we have $S=\{\alpha \sigma^i(\beta)-\beta\sigma^{-i}(\alpha): \alpha,\beta\in\Wbar\}=
\{\sigma^i(\gamma)-\gamma: \gamma\in \Wbar\}$. If in addition $i$ is coprime to $d$, then
$\sigma^i$ is a generator of $\Gal(\Wbar/\Fbar)$, so by Hilbert's Theorem~90
$S=\sl(\Wbar)$, and thus $S=\rho_{i+j}(G_{i+j})$. Finally, if $d\nmid (i+j)$,
$S=\Wbar$, as proved in \cite{Ri}.
\skv

(c) We first prove the second assertion. Assume that $i>\frac{de}{p-1}$ and take any $g\in G_i\setminus G_{i+1}$.
Thus $g=1+\pi^i a+b$ for some $a\in O_D^{*}=O_D\setminus \pi O_D$ and $b\in \pi^{i+1}O_D$. Then
$g^p=1+\sum\limits_{j=1}^p {p\choose j}(\pi^i a+b)^j=1+p\pi^{i}a+c$ where $c\in \pi^{pi} O_D + p\pi^{i+1}O_D$. Let $v:O_D\to \dbZ_{\geq 0}$
be the (additive) valuation on $D$ such that $v(\pi)=1$, that is, $v(x)=i$ if and only if $x\in \pi^i O_D\setminus \pi^{i+1}O_D$.
By definition of ramification index we have $v(p)=e\cdot v(\tau)=e\cdot v(\pi^d)=de$. Thus,
$v(p\pi^{i}a)=i+de$ and $v(c)\geq \min\{pi,i+1+de\}$. Since $i>\frac{de}{p-1}$, we have $pi=i+(p-1)i>i+de$, so
$i+de=v(p\pi^{i}a)<v(c)$ and hence
$v(g^p-1)=v(p\pi^{i}a+c)=i+de$. Therefore, $g^p\in G_{i+de}\setminus G_{i+de+1}$.
\skv

Now we prove the first assertion. We already proved that $G_{i}^p\subseteq G_{i+de}$.
The computation in the previous paragraph also shows that the map $\rho_{i}(G_i)\to \rho_{i+de}(G_i^p)$
induced by the $p$-power map $x\mapsto x^p$ is given by multiplication by a nonzero element of $\Wbar$ (namely
$\frac{p}{\tau^e}\,\,{\rm mod}\,\,\ \pi O_D$). In particular this means
that $\rho_{i}(G_i)=\rho_{i+de}(G_i^p)$. On the other hand, $\rho_{i+de}(G_{i+de})=\rho_{i}(G_i)$ by Lemma~\ref{lem:Riehm6}.

Thus, $\rho_{i+de}(G_i^p)=\rho_{i+de}(G_{i+de})$, so $G_{i+de}\subseteq G_i^p G_{i+de+1}$. Since this inclusion holds
for all $i> \frac{de}{p-1}$, arguing inductively, we deduce that $G_{i+de}\subseteq G_i^p G_{k}$ for all $k\in\dbN$. Since
$\{G_k\}$ is a base of neighborhoods of identity in $G$ (and since $G_i^p$ is closed), we deduce that
$G_{i+de}\subseteq G_i^p$ as desired.
\end{proof}

\skv
\paragraph{\bf Lie algebras of congruence subgroups and their quotients.}
If $n\geq de$, the group $G_n$ is powerful and torsion-free
by Proposition~\ref{sl1basic}(b)(c), so we can consider the Lie
algebra $\log(G_n)$ using Weigel's $\log$ functor.
Also for any $m,n\in\dbN$ such that $n\leq m\leq pn$,
we have $G_n/G_m\in Ob(\grG_{<p})$ by Proposition~\ref{sl1basic}(b), 
so we can consider $\log(G_n/G_m)$, this time using Lazard's $\log$ functor.

As one might expect, the Lie algebras obtained in this way are isomorphic
to $\grg_n$ in the first case and $\grg_n/\grg_m$ in the second case where
$$\grg_n=\sl(\pi^{n}O_D))=\{x\in \pi^{n}O_D: \Tred(x)=0\}.$$
Explicit isomorphisms between these Lie algebras are defined in our next proposition:

\begin{Proposition}
\label{prop:LA}
The following hold:
\begin{itemize}
\item[(a)] Let $n\geq de$. Let us define the maps $\Exp:\pi^n O_D\to 1+\pi^n O_D$ and 
$\Log:1+\pi^n O_D\to \pi^n O_D$ by 
$$\Exp(x)=\sum\limits_{i=0}^{\infty}\frac{x^i}{i!}
\mbox{ and } \Log(1+x)=\sum\limits_{i=1}^{\infty}\frac{(-1)^{i-1}x^i}{i}.$$
Then $\Exp$ and $\Log$ are mutually inverse bijections and $\Exp(\grg_n)=G_n$.
Moreover, if $\log(G_n)$ is the Lie algebra of $G_n$ where $\log$ is Weigel's functor
(recall that $\log(G_n)=G_n$ as sets), the map $x\mapsto \Exp(x)$ is a Lie algebra isomorphism
between $\grg_n$ and $\log(G_n)$.
\skv

\item[(b)] Assume that $1\leq n\leq m\leq pn$. Define the truncated exponential and logarithm mappings
$\Exp':\pi^n O_D/\pi^m O_D \to (1+\pi^n O_D)/(1+\pi^m O_D)$ and
$\Log':(1+\pi^n O_D)/(1+\pi^m O_D)\to \pi^n O_D/\pi^m O_D$ by 
$$\Exp'(x)=\sum\limits_{i=0}^{p-1}\frac{x^i}{i!} \mbox{ and } \Log'(1+x)=\sum\limits_{i=1}^{p-1}\frac{(-1)^{i-1}x^i}{i}.$$
Then $\Exp'$ and $\Log'$ are mutually inverse bijections and $\Exp(\grg_n/\grg_m)=G_n/G_m$.
Moreover, if $\log(G_n/G_m)$ is the Lie algebra of $G_n/G_m$ where $\log$ is Lazard's functor, 
the map $x\mapsto \Exp'(x)$ is a Lie algebra isomorphism between $\grg_n/\grg_m$ and $\log(G_n/G_m)$.
\skv

\item[(c)] Suppose that either $H=G_n$ and $\grh=\grg_n$ with $n\geq de$ or
$H=G_n/G_m$ and $\grh=\grg_n/\grg_m$ with $1\leq n\leq m\leq pn$. In both cases
the conjugation action of $G$ on $H$ induces an action of $G$ on $\log(H)$
by Proposition~\ref{exthom1}(a). Then under the isomorphisms $\grh\to \log(H)$
described in (a) and (b) this action corresponds to the conjugation action of
$G$ on $\grh$ given by
\begin{equation}
\label{LAaction1}
u^g=g^{-1}ug \mbox{ for } u\in\grg_n \mbox{ and }g\in G
\end{equation}
(and the corresponding induced action on $\grg_n/\grg_m$).
\end{itemize}
\end{Proposition}
\begin{Remark} Note that we capitalize the first letters of the maps $\Exp,\Log,\Exp'$ and $\Log'$ defined in (a) and (b)
to avoid confusion with the $\exp$ and $\log$ functors.
Also in the proof below we will use $\Exp,\Log,\Exp'$ and $\Log'$ not just for the maps
with the domain and codomain specified in Proposition~\ref{prop:LA}, but any maps given by the same formulas.
\end{Remark}
\begin{proof}
(a) Since $n\geq de$, it is easy to check that both $\Exp$ and $\Log$ are well defined, that is, the series converge and
all terms take values in the appropriate ideals (in fact, for this it would be enough to assume that $n>\frac{de}{p-1}$).
We know that $\Exp$ and $\Log$ are mutually inverse as power series under substitution, and a routine convergence check shows that
they are also mutually inverse as functions. To prove that $\Exp(\grg_n)=G_n$ it suffices to check
that $\Nred(\Exp(x))=\Exp(\Tred(x))$ and $\Tred(\Log(1+x))=\Log(\Nred(1+x))$ for all $x\in \pi^{de} O_D$.
By definition of reduced trace and reduced norm this amounts to showing that $\det(\Exp(A))=\Exp({\rm tr}(A))$ and
${\rm tr}(\Log(1+A))=\Log({\rm det}(1+A))$ for any $d\times d$ matrix $A$ over $\tau^d O_W$.
Both equalities are well-known over $\dbC$. However, since both sides of either equality are clearly power series in the entries of $A$,
they must coincide as power series, and hence the equalities holds over $\tau^d O_W$ as well.

Finally, the assertion that $x\mapsto \Exp(x)$ is a Lie algebra isomorphism between $\grg_n$ and $\log(G_n)$
is a direct consequence of Theorem~\ref{CHF} (recall that the group operation on $\log(G_n)$ is defined using the
BCH series $\Phi$).

\skv
(b) This time we do not have to worry about convergence, and since the coefficients in the expressions for
$\Exp'$ and $\Log'$ are $p$-adic integers, $\Exp'$ and $\Log'$ are well defined functions between
$\pi^n O_D/\pi^m O_D$ and $(1+\pi^n O_D)/(1+\pi^m O_D)$. 

Let us check that $\Exp'$ and $\Log'$ are mutually inverse.
It is clear from definitions that the polynomial $\Log'\circ \Exp'(x)$ coincides with the power series
$\Log\circ \Exp(x)=x$ up to (and including) degree $p-1$. 
Therefore 
$\Log'\circ \Exp'(x)=x+\sum\limits_{k=p}^{(p-1)^2}c_k x^k$
where $c_k$ are $p$-adic integers. Hence for all $x\in \pi^n O_D$ we have
 $\Log'\circ \Exp'(x)-x\in \pi^{pn} O_D\subseteq \pi^m O_D$ since $m\leq pn$. Thus,
 $\Log'\circ \Exp'(x)=x$ for all $x\in \grg_n/\grg_m$, as desired. The proof for the composition
in reverse order is analogous.
\skv
Next we prove that $\Exp'(\grg_n/\grg_m)=G_n/G_m$. The argument in (a) shows that if $A$ is a $d\times d$ matrix with generic entries 
$(a_{ij})$, then $\det(\Exp(A))=\Exp({\rm tr}(A))$ and ${\rm tr}(\Log(1+A))=\Log({\rm det}(1+A))$ as power series 
(say, over $\dbQ$) in formal commuting variables $\{a_{ij}\}$. If we now replace $\Exp$ and $\Log$
by $\Exp'$ and $\Log'$, both sides of either equality above will not change modulo $I^p$
where $I$ is the ideal of $\dbQ[[\{a_{ij}\}]]$ generated by $\{a_{ij}\}$.
Therefore $\det(\Exp'(A))\equiv\Exp'({\rm tr}(A)){\,\,\rm mod\,\,} I^p$ and  
${\rm tr}(\Log'(1+A))\equiv \Log'({\rm det}(1+A)){\,\,\rm mod\,\,} I^p$.

It follows that $\det(\Exp'(A))=\Exp'({\rm tr}(A))$ and
${\rm tr}(\Log'(1+A))=\Log'({\rm det}(1+A))$ for all
$A\in Mat_d(\tau^r O_D/\tau^s O_D)$ whenever $1\leq r\leq s\leq pr$. As in the proof of (a) this implies
that $\Exp'(\grg_n/\grg_m)=G_n/G_m$ whenever $1\leq n\leq m\leq pn$.
\skv
Let us now prove the last assertion of (b). Recall that the group
operation on $\exp(\grg_n/\grg_m)$ is defined by
$u\cdot v=\Phi'(u,v)$ where $\Phi'$ is the BCH series truncated after degree $p-1$.
By Theorem~\ref{CHF} we have $\Exp(x)\cdot \Exp(y)=\Exp (\Phi(x,y))$ where both sides
are treated as formal power series in non-commuting variables $x$ and $y$.
If we replace $\Exp$ by $\Exp'$ and $\Phi$ by $\Phi'$, both
sides will still be congruent modulo the terms of degree $\geq p$.
Since $\Exp'$ and $\Phi'$ are polynomials whose coefficients are $p$-adic integers,
arguing as the first part of the proof of (b), we deduce that
$\Exp'(u)\cdot \Exp'(v)=\Exp' (\Phi'(u,v))$ for all 
$u,v\in \grg_n/\grg_m$. This means precisely that
$\Exp'$ is a group isomorphism from
$\exp(\grg_n/\grg_m)$ to $G_n/G_m$. Finally, 
by Theorem~\ref{Weigel0} we deduce that $\Exp'$
is also a Lie algebra isomorphism from
$\log(\exp(\grg_n/\grg_m))=\grg_n/\grg_m$ to
$\log(G_n/G_m)$, as desired.

\skv
(c) follows immediately from definitions and the fact
the conjugation actions of $G$ on $G_n$, $G_n/G_m$,
$\grg_n$ and $\grg_n/\grg_m$ commute with 
the maps $\Exp,\Log,\Exp'$ and $\Log'$ (which is clear from the explicit formulas for these maps).
\end{proof}

\subsection{More on local fields.}
We finish this section with a collection of elementary facts which will come very handy when 
we compute the cohomology of $\grg_n$ in \S~\ref{sec:cohomology}.

Given a local field $K$, we denote by $\bfk$ the abelian group $K/O_K$ (and use similar notations with $K$ replaced by other capital Latin letters).
Clearly, $\bfk$ is an $O_K$-module. Since all automorphisms of local fields are continuous, any field automorphism $\phi:K\to K$ preserves $O_K$ and thus induces a map $\bfk\to\bfk$ which will also be denoted by $\phi$.

\begin{Claim}
Let $L/K$ be an unramified extension of local fields. Define the trace map
$\tr_{\bfl/\bfk}: \bfl\to \bfk$ by $$\tr_{\bfl/\bfk}(x+O_L)=\tr_{L/K}(x)+O_K$$
Let $\sl(L)=\{x\in L: \tr_{L/K}(x)=0\}$ and 
$\sl(\bfl)=\{x\in \bfl: \tr_{\bfl/\bfk}(x)=0\}$. 

The following hold:
\begin{itemize}
\label{torring}
\item[(a)] $\tr_{\bfl/\bfk}(\phi(\alpha))=\tr_{\bfl/\bfk}(\alpha)$ for all $\alpha\in \bfl$ and $\phi\in \Gal(L/K)$.
\item[(b)] Let $\alpha\in \bfl$. Then $\alpha$ is fixed by all $\phi\in \Gal(L/K)$
if and only if $\alpha\in \bfk$.
\item[(c)] Let $\alpha\in \bfl$ and $\sigma$ a generator of $\Gal(L/K)$.
Then $\alpha\in \sl(\bfl)$ if and only if $\alpha=\beta-\sigma(\beta)$ for some $\beta\in\bfl$.
\item[(d)] Let $\pi$ be a uniformizer of $K$.
Then $\tr_{\bfl/\bfk}$ maps $\frac{1}{{\pi}^n}O_L/O_L$ onto $\frac{1}{{\pi}^n}O_K/O_K$ for every $n\in\dbN$.
\item[(e)] If $\alpha\in \bfl$ is such that $\tr_{\bfl/\bfk}(\alpha u)=0$ for all $u\in O_L$, then
$\alpha=0$.
\item[(f)] Let ${\mathbf n}=\{\alpha\in \bfl: \tr_{\bfl/\bfk}(\alpha uv)=0 \mbox{ for all }u,v\in \sl(O_L)\}$.
If $d>2$, then ${\mathbf n}=\{0\}$, and if $d=2$, then ${\mathbf n}=\sl(\grl)$.
 \item[(g)] Let $M$ be another local field with $M\subseteq K$ such that $K/M$ is unramified, and let 
$\alpha\in \bfl$ be such that $\tr_{\bfl/\bfm}(\alpha u)=0$ for all $u\in O_L$ with $\tr_{L/K}(u)=0$.
Then $\alpha\in \bfk$.
\end{itemize}
\end{Claim}
\begin{proof}All parts of this claim can be proved using the same principle. For every $n\in\dbN$ and a local field $M$ we set 
$R_n(M)=\frac{1}{{\pi}^n}O_M/O_M$ where $\pi$ is a uniformizer of $M$. Then each $R_n(M)$ is an $O_M$-submodule and each quotient $Q_n(M)=R_n(M)/R_{n-1}(M)$ can be naturally identified with
$\Mbar$, the residue field of $M$. Also, under such identification, for any unramified extension of local fields $L/K$,
the induced map $\Gal(L/K)\to \Aut (Q_n(L))$ is an isomorphism onto $\Gal(\Lbar/\Kbar)$.

This way we can reduce each of the parts (a)-(g) to the corresponding result about finite fields. The finite field analogues of (a)-(e) are well known. We will give a proof for the finite field counterpart of (g), with the proof of (f) being similar.
But first let us give a full argument for (c) as an illustration of the general principle. 

The backwards direction in (c) holds by (b). Suppose now that $\alpha\in \sl(\bfl)$.
We know that $\alpha\in R_n(L)$ for some $n$. Since $\alpha\in \sl(\bfl)$,
the projection of $\alpha$ to $Q_n(L)$, call it $\pi_n(\alpha)$, satisfies $\tr_{\Lbar/\Kbar}(\pi_n(\alpha))=0$. By Hilbert's Theorem~90
(applied to the cyclic extension $\Lbar/\Kbar$) there exists $\beta\in R_n(L)$ such that $\pi_n(\beta-\sigma(\beta))=\pi_n(\alpha)$.
Let $\alpha'=\alpha-(\beta-\sigma(\beta))$. Then $\alpha'$ lies in $ \sl(\bfl)\cap R_{n-1}(L)$.
Applying the same argument to $\alpha'$ and repeating the process $n-1$ times, we deduce that $\alpha=\gamma-\sig(\gamma)$
for some $\gamma\in \bfl$.
\skv

Let us now prove the analogue of (g) for finite fields which asserts the following. 
\it Let $F_1\subseteq F_2\subseteq F_3$ be a tower of finite fields. Let $\alpha\in F_3$ be such that $\tr_{F_3/F_1}(\alpha\beta)=0$ for all $\beta\in \sl_{F_2}(F_3)$ where
$\sl_{F_2}(F_3)=\{\beta\in F_3: \tr_{F_3/F_2}(\beta)=0\}$. Then $\alpha\in F_2$.\rm
\skv

We first claim that $\tr_{F_3/F_1}(u\beta)=0$ for all $\beta\in \sl_{F_2}(F_3)$ and $u\in F_2\alpha+F_2$.
Indeed, write $u=f\alpha+f'$ with $f,f'\in F_2$. Then $f'\beta, f\beta \in  \sl_{F_2}(F_3)$ as well, so
$\tr_{F_3/F_1}(f'\beta)=\tr_{F_2/F_1}(\tr_{F_3/F_2}(f'\beta))=0$,  $\tr_{F_3/F_1}(\alpha\cdot f\beta)=0$ by assumption on $\alpha$ and hence
$\tr_{F_3/F_1}(u\beta)=\tr_{F_3/F_1}(f\alpha\beta +f'\beta)=0$ as well.

Since $\tr_{F_3/F_1}:F_3\times F_3\to F_1$ is a non-degenerate bilinear form and $\dim_{F_1}(\sl_{F_2}(F_3))=[F_3:F_1]-[F_2:F_1]$,
its orthogonal complement must have $F_1$-dimension $[F_2:F_1]$. Hence $\dim_{F_1}(F_2\alpha+F_2)\leq [F_2:F_1]$
whence $F_2\alpha =F_2$ and so $\alpha\in F_2$.
\end{proof}

\section{Group-theoretic structure of central extensions of $SL_1(D)$}
For the next three sections we fix a $p$-adic field $F$ and a central
division algebra $D$ over $F$. We preserve all notations from
\S~\ref{sec:norm1}. Recall that $G=SL_1(D)$, $d$ is the degree of $D$
and $e$ is the ramification index of $F$. Throughout this section we shall assume 
that $(p,d)\neq (2,2)$.

Let $A_{\infty}$ denote the group $\dbQ_p/\dbZ_p$ (note that $\dbQ_p/\dbZ_p$
is isomorphic to the $p$-primary component of $\dbQ/\dbZ$). Given $n\in\dbN$,
let $A_n$ be the group of elements of order $\leq p^n$ in $A_{\infty}$
(of course, $A_n$ is simply a cyclic group of order $p^n$,
but it will be convenient to think of it as a subgroup of $\dbQ_p/\dbZ_p$).
The symbol $A$ will denote $A_n$ for some $n$ when the value of $n$
is not important.

The following key observation appears in \cite[2.2]{PR2}: 

\begin{Lemma}
\label{lem:PR}
Let $k\in\dbN$, and let $\iota_k: H^2(G,A_k)\to H^2(G,A_{\infty})$
be the homomorphism induced by the embedding map $A_k\to A_{\infty}$.
Then 
\begin{itemize}
\item[(a)] each $\iota_k$ is injective  
\item[(b)] $\iota_k(H^2(G,A_k))$ is the subgroup
of elements of order $\leq p^k$ in $H^2(G,A_{\infty})$.
\item[(c)] $H^2(G,A_{\infty})=\bigcup\limits_{k=1}^{\infty}\iota_k(H^2(G,A_k)).$
\end{itemize}
\end{Lemma}
\begin{proof}
We only need to use two facts about $G$:
\begin{itemize}
\item[(i)] $G$ is profinite.
\item[(ii)] $G/[G,G]$ is finite of order prime to $p$ (this holds by Proposition~\ref{sl1basic}(a)).
\end{itemize}
In view of (ii), (a) and (b) follow from the long exact
sequence of cohomology groups associated to the short exact sequence
$1\to A_{k}\to A_{\infty}\str{\times p^{k}}{\lra} A_{\infty}\to 1$.

Since $G$ is profinite and $A_{\infty}$ is discrete, any continuous map $G\times G\to A_{\infty}$
has finite image and thus lies in $A_k$ for some $k$. This implies (c).
\end{proof}
\begin{Remark}\rm 
In view of Lemma~\ref{lem:PR} from now on we will identify each $H^2(G,A_k)$ with the subgroup
of elements of order $\leq p^k$ in $H^2(G,A_{\infty})$.
\end{Remark}
\skv

The main result of \cite{PR2} asserts
that $H^2(G,A_{\infty})$ is a finite cyclic group of $p$-power order.

\begin{Corollary}
\label{cor:PR}
Each group $H^2(G,A_k)$ is cyclic of $p$-power order and (with the above identification)
$H^2(G,A_{\infty})= H^2(G,A_k)$ for all sufficiently large $k$.
\end{Corollary}
\skv
In this section we study group-theoretic properties of central extensions of $G$ 
by $A_k$ for $k\in\dbN$.
Throughout this section we write $G_n=SL_1^n(D)$ for $n\geq 1$ and set $S=G_1$.
The use of the letter $S$ is ``justified'' by the fact that $S$ is the Sylow pro-$p$ subgroup of $G$.
Recall that $G_n=\gamma_n S$ for $n\geq 1$ by Proposition~\ref{sl1basic}(b).
\vskip .12cm

The following proposition describes the basic power-commutator structure in
extensions of $G$ by $A$:

\begin{Proposition}
\label{powercomm}
Let $c\in H^2(G,A)$ with $\Ext(c)=\ex{A}{\Ghat}{G}{\iota}{\phi}$.
Let $\Shat=\phi^{-1}(S)$ and $\Ghat_k=\phi^{-1}(G_k)$ for $k\in\dbN$.
The following hold:
\begin{itemize}
\item[(a)] $\gamma_{k+de}\Shat=(\gamma_k\Shat)^p$ for any $k>\frac{de}{p-1}+1$;

\item[(b)] For any $k\geq 1$ we have $\gamma_{2k+1+\delta}\Shat\subseteq \gamma_2 \Ghat_k\subseteq \gamma_{2k}\Shat$
where $\delta=0$ if $d\nmid k$ and $\delta=1$ if $d\mid k$.

\item[(c)] Let $x\in G_k\backslash G_{k+1}$ for some $k>\frac{de}{p-1}+1$, and choose
any $\xhat\in \gamma_k \Shat$ such that $\phi(\xhat)=x$. Then
$\xhat^{p^n}\in \gamma_{k+nde} \Shat\backslash \gamma_{k+1+nde} \Shat$ for any $n\geq 0$.
\end{itemize}
\end{Proposition}
\proofs
The following property will be used several times in the computation below.
\emph{If $U$ and $V$ are subgroups of $\Ghat$ such that $\phi(U)=\phi(V)$,
and $W$ is another subgroup of $\Ghat$, then $[U,W]=[V,W]$.} This holds since by assumption
$\Ker\phi$ is central, so the equality $\phi(U)=\phi(V)$ implies that $U\subseteq V\cdot \Ker\phi\subseteq V\cdot Z(G)$
and similarly $V\subseteq U\cdot Z(G)$.

We will also use the well-known Hall-Petrescu formula~\cite[Appendix~A]{DDMS} which asserts the following:
let $a$ and $b$ be elements of the same group and $H=\la a,b\ra$. Then for any integer $m\geq 2$
there exist elements $c_2,\ldots, c_m$ (all depending on $m$) such that $c_i\in \gamma_i H$ for all $i$ and
$(ab)^m=a^m b^m \prod\limits_{i=2}^m c_i^{m\choose i}$.

If $m$ is prime, letting $a=x^{-1}$ and $b=y^{-1}xy$ where $x$ and $y$ are elements of the same group,
and using the fact that $\la a,b\ra=\la x, [x,y]\ra$
we deduce the following:
\begin{equation}
\label{eq:HP}
[x,y]^p \equiv [x^p,y] \mod K
\end{equation}
where $K$ is the subgroup normally generated by the $p^{\rm th}$ powers of commutators in $x$ and $y$ involving at least
two occurrences of $x$ and commutators in $x$ and $y$ involving at least $p$ occurrences of $x$.
\skv

We now proceed with the proof.

(a) By the Hall-Petrescu formula, every element of $(\gamma_{k} \Shat)^p$ is a product of 
an element of $\gamma_{pk}\Shat$ and elements of the form $[u,v]^p$ with $u\in\gamma_{k-1}\Shat$ and $v\in\Shat$,
and by \eqref{eq:HP}, every such $[u,v]^p$ lies in $[(\gamma_{k-1}\Shat)^p,\Shat](\gamma_{2k-1}\Shat)^p \gamma_{p(k-1)+1}\Shat$. Thus,
$$(\gamma_{k} \Shat)^p\subseteq  [(\gamma_{k-1}\Shat)^p,\Shat](\gamma_{2k-1}\Shat)^p \gamma_{p(k-1)+1}\Shat.$$
Replacing $k$ by $2k-1$ in the above inclusion, then by $2(2k-1)-1=4k-3$ etc., after finitely many steps
we conclude that
\begin{equation}
(\gamma_{k} \Shat)^p\subseteq  [(\gamma_{k-1}\Shat)^p,\Shat]\,\gamma_{p(k-1)+1}\Shat.
\label{HP}
\end{equation}
By Proposition~\ref{sl1basic}, $(\gamma_{k-1}S)^p=G_{k-1}^p=G_{k+de-1}=\gamma_{k+de-1} S$.
Hence $\phi((\gamma_{k-1}\Shat)^p)=\phi(\gamma_{k+de-1} \Shat)$, so 
$ [(\gamma_{k-1}\Shat)^p,\Shat]=[\gamma_{k+de-1}\Shat,\Shat]=\gamma_{k+de}\Shat$. 
It follows from (\ref{HP}) that
$(\gamma_{k} \Shat)^p\subseteq \gamma_{\min(k+de, p(k-1)+1)}\Shat$.
Since $k> \frac{de}{p-1}+1$, we have $k+de<p(k-1)+1$, whence
$(\gamma_{k} \Shat)^p\subseteq\gamma_{k+de}\Shat$.

The reverse inclusion
$\gamma_{k+de}\Shat\subseteq (\gamma_{k} \Shat)^p$ is proved in a similar fashion:
as we already showed, $\gamma_{k+de}\Shat= [(\gamma_{k-1}\Shat)^p,\Shat]$,
and by the Hall-Petrescu formula and \eqref{eq:HP} we have
$$
\gamma_{k+de}\Shat= [(\gamma_{k-1}\Shat)^p,\Shat]\subseteq (\gamma_{k} \Shat)^p (\gamma_{2k-1}\Shat)^p \gamma_{p(k-1)+1}\Shat=
(\gamma_{k} \Shat)^p \gamma_{p(k-1)+1}\Shat.
$$
Since $p(k-1)+1>k+de$, we showed that the subgroups $U=\gamma_{k+de}\Shat$ and $V=(\gamma_{k} \Shat)^p$ of $\Shat$
satisfy $U\subseteq V[U,\Shat]$. Since $V$ is closed and $\Shat$ is pro-$p$, this easily implies
that $U\subseteq V$, as desired.
\vskip .12cm

(b) Since $\phi(\Ghat_k)=\phi(\gamma_k \Shat)$, we have
$\gamma_2\Ghat_k= [\Ghat_k,\gamma_k \Shat]= [\gamma_k \Shat,\gamma_k \Shat]$,
whence $\gamma_2 \Ghat_k\subseteq \gamma_{2k}\Shat$.
By Proposition~\ref{sl1basic}(b) we have $[G_k,G_k]=G_{2k+\delta}$, whence
$\phi( [\Ghat_k,\Ghat_k])= [\gamma_k S,\gamma_k S]=
\gamma_{2k+\delta}S=\phi(\gamma_{2k+\delta} \Shat)$. Therefore,
$$\gamma_2\Ghat_k\supseteq  [\gamma_2\Ghat_k, \Shat]=
 [\gamma_{2k+\delta} \Shat,\Shat]= \gamma_{2k+1+\delta} \Shat.$$
\vskip .12cm
(c) First note that an element $\xhat$ with required properties always exists
since $G_k=\gamma_k S=\phi(\gamma_k \Shat)$.
By part (a) we have $\xhat^{p^n}\in \gamma_{k+nde} \Shat$.
Now suppose that $\xhat^{p^n}\in \gamma_{k+1+nde} \Shat$. Then 
$x^{p^n}=\phi(\xhat^{p^n})\in \phi(\gamma_{k+1+nde} \Shat)=G_{k+1+nde}$
which contradicts Proposition~\ref{sl1basic}(c) since $x\not\in G_{k+1}$ by assumption.
\QED
\skv
\paragraph{\bf Depth and commutator breaks.}
\skv
Given $c\in H^2(G,A)$, there are two natural ways to measure
the ``complexity'' of the associated extension which lead to
the notions of inflation depth and commutator depth of $c$. 
However, we will show (see Proposition~\ref{depths} below) 
that the two notions of depth always coincide.

\begin{Definition}\rm
Let $A=A_k$ for some $k\geq 1$.
Let $c$ be an element of $H^2(G,A)$ with
$\Ext(c)=\ex{A}{\Ghat}{G}{\iota}{\phi}$.
Let $\Shat=\phi^{-1}(S)$.

\begin{itemize}
\item
The \it{inflation depth }\rm of $c$, denoted by $\infdep(c)$,
is the smallest integer $m$ such that $c$ lies in the image
of the inflation map $\inf: H^2(G/G_m,A)\to H^2(G,A)$.
\item
An integer $m>1$ will be called a \it{commutator break }\rm of $c$
if $$\Ker\phi \cap\gamma_m \Shat\neq \Ker\phi \cap\gamma_{m+1} \Shat.$$
\item
The \it{commutator depth }\rm of $c$, denoted by $\comdep(c)$,
is the largest integer $m$ such that $\Ker\phi \cap\gamma_m \Shat\neq \{1\}$.
Thus, $\comdep(c)$ is the largest commutator break of $c$ if
there is at least one break, and $\comdep(c)=1$ if $c$ has no breaks.
\end{itemize}
\end{Definition}

\begin{Proposition}
\label{depths}
For any $c\in H^2(G,A)$
we have $\infdep(c)=\comdep(c)$.
\end{Proposition}

\begin{Remark}\rm 
1. Recall that we identify $H^2(G,A_r)$ with its image in $H^2(G,A_s)$
for any $s\geq r$. This leads to a potential ambiguity in the definitions of commutator
breaks and inflation depth, but fortunately the problem does not occur. The fact that
commutator breaks do not depend on the choice of $s$ above will be proved in Lemma~\ref{comdepinv}
below, and the analogous result for the inflation depth is a direct consequence of
Lemma~\ref{comdepinv} and Proposition~\ref{depths}. Until such independence has been established,
we will occasionally use the notation $\infdep(c,A_s)$ for the inflation depth of
some $c\in H^2(G,A_r)$ considered as an element of $H^2(G,A_s)$ with $s\geq r$. It is clear
from the definitions that $\infdep(c,A_s)\geq \infdep(c,A_{s'})$ if $s\leq s'$.
\skv

2. The definitions of commutator and inflation depth make perfect sense for any profinite group $G$
with an open pro-$p$ subgroup $S$, where in the definition of inflation depth we set $G_i=\gamma_i S$.
As our proof of Proposition~\ref{depths} will show, the inequalities $\infdep(c)-1\leq \comdep(c)\leq \infdep(c)$
will hold in general, but in order to rule out the possibility that $\comdep(c)=\infdep(c)-1$ we will need to use
rather specific information about $G=SL_1(D)$. 
\end{Remark}
\skv

In order to prove Proposition~\ref{depths}, we will use the following deep result from \cite{PR2}. 

\begin{Lemma}{\rm (see \cite[Theorem~7.1]{PR2})}.
If $F$ has no primitive $p^{\rm th}$ root of unity, the group $H^2(G,A_{1})$ is trivial
(and hence $H^2(G,A_{\infty})$ is also trivial by Lemma~\ref{lem:PR}(b)(c)).
Otherwise, $H^2(G,A_{1})$ has order $p$ and for any non-trivial $c\in H^2(G,A_1)$ one has $$\infdep(c,A_s)= \frac{pde}{p-1}$$ 
for sufficiently large $s$.
\label{infdep} \QED
\end{Lemma} 

\begin{Remark}\rm 
1. The existence of a primitive $p^{\rm th}$ root of unity in $F$
implies that $p-1$ divides $e$. 

2. The statement in \cite{PR2} asserts that any non-trivial $c\in H^2(G,A_1)$ has inflation depth $\frac{pde}{p-1}$
as an element of $H^2(G,A_{\infty})$. However, the above formulation is clearly equivalent since
any element of $H^2(G/G_m,A_{\infty})$ comes from $H^2(G/G_m,A_{s})$ for some $s\in\dbN$, as $G/G_s$ is finite
(see the proof of Lemma~\ref{lem:PR}(c)).

3. For our proof of Proposition~\ref{depths} to go through it would be sufficient to know that $\infdep(c)>2$
for any $0\neq c\in H^2(G,A_1)$. However, we will need the full power of Lemma~\ref{infdep} later in the paper anyway.
\end{Remark}
\skv

\begin{Corollary}
\label{corcom}
Let $c$ be a non-trivial element of $H^2(G,A_r)$ for some $r$.
Then $$\infdep(c,A_s)\geq \frac{pde}{p-1} \mbox{ for all }s\geq r.$$
\end{Corollary}

\proofs
Suppose that $c$ has order $p^k$. Fix $s\geq r$, and let $m=\infdep(c,A_s)$. Thus $c$ is the inflation image of some $c_1\in H^2(G/G_m,A_s)$. Then $p^{k-1}c_1$ maps to $p^{k-1}c$, so $\infdep(p^{k-1}c,A_s)\leq m$ and hence $\infdep(p^{k-1}c,A_t)\leq m$ for sufficiently large $t$.

On the other hand,  $p^{k-1}c$ has order $p$ and thus lies in $H^2(G,A_1)$. Therefore $\infdep(p^{k-1}c,A_t)=\frac{pde}{p-1}$ for sufficiently
large $t$ by Lemma~\ref{infdep}, so $m\geq\frac{pde}{p-1}$, as desired.
\QED
\skv

\proofs[Proof of Proposition~\ref{depths}]

Let $\ex{A}{\Ghat}{G}{\iota}{\phi}$ be the extension corresponding to $c$.
\skv

We first prove that $\comdep(c)\leq \infdep(c)$.
Let $n=\infdep(c)$. By Lemma~\ref{cohbasic}(a), there exists a continuous section $\psi:G\to\Ghat$ such that
$\psi(G_n)$ is a normal subgroup of $\Ghat$. Since $G_n=\gamma_n S$, we have
$\phi(\psi(G_n))=\phi(\gamma_n \Shat)$. Since $\Ker\phi$ is central in $\Ghat$, we have
$[\psi(G_n),\Shat]=[\gamma_n \Shat,\Shat]=\gamma_{n+1}\Shat$, and since $\psi(G_n)$ is normal,
we conclude that $\gamma_{n+1}\Shat\subseteq \psi(G_n)$. Since $\Im\psi$ cannot contain non-trivial
elements of $\Ker\phi$, it follows that $\gamma_{n+1}\Shat\cap \Ker\phi=\{1\}$. Hence
$\comdep(c)\leq n$ by the definition of commutator depth.
\skv

Next we prove that $\comdep(c)\geq \infdep (c)-1$.
Let $m=\comdep(c)$.
Then
$\Ker\phi\,\cap\,\gamma_{m+1} \Shat=\{1\}$, whence
$\phi$ maps $\gamma_{m+1} \Shat$ isomorphically
onto $\gamma_{m+1} S=G_{m+1}$. Therefore, $\phi$ has
a continuous section $\psi$ such that $\psi(G_{m+1})=\gamma_{m+1} \Shat$.
Since $\Shat$ is normal in $\Ghat$, so is $\gamma_{m+1} \Shat$
and therefore $\infdep(c)\leq m+1$ by Lemma~\ref{cohbasic}(b).
Thus, we showed that $m\geq \infdep(c)-1$.

\skv
Finally, we will show that $\comdep(c)\geq \infdep (c)$ (thereby finishing the proof) using more specific information about $G=SL_1(D)$.
If $c=0$, Proposition~\ref{depths} is trivially true, so
from now on we assume that $c\neq 0$.  
By Corollary~\ref{corcom} we have $m\geq \frac{pde}{p-1}-1$,
whence $(\gamma_m \Shat)^p\subseteq \gamma_{m+1} \Shat$
by Proposition~\ref{powercomm}(a).
We can consider $\gamma_m\Shat/\gamma_{m+1} \Shat$ as a vector
space over $\Fp$, with the action of $\Delta$, identifying
$\Delta$ with $G/S\cong \Ghat/\Shat$ (recall from Proposition~\ref{sl1basic}(a) that
$\Delta$ is the group of roots of unity of order prime to $p$ in $G\cap W^*$).

Let $\Kbar$ be the subspace $(\gamma_m\Shat\cap \Ker\phi)\gamma_{m+1} \Shat/\gamma_{m+1} \Shat$
of $\gamma_m\Shat/\gamma_{m+1} \Shat$. Clearly, $\Kbar$ is $\Delta$-invariant,
and since $\Delta$ is a finite group of order prime to $p$, we can find
a $\Delta$-invariant vector subspace $\Lbar$ of $\gamma_m\Shat/\gamma_{m+1} \Shat$
such that $\gamma_m\Shat/\gamma_{m+1} \Shat=\Lbar\oplus\Kbar$.
Let $H$ be the full preimage of $\Lbar$ under the projection $\gamma_m\Shat\to \gamma_m\Shat/\gamma_{m+1} \Shat$.
Then $H$ lies between $\gamma_m \Shat$ and $\gamma_{m+1} \Shat$,
so $H$ is automatically normal in $\Shat$. Moreover, $H$ is $\Delta$-invariant,
so $H$ is normal in $\Ghat$. By construction,
$H\cap\Ker\phi=\{1\}$ and $\phi(H)=\gamma_m S=G_m$, so $\phi$ has
a continuous section $\psi'$ such that $\psi'(G_{m})=H$.
Applying Lemma~\ref{cohbasic}(b) as earlier in this proof, we conclude that $\infdep(c)\leq m=\comdep(c)$.
\QED
\skv
The following lemma combined with Proposition~\ref{depths} ensures that there is no ambiguity in the definition
of commutator breaks and inflation depth:

\begin{Lemma}
\label{comdepinv}
Let $c\in H^2(G,A_r)$ for some $r$, and let $s>r$. 
If $c'\in H^2(G,A_s)$ is the image of $c$ under the natural
mapping $H^2(G,A_r)\to H^2(G,A_s)$,
then $c$ and $c'$ have the same set of commutator breaks.
In particular, $\comdep(c')=\comdep(c).$
\end{Lemma}
\proofs
Let $Z\in Z^2(G,A_r)$ be a cocycle representing $c$.
Suppose that
$\Ext(c)=\ex{A_r}{\Ghat}{G}{}{\phi}$ and $\Ext(c')=\ex{A_s}{\Ghat'}{G}{}{\phi'}$.
We can assume that $\Ghat=A_r\times G$, $\Ghat'=A_s\times G$ as topological spaces, both $\phi$ and $\phi'$
are the projections onto the second component and the multiplication in both
$\Ghat$ and $\Ghat'$ is given by \eqref{eq:cocycle} (using the same cocycle $Z$
in both cases). 

If $\Shat=\phi^{-1}(S)$ and $\Shat'=\phi'^{-1}(S)$, then 
$\Shat=A_r\times S$ and $\Shat'=A_s\times S$, so $\Shat$ is a subgroup
of $\Shat'$, and $\Shat'$ is the product of $\Shat$ and 
the center of $\Shat'$. Thus for all $m\geq 2$ we have
$\gamma_m \Shat'=\gamma_m \Shat$ and hence
$\Ker\phi'\cap \gamma_m \Shat'=\Ker\phi \cap \gamma_m \Shat$
(since $\Ker\phi=\Ker\phi'\cap \Ghat$). Hence by definition
$c$ and $c'$ have the same commutator breaks.
\QED

\skv
Now we are ready to prove a formula for commutator breaks.

\begin{Proposition}
\label{combreaks}
Let $c\in H^2(G,A_s)$ for some $s\in\dbN$, and suppose that $ord(c)=p^n$.
Then $c$ has exactly $n$ commutator breaks $b_1<\ldots< b_n$, and moreover
$b_i=de(i+\frac{1}{p-1})$ for all $i$. In particular, $\comdep(c)=de(n+\frac{1}{p-1})$.
\end{Proposition}
\begin{Remark}\rm An essentially equivalent statement was proved earlier by Prasad (unpublished)
using a different method.
\end{Remark} 
\proofs
By Lemma~\ref{comdepinv} and Lemma~\ref{lem:PR}(b) we can assume that $s=n$.

We will first show that $b_1=de(1+\frac{1}{p-1})=\frac{pde}{p-1}$, then prove that
$c$ has exactly $n$ commutator breaks and finally prove the formulas for the remaining breaks.
In all parts of the proof we let $\ex{A_n}{\Ghat}{G}{\iota}{\phi}$ be the extension corresponding to $c$ and 
$\Shat=\phi^{-1}(S)$.
\skv

{\it Step 1:} Consider the exact sequence $1\to A_{n-1}\to A_{n}\str{\times p^{n-1}}{\lra} A_{1}\to 1$,
let $H^2(G,A_{n-1})\to H^2(G,A_{n})\to H^2(G,A_{1})$ be the induced maps between the cohomology groups,
and let $\cbar$ be the image of $c$ in $H^2(G,A_1)$. Since $c$ has order $p^n$, it does not come from 
$H^2(G, A_{n-1})$, and hence $\cbar$ is non-trivial. Hence
$\comdep(\cbar)=\infdep(\cbar)=\frac{pde}{p-1}$ by Proposition~\ref{depths} and Lemma~\ref{infdep}.

It is easy to check that $\Ext(\cbar)=\ex{A_1}{\Ghat/\iota(A_{n-1})}{G}{\bar\iota}{\bar\phi}$,
where $\Ker\bar\phi=\bar\iota(A_1)=\iota(A_n)/\iota(A_{n-1})$. Also note that $(\bar\phi)^{-1}(S)=\Shat/\iota(A_{n-1})$.
We just showed that $b=\frac{pde}{p-1}$ is the unique commutator break of $\cbar$, so
$\iota(A_n)/\iota(A_{n-1})$ is contained in $\gamma_{b}(\Shat/\iota(A_{n-1}))$
and is not contained (and hence has trivial intersection with) $\gamma_{b+1}(\Shat/\iota(A_{n-1}))$.
Since $\iota(A_{n-1})$ is the unique maximal subgroup of $\iota(A_n)$, it follows that
$\iota(A_n)\subseteq \gamma_{b}\Shat$ while $\iota(A_n)\not\subseteq \gamma_{b+1}\Shat$,
so by definition the first commutator break of $c$ is equal to $b=\frac{pde}{p-1}$. 

\skv
{\it Step 2:} Now let $b_1=\frac{pde}{p-1}<\ldots<b_k$ be the commutator breaks of $c$.
We want to show that $k=n$.

Since $\iota(A_n)\cap\gamma_{b_i}\Shat\neq  \iota(A_n)\cap\gamma_{b_i+1}\Shat$ for any
$1\leq i\leq k$ and $|A_n|=p^n$, it is clear that $k\leq n$. Suppose now that $k<n$. 
Then for some $i$ we must have 
$|\iota(A_n)\cap\gamma_{b_i}\Shat/  \iota(A_n)\cap\gamma_{b_i+1}\Shat|\geq p^2$. This
would imply that the group $\gamma_{b_i}\Shat/\gamma_{b_i+1}\Shat$ 
contains an element of order $\geq p^2$, namely the image of a generator
of $\iota(A_n)\cap\gamma_{b_i}\Shat$. Since $b_i\geq \frac{pde}{p-1}>\frac{de}{p-1}+1$,
we obtain a contradiction with Proposition~\ref{powercomm}(a).

The above argument also shows that for each $1\leq i\leq n$ we have
\begin{equation}
\label{eq:saturated}
\iota(A_n)\cap\gamma_{b_i}\Shat=\iota(A_{n+1-i}) \mbox{ and } \iota(A_n)\cap\gamma_{b_i+1}\Shat=\iota(A_{n-i}) 
\end{equation}
We will use this observation in Step~3.
\skv

{\it Step 3:} Finally we prove that $b_i=de(i+\frac{1}{p-1})$ for $1\leq i\leq n$
by induction on $i$. We already know that the result holds for $i=1$.

Suppose now that $b_i=de(i+\frac{1}{p-1})$ for some $1\leq i\leq n-1$.
By the first equation in \eqref{eq:saturated} $\iota(A_{n+1-i})\subset \gamma_{b_i}\Shat$
whence $\iota(A_{n-i})=\iota(A_{n+1-i})^p\subset(\gamma_{b_i}\Shat)^p
= \gamma_{b_i+de}\Shat$,
where the last equality holds by Proposition~\ref{powercomm}(a).
On the other hand, applying the second equation in \eqref{eq:saturated} with $i$ replaced by $i+1$, 
we get $\iota(A_{n-i})\not\subset \gamma_{b_{i+1}+1}\Shat$. It follows that $b_{i+1}+1>b_i+de$, 
whence $b_{i+1}\geq b_i+de= (i+1+\frac{1}{p-1})de$.

Let us now prove the opposite inequality $b_{i+1}\leq(i+1+\frac{1}{p-1})de$.
Let $x$ be a generator of $A_{n-i}$; then $x\in \gamma_{b_{i+1}}\Shat$
by \eqref{eq:saturated}. We already know that $b_{i+1}>2de$, so $\gamma_{b_{i+1}}\Shat=(\gamma_{b_{i+1}-de}\Shat)^p$ by Proposition~\ref{powercomm}(a)(b).
Moreover, $\gamma_{b_{i+1}-de}\Shat$ is powerful, so 
$x=y^p$ for some $y\in\gamma_{b_{i+1}-de}\Shat$ by Theorem~\ref{powchar}. 

We claim that $y\in \iota(A_n)$. Indeed, $\phi(y)^p=\phi(y^p)=\phi(x)=1$ since $x\in \iota(A_n)$.
On the other hand, $\phi(y)\in G_{b_{i+1}-de}$, and $G_{b_{i+1}-de}$ is torsion-free by 
Proposition~\ref{sl1basic}(c). Therefore, $\phi(y)=1$, whence $y\in \iota(A_n)$.

Since $x$ is a generator of $\iota(A_{n-i})$ and $x=y^p$, it follows that
$y\not\in \iota(A_{n-i})$, so $y\not\in \iota(A_n)\cap\gamma_{b_{i+1}}\Shat$ by \eqref{eq:saturated}.
Since $b_i$ and $b_{i+1}$ are consecutive commutator breaks,  $\iota(A_n)\cap\gamma_{b_{i}+1}\Shat=\iota(A_n)\cap\gamma_{b_{i+1}}\Shat$,
so $y\not\in\gamma_{b_{i}+1}\Shat$. 
Since we previously established that $y\in\gamma_{b_{i+1}-de}\Shat$, 
it follows that $b_{i}+1>b_{i+1}-de$ whence $b_{i+1}\leq b_i+de=(i+1+\frac{1}{p-1})de$, as desired.
\QED
\skv
The final result of this section is concerned with the kernel of the restriction map
$H^2(G,A)\to H^2(G_n,A)$ where $A=A_s$ for some $s$.

\begin{Proposition}
Let $n\in\dbN$. Let $K$ be the kernel of the restriction map
$H^2(G,A)\to H^2(G_n,A)$, and let $m=\log_p |K|$. Then
$m\leq \max\{0, \frac{2n+1}{de}-\frac{1}{p-1}\}$.
\label{infkernel}
\end{Proposition}
\proofs
Suppose that $m>0$, and let $c$ be an element of $K$ of order $p^m$
(recall that $K$ is cyclic).
Let $\ex{A}{\Ghat}{G}{\iota}{\phi}$ be the 
extension determined by $c$ and $\Ghat_n=\phi^{-1}(G_n)$.
Since $c\in K$, the extension $\ex{A}{\Ghat_n}{G_n}{\iota}{\phi}$
splits, whence $\gamma_2\Ghat_n\cap \iota(A)=\{1\}$.
On the other hand, $\gamma_2\Ghat_n\supseteq \gamma_{2n+2}\Shat$
by Proposition~\ref{powercomm}(c), whence $\gamma_{2n+2}\Shat\,\cap\,\iota(A)=\{1\}$,
and therefore $\comdep(c) \leq 2n+1$. On the other hand,
$\comdep(c)=(m+\frac{1}{p-1})de$ by Proposition~\ref{combreaks},
whence $m\leq \frac{2n+1}{de}-\frac{1}{p-1}$.
\QED
\vskip .5cm

\section{Reduction to Lie algebras}

\bf{Notations. }\rm Recall that $A_s$ denotes the cyclic group of order $p^s$
for $s\in\dbN$. We set $\gra_s=\log(A_s)$; thus, $\gra_s\cong \dbZ/p^s\dbZ$
considered as an abelian Lie algebra. In analogy with the previous section,
we will use the symbol $A$ (resp. $\gra$) to denote $A_s$ 
(resp. $\gra_s$) for some $s\in\dbN$ when the value of $s$ is not important.
\vskip .1cm

Also recall the basic notations from \S~4: $F$ is a $p$-adic field,
$e$ is the ramification index of $F$, $D$ is a finite-dimensional
central division algebra of degree $d$ over $F$, $G=SL_1(D)$, and
for $n\in\dbZ_{\geq 0}$ we let $\grg_n=\sl(\pi^{n}O_D)$.
\skv

The following two results on Lie algebra cohomology will be established in the next section.
In both theorems $\overline F$ denotes the residue field of $F$.

\begin{Theorem}
\label{Thm61}
Let $p\geq 3$ with $(p,d,|\overline F|)\neq (3,2,3)$, and
let $n=de+1$. 
Then the group $H^2(\grg_n,\gra)^{G}$ has exponent $\leq p^{w+4}$,
where as before $p^w$ is the largest power of $p$ dividing $e$.
\end{Theorem}

\begin{Theorem}
\label{Thm62}
Let $p\geq 3$ with $(p,d,|\overline F|)\neq (3,2,3)$, and assume that $n\equiv l\equiv 1\mod d$ and $l>2n$. Let $c\in H^2(\grg_n,\gra)^G$ and let $c_1$ be the image of $c$ in 
$H^2(\grg_{l},\gra)$. Then $ord(c_1)\leq p^{w+1}$. 
Furthermore, if $m\geq l+(w+1)de$, there exists 
$c_2\in H^2(\grg_{l}/\grg_m,\gra)$ such that $ord(c_2)\leq p^{w+1}$ and $c_2$
maps to $c_1$ under the inflation map $H^2(\grg_{l}/\grg_m,\gra)\to H^2(\grg_{l},\gra)$. 
\end{Theorem}

In this section we will deduce parts (a) and (b) of Theorem~\ref{Thmain} from
Theorem~\ref{Thm61} and Theorem~\ref{Thm62}, respectively. 
We start with the less technical proof of part (a).

\begin{Lemma}
Let $p\geq 5$, $n\geq de+1$, and let $H^2(G_n,A)^{\#}$ be the image
of the restriction map $H^2(G,A)\to H^2(G_n,A)$. Then
$H^2(G_n,A)^{\#}\subseteq H^2_{\grG}(G_n,A)^{G}$ where $\grG=\grG_{ppc}$.
\label{LemmaA}
\end{Lemma}
\proofs
Let $c\in H^2(G,A)$, and write $\Ext(c)=(\ex{A}{\Ghat}{G}{}{\phi})$.
Let $c_1\in H^2(G_n,A)$ be the restriction of $c$; then
$\Ext(c_1)=(\ex{A}{\widehat{G_n}}{G_n}{}{\phi})$ where 
$\widehat{G_n}=\phi^{-1}(G_n)$. We need to prove that
$\Ext(c_1)\in \Ext_{\grG}(G_n,A)^G$ which amounts to showing
that $\widehat{G_n}\in Ob(\grG)$ and $\Ext(c_1)$ is $G$-equivariant.

Let $\Shat=\phi^{-1}(S)$.
By Proposition~\ref{powercomm}(b)(c) we have 
$\widehat{G_n}^p= \gamma_{n+de}\Shat$ and
$\gamma_2\widehat{G_n}\subseteq \gamma_{2n}\Shat$. Since 
$n\geq de$, we conclude that $\gamma_2\widehat{G_n}\subseteq \widehat{G_n}^p$,
so $\widehat{G_n}$ is powerful. Since $G_n$ is torsion-free, $\widehat{G_n}$ is 
automatically $p$-central, so $\widehat{G_n}\in Ob(\grG)$.
Finally, $G$-equivariance of $\Ext(c_1)$ is clear:
the desired action of $G$ on $\widehat{G_n}$ is induced by the 
conjugation action of $\Ghat$.
\QED
\vskip .1cm
Now we prove Theorem~\ref{Thmain}(a) whose statement is recalled below.
\begin{thma} Assume that $p\geq 5$. Then $|H^2(G,\dbQ_p/\dbZ_p)|\leq p^{w+6}$.
\end{thma}
\proofs
By Lemma~\ref{infdep} (and the remark after it) we can assume that $p-1$ divides $e$.
By Corollary~\ref{cor:PR} we just need to show that for any $s\in\dbN$
any element of $H^2(G,A_s)$ has order $\leq p^{w+6}$. So let $A=A_s$
for some $s$, $\gra=\log(A)$, and take any $C\in H^2(G,A)$. 

Let $n=de+1$, and let $C_1$ be the image of $C$ under the restriction
map $H^2(G,A)\to H^2(G_n,A)$. The kernel of this map has order $\leq p^2$
by Proposition~\ref{infkernel} since
$\frac{2n+1}{de}-\frac{1}{p-1}=2+\frac{3}{de}-\frac{1}{p-1}<3$.
Hence $ord(C_1)\geq \frac{ord(C)}{p^2}$. 

By Lemma~\ref{LemmaA},
$C_1\in H^2_{\grG}(G_n,A)^G$ where $\grG=\grG_{ppc}$.
Recall that $G_n$ is torsion-free by Proposition~\ref{sl1basic}(c). Hence by Proposition~\ref{exthom1}(b)
$\Ext_{\grG}(G_n,A)^G$ and $\Ext_{\grL}(\grg_n,\gra)^G$ are both groups isomorphic to each other. Therefore, $H^2_{\grG}(G_n,A)^G$ and $H^2_{\grL}(\grg_n,\gra)^G$ are isomorphic groups, so in particular $H^2_{\grG}(G_n,A)^G$ is isomorphic to a subgroup of $H^2(\grg_n,\gra)^G$. By Theorem~\ref{Thm61}, $H^2(\grg_n,\gra)^{G}$ has exponent $\leq p^{w+4}$.
Hence, $ord(C_1)\leq p^{w+4}$, whence $ord(C)\leq p^{w+6}$, as desired.
\QED
\vskip .2cm

Now we turn to Theorem~\ref{Thmain}(b). The idea
of the proof is similar to that of part (a) except that instead of Weigel's
$\log$ functor we shall work with Lazard's $\log$ functor which
will be applied to appropriate congruence quotients of $G$.

\begin{Lemma}
Let $m,n\in\dbN$ be such that $n\leq m\leq (p-1)n$. Let $A=A_s$ and $\gra=\gra_s$
for some $s$. Then there is a natural isomorphism
$$\log_{n,m}: \Ext(G_n/G_m,A)\to \Ext(\grg_n/\grg_m,\gra)$$
Moreover, $\log_{n,m}$ maps $\Ext(G_n/G_m,A)^G$ onto $\Ext(\grg_n/\grg_m,\gra)^G$.
\label{LemmaB}
\end{Lemma}
\begin{proof}
Since $m\leq (p-1)n$, the nilpotency class of the group $G_n/G_m$ is at most $p-2$.
Thus, if $1\to A\to \Hhat\to G_n/G_m\to 1$ is any central extension, then
$\Hhat$ has nilpotency class $\leq p-1$.
It follows that $\Ext(G_n/G_m,A)=\Ext_{\grG}(G_n/G_m,A)$ where $\grG=\grG_{<p}$, and
similarly $\Ext(\grg_n/\grg_m,\gra)=\Ext_{\grL}(\grg_n/\grg_m,\gra)$ where
$\grL=\grL_{<p}$. Thus Lemma~\ref{LemmaB} follows directly from 
Propositions~\ref{exthom1} which establishes an isomorphism
$\log:\Ext(G_n/G_m,A)\to \Ext(\log(G_n/G_m),\gra)$ mapping
$\Ext(G_n/G_m,A)^G$ onto $\Ext(\log(G_n/G_m),\gra)^G$
and Proposition~\ref{prop:LA}(b) which establishes a $G$-equivariant isomorphism
between $\log(G_n/G_m)$ and $\grg_n/\grg_m$.
\end{proof}
\skv

We are now ready to prove Theorem~\ref{Thmain}(b). As with Theorem~\ref{Thmain}(a),
we recall the statement.

\begin{thmb}
Let $p^w$ be the largest power of $p$ dividing $e$. Suppose that
$4w+15\leq p$. Then $|H^2(G,A)|\leq p^{w+1}$.
\end{thmb}

\proofs Again we can assume that $p-1$ divides $e$ by Lemma~\ref{infdep}, and 
by Corollary~\ref{cor:PR} we just need to show
that for any $s$, any element of $H^2(G,A_s)$ has order $\leq p^{w+1}$.

For the rest of the proof we fix $A=A_s$ for some $s$ and set
$H^2(H)=H^2(H,A)$ for any profinite group $H$. 
\skv

Step 1: We claim that there exist $n,l,m\in\dbN$ such that
\begin{itemize}
\item[(i)] $n$ and $l$ satisfy the hypotheses of Theorem~\ref{Thm62}, that is,
$n\equiv l\equiv 1\mod d$, $l>2n$ and $m\geq l+(w+1)de$, and

\item[(ii)] $2l+1<\frac{pde}{p-1}$, $m\leq n(p-1)$ and $m\geq (w+2+\frac{1}{p-1})de$. 
\end{itemize}
Indeed, first take $l\equiv 1\mod d$ such that $\frac{pde}{2(p-1)}-d-\frac{1}{2}
\leq l<\frac{pde}{2(p-1)}-\frac{1}{2}$,
then take $n\equiv 1\mod d$ such that $\frac{l}{2}-d\leq n<\frac{l}{2}$, and set $m=n(p-1)$.
We have 
\begin{multline*}
m\geq \left(\frac{l}{2}-d\right)(p-1)\geq \left(\frac{pde}{4(p-1)}-\frac{3d}{2}-\frac{1}{4}\right)(p-1)
\\
=
de\left(\frac{p}{4}-\frac{3(p-1)}{2e}-\frac{p-1}{4de}\right)\geq 
de\left(w+\frac{15}{4}-\frac{3(p-1)}{2e}-\frac{p-1}{4de}\right).
\end{multline*}
Since $p-1$ divides $e$ and $p\geq 19$, we have $\frac{15}{4}-\frac{3(p-1)}{2e}-\frac{p-1}{4de}
\geq \frac{15}{4}-\frac{3}{2}-\frac{1}{8}=
\frac{17}{8}> 2+\frac{1}{p-1}$. Thus all inequalities in part (ii) hold.
The remaining inequality $m\geq l+(w+1)de$ in part (i) holds as well since
$l+(w+1)de<\frac{pde}{p-1}+(w+1)de=(w+2+\frac{1}{p-1})de$.
\vskip .12cm

Now consider the following commutative
diagram. All vertical arrows are restriction maps,
horizontal arrows without labels are inflation maps, and
the two labeled arrows are maps defined in Lemma~\ref{LemmaB} (with $\Ext$ groups
replaced by the corresponding $H^2$ groups).
\vskip .1cm
\begin{equation}
\begin{CD}
 @.  @. H^2 (G/G_m) @ >>> H^2(G)\\
@. @. @ VVV @ VVV\\
H^2(\grg_n)@ <<< H^2(\grg_n/\grg_m)@ <\log_{n,m}<< H^2(G_n/G_m)@ >>> H^2(G_n)\\
@ VVV @ VVV @ VVV @ VVV\\
H^2(\grg_{l})@ <<< H^2(\grg_{l}/\grg_m)@ <\log_{l,m}<< H^2(G_{l}/G_m)@ >>> H^2(G_{l})\\
\end{CD}
\label{eqCD}
\end{equation}
\vskip .1cm
Step 2: We will argue by contradiction, so assume that $H^2(G)$ contains an element $C$
of order $p^{w+2}$. By Proposition~\ref{combreaks},
$\infdep(C)=(w+2+\frac{1}{p-1})de$, so by our assumptions
$m\geq \infdep(C)$. Thus $C$ is the inflation image of some $C_1\in H^2(G/G_m)$.
Now let $C_2\in H^2(G_n/G_m)$, $C_5\in H^2(G_{l}/G_{m})$, $C_3\in H^2(G_n)$ and $C_4\in H^2(G_l)$
be the images of $C_1$ in the commutative diagram \eqref{eqCD}.
Since $2l+1<\frac{pde}{p-1}$, the map $H^2(G)\to H^2(G_l)$ is injective
by Proposition~\ref{infkernel}.
Therefore, $ord(C_2)\geq ord(C_5)\geq ord(C_4)=ord(C)= p^{w+2}$.
\vskip .2cm
$$
\begin{CD}
 @.  @. C_1 @ >>> C\\
@. @. @ VVV @ VVV\\
c_3@ <<< c_2@ <\log_{n,m}<< C_2@ >>> C_3 \\
@ VVV @ VVV @ VVV @ VVV\\
c_4@ <<< c_5@ <\log_{l,m}<< C_5@ >>> C_4\\
\end{CD}
$$
\vskip .5cm
Step 3: Let $c_2=\log_{n,m}(C_2)\in H^2(\grg_n/\grg_m)$ and 
$c_5=\log_{l,m}(C_5)\in H^2(\grg_l/\grg_m)$ (recall that the $\log$ maps are defined
in Lemma~\ref{LemmaB}). Since both $\log$ maps are isomorphisms, we have
$ord(c_2)\geq ord(c_5)\geq p^{w+2}$. Let $c_3\in H^2(\grg_n)$ and
$c_4\in H^2(\grg_l)$ be the images of $c_2$ in \eqref{eqCD}.

We claim that $c_3$ and $c_4$ are $G$-equivariant, that is, $c_3\in H^2(\grg_n)^G$ and
$c_4\in H^2(\grg_l)^G$. Indeed, $C_2$ and $C_5$
are $G$-equivariant by Lemma~\ref{equiv_subgp} since they are restriction images of an element of $H^2(G/G_m,A)$.
Hence, by Lemma~\ref{LemmaB}, the elements $c_2$ and $c_5$ are also $G$-equivariant.
Finally, $c_3$ and $c_4$ are $G$-equivariant by Lemma~\ref{equiv_subgp} being inflation images of $c_2$ and $c_5$,
respectively.

Step 4: By Theorem~\ref{Thm62}, the image of $H^2(\grg_{n})^G$ 
in $H^2(\grg_{l})$ has exponent $\leq p^{w+1}$, and
every element of $H^2(\grg_{l})^G$ is inflated from some
element of $H^2(\grg_{l}/\grg_m)$ of order $\leq p^{w+1}$.
Since $c_4\in H^2(\grg_{l})^G$ is the restriction image of $c_3\in H^2(\grg_{n})^G$,
it follows that $ord(c_4)\leq p^{w+1}$ and there
exists $c_5'\in H^2(\grg_{l}/\grg_m)$ such that $ord(c_5')\leq p^{w+1}$, 
and $c_5'$ and $c_5$ both inflate to  $c_4$.
\vskip .12cm

Step 5: Let $K$ be the kernel of the inflation map
$H^2(\grg_{l}/\grg_m,\gra)\to H^2(\grg_{l},\gra)$ where $\gra=\log(A)$.
By the inflation-restriction sequence (see, e.g. \cite[7.5.3, p.~333]{Wb}), $K$
is equal to the transgression image of
$H^1(\grg_{m},\gra)^{\grg_{l}}=\Hom(\grg_m/[\grg_m,\grg_{l}], \gra)$.

Since $l\equiv 1\mod d$ we have $[\grg_m,\grg_{l}]=\grg_{m+l}$ (this can be proved similarly to Proposition~\ref{sl1basic}(b)).
Since $l<de$ by assumption, $[\grg_m,\grg_{l}]\supset \grg_{m+de}=p\,\grg_m$. It follows that
$\grg_m/[\grg_m,\grg_{l}]$ and hence also $K$ has exponent $\leq p$.

On the other hand, $K$ contains the element $c_5-c_5'$ 
which has order $\geq p^{w+2}$ since $ord(c_5)\geq p^{w+2}$ by Step~3 while
$ord(c_5')\leq p^{w+1}$ by Step~4. The obtained contradiction finishes the proof.
\QED

\section{Cohomology of Lie algebras}\label{sec:cohomology}
In this section we identify $\gra_n$ with the abelian
Lie algebra $\frac{1}{p^n}\dbZ_p/\dbZ_p$ for $n\in\dbN$. 
As before, $\gra$ will denote $\gra_n$ for some $n$
when the value of $n$ is not important. We also set
$\gra_{\infty}=\cup_{n=1}^{\infty}\gra_n$. Note that $\gra_{\infty}$ can be naturally
identified with $\dbQ_p/\dbZ_p$.

The goal of this section is to prove Theorems~\ref{Thm61} and \ref{Thm62}.
The main part of the proof consists of describing
$\gra_{\infty}$-valued $G$-invariant cocycles of the Lie algebras 
$\grg_N$ for $N\equiv 1 \mod d$. Once this is achieved, both Theorems~\ref{Thm61} 
and \ref{Thm62} follow very easily. 

Throughout the section we will assume that $p\geq 3$, and if $p=3$, we will assume that
$(d,|\overline F|)\neq (2,3)$ where $\overline F$ is the residue field of $F$. We expect that
Theorems~\ref{Thm61} and \ref{Thm62} are also valid for $p=2$ except when $d=2$
or $(d,|\overline F|)= (3,4)$ (apart from these 2 cases, the only step which does not work for $p=2$ is the proof of Theorem~\ref{G-inv4}).
However, we decided not to deal with $p=2$ since this would add extra technicalities
to the proofs while all the parts of Theorem~\ref{Thmain} require $p\geq 5$ anyway.
\vskip .12cm

Let $\grh=\grg_{N}$ for some $N\in\dbN$ (eventually we will assume that
$N\equiv 1\mod d$, but everything before Proposition~\ref{G-inv3} will be valid for
all $N$). Since $\grh$ is a free $\dbZ_p$-module,  we can realize $H^2(\grh,\gra)$ as 
$Z^2(\grh,\gra)/B^2(\grh,\gra)$ where $Z^2(\grh,\gra)$ and $B^2(\grh,\gra)$
are defined as in \S~2.2.
\vskip .1cm

Now let $\calE=(\ex{\gra}{\widehat{\grh}}{\grh}{\iota}{})$ be an element of $\Ext(\grh,\gra)$.
Given a linear section $\psi:\grh\to\widehat{\grh}$, let $Z_{\psi}$ be the $\gra$-valued cocycle
of $\grh$ corresponding to $\psi$. If $\gra$ is identified with $\iota(\gra)$, the formula
for $Z_{\psi}$ given by \eqref{LA:coc} becomes
$$Z_{\psi}(u,v)=\psi([u,v])-[\psi(u),\psi(v)].$$

Suppose that $\calE$ is a $G$-equivariant extension. Can we always choose $\psi$ such that $Z_{\psi}$ is $G$-invariant?
We do not know the answer to this question; however, it is certainly
possible to make $Z_{\psi}$ invariant under the action of the smaller
group $\Delta$ (defined in Proposition~\ref{sl1basic}(a)) which, as we recall here, 
consists of roots of unity in $W^{*}\cap G$ of order prime to $p$.

\begin{Proposition}
\label{G-inv}
Let $\grh$ and $\calE$ be as above.
The section $\psi$ can be chosen in such a way that
the cocycle $Z=Z_{\psi}$ is $\Delta$-invariant, that is,
$Z(u^g,v^g)=Z(u,v)$ for any $u,v\in\grh$ and $g\in\Delta$
(recall that the action of $G$ on $\gra$ is trivial).
\end{Proposition}
\proofs
Let $\psi$ be some section, and define
$\grz:\Delta\to \Hom(\grh,\gra)$ by setting
$\grz(g) (u)=\psi(u)^g - \psi(u^g)$. Define
the left action of $\Delta$ on $\Hom(\grh,\gra)$ by
setting $g*l(u)=l(u^g)$ (where $l\in\Hom(\grh,\gra)$
and $u\in\grh$). Then it is easy to check that $\grz$
is a $\Hom(\grh,\gra)$-valued $1$-cocycle of $\Delta$.
Since the order of $\Delta$ is prime to $p$
and $\Hom(\grh,\gra)$ has $p$-power order,
the cohomology group $H^1(\Delta,\Hom(\grh,\gra))$
is trivial, whence $\grz$ is a coboundary.
Hence, $\grz(g)(u)=l(u)-l(u^g)$ for some $l\in \Hom(\grh,\gra)$,
so 
\begin{equation}
\label{eq:coh}
\psi(u)^g - \psi(u^g)=l(u)-l(u^g).
\end{equation}
Now define $\psi':\grh\to\grhat$ by $\psi'(u)=\psi(u)-l(u)$.
Clearly, $\psi'$ is also a section of $\calE$.
Note that $l(u)^g=l(u)$ for any $u\in \grh$
and $g\in\Delta$, since the action of $\Delta$ on $\gra$ is trivial.
Therefore, using \eqref{eq:coh} for any $u\in \grh$
and $g\in\Delta$
we have
$\psi'(u^g)=\psi(u^g)-l(u^g)=\psi(u)^g-l(u)=(\psi(u)-l(u))^g=\psi'(u)^g$ , and it follows that
$Z_{\psi'}$ is $\Delta$-invariant.
\QED
\skv
Our next goal is to determine all bilinear $\Delta$-invariant
maps from $\grg_N\times\grg_N$ to $\gra_{\infty}$.
Since $\grg_N$ is a finitely generated $\dbZ_p$-module,
the image of such map lies in $\gra_n$ for some $n$.
But first we introduce some additional notations.
\skv
Let $\sig$ and $\pi$ be as in \S~4.
For a subset $U$ of $W$, we set $$\sl(U)=\{a\in U\colon \tr_{W/F}(a)=0\}.$$
Let $W_{ur}$ (resp. $F_{ur}$) be the maximal unramified extension of $\dbQ_p$ in
$W$ (resp. $F$).
Then $F_{ur}=W_{ur}\cap F$, and the restriction map
$\Gal(W/F)\to\Gal(W_{ur}/F_{ur})$ is an isomorphism (recall that 
the extension $W/F$ is unramified by construction).

\skv
Let $O=O_{W_{ur}}$ be the ring of integers of $W_{ur}$. The extension $W/W_{ur}$
is totally ramified of degree $e$ with uniformizer $\tau=\pi^d$, so
$O_W=\bigoplus\limits_{i=0}^{e-1}\tau^i O=\bigoplus\limits_{i=0}^{e-1}\pi^{di} O$.
Combining this with the decomposition $O_D=\bigoplus\limits_{j=0}^{d-1}\pi ^j O_W$
given by \eqref{eq:directsum} we deduce that
$$O_D=\bigoplus\limits_{i=0}^{de-1}\pi^{i} O$$
Further if we set 
\begin{equation}
O_i=
\left\{
\begin{array}{ll}
O&\mbox{ if } d\nmid i\\
\sl(O)& \mbox{ if } d\mid i
\end{array}
\right.
\label{def_O}
\end{equation}
then $\grg_0=\bigoplus\limits_{i=0}^{de-1}\pi^{i} O_i$ and moreover
\begin{equation}
\label{eq:O}
\grg_N=\bigoplus\limits_{i=N}^{N+de-1}\pi^{i} O_i\mbox{ for all }N\in\dbN.
\end{equation}
Finally if we turn $O_i$ into a $\Delta$-module by setting
$$ \alpha^g=\alpha\frc{g}{\sig^i(g)}
\mbox{ for any } \alpha\in O_i \mbox{ and } g\in\Delta,  $$
it is easy to see that the map from $O_i$ to $\grg_0$
given by $\alpha\mapsto \alpha \pi^i$ is a monomorphism of $\Delta$-modules.

\begin{Definition}\rm
Let $C:\grg_N\times\grg_N\to\gra_{\infty}$ be a bilinear map.
Given $i,j\geq n$, define $C_{i,j}:O_i\times O_j\to\gra_{\infty}$
by $C_{i,j}(\alpha,\beta)=C(\alpha \pi^i,\beta \pi^j)$.
\end{Definition}
\begin{Remark}
Equation \eqref{eq:O} implies that a bilinear map $C:\grg_N\times\grg_N\to\gra_{\infty}$ 
is completely determined by the corresponding $C_{i,j}$ for $N\leq i,j\leq N+de-1$, and one
can think of $C_{i,j}$ for $i$ and $j$ in this range as {\it homogeneous components} of $C$.
However, it will be convenient to keep track of the maps $C_{i,j}$ for all $i,j\geq N$.
\end{Remark}
\skv

It is clear that a bilinear map $C:\grg_N\times\grg_N\to\gra_{\infty}$ is $\Delta$-invariant if and only if
each $C_{i,j}$ is $\Delta$-invariant. A complete description
of $\Delta$-invariant maps from $O_i\times O_j$ to $\gra_{\infty}$
is given by the following proposition. In order to state this and all subsequent results,
we use the maps introduced in Claim~\ref{torring} as well as
the following shortcut notations:\footnote{These notations are slightly different from \S~4.2, e.g.
$\grw$ denotes $W_{ur}/O_{W_{ur}}$ and not $W/O_W$.}
$\grw=W_{ur}/O_{W_{ur}}$, $\grf=F_{ur}/O_{F_{ur}}$,
$\Tr=\tr_{\grw/\grf}$ and $\tr=\tr_{\grw/\grq_p}$, where
$\grq_p=\dbQ_p/\dbZ_p$. Recall that $\grq_p$ can be
naturally identified with $\gra_{\infty}$. 

\begin{Proposition}
\label{G-inv2}
Fix $i,j\in\dbN$ and let $E:O_i\times O_j\to\gra_{\infty}$
be bilinear and $\Delta$-invariant.
\begin{itemize}
\item[(a)] If $d\nmid (i+j)$, then $E=0$.
\item[(b)] If $d\mid (i+j)$ and $d\nmid i$ (hence $d\nmid j$ as well),
there exists $\lam\in\grw $
such that $$E(\alpha,\beta)=
\tr(\lam\alpha\sig^i(\beta))
\mbox{ for all }\alpha\in O_i\mbox{ and }\beta\in O_j.$$
\end{itemize}
\end{Proposition}
\proofs
First, without loss of generality we can assume that $d\nmid i$.

Let $k$ be the smallest integer such that the image of $E$ lies in $\gra_k$.
We will prove (a) and (b) simultaneously
by induction on $k$ (with the case $k=0$ being obvious).

Let $\gro_n=O_n/p O_n$, with the induced $\Delta$-module
structure. If $d\nmid n$, we can identify $\gro_n$ with $\overline W$, the residue field of $W$ (which is naturally isomorphic
to the residue field of $W_{ur}$), and if $d\mid n$, we can identify $\gro_n$
and with $\sl(\overline W)=\{x\in \overline W: \tr_{\overline W/\overline F}(x)=0\}$.
Define the map $\Ebar:\gro_i\times\gro_j\to \frac{1}{p}\dbZ_p/\dbZ_p$ by
setting $$\Ebar(\alpha+pO_i,\beta+pO_j)=p^{k-1}E(\alpha,\beta).$$
Clearly, $\Ebar$ is $\Delta$-invariant as well.

\begin{Claim}
\label{PR:Delta_invaraint}
There exists $\mu\in\overline W$
such that $\Ebar(\alpha,\beta)=\tr_{\overline W/\overline F}(\mu\alpha\sig^i(\beta))$
for any $\alpha\in \gro_i$ and $\beta\in \gro_j$. Moreover, $\mu=0$ if $d\nmid (i+j)$.
\end{Claim}

Claim~\ref{PR:Delta_invaraint} is essentially established in \cite{PR2}; however since we cannot
formally deduce it from the results stated in \cite{PR2}, we will provide a proof.

\begin{proof}[Proof of Claim~\ref{PR:Delta_invaraint}]
First, since by our assumptions the pair $(d,|\overline F|)$ is different from $(3,4)$ and $(2,3)$ and $d\nmid i$,
by \cite[Theorem~7]{Ri} $\gro_i$ is a simple $\overline F[\Delta]$-module. This implies that any $\Delta$-invariant
bilinear map $\gro_i\times \gro_j\to \overline F$ is either the zero map or is non-degenerate on the left. In the latter case this map induces a non-trivial $\overline F[\Delta]$-module homomorphism from $\gro_i$ to the dual of $\gro_j$, which
by \cite[1.8]{PR2} is only possible if $d\mid (i+j)$. Thus, $\Ebar=0$ if $d\nmid (i+j)$.

Suppose now that $d\mid (i+j)$ (and hence $d\nmid j$ as well). 
For each $\mu\in \overline W$ define $E_{\mu}:\overline W\times \overline W\to \overline F$ by
$E_{\mu}(\alpha,\beta)=\tr_{\overline W/\overline F}(\mu\alpha\sig^i(\beta))$. It is straightforward
to check that each $E_{\mu}$ is $\Delta$-invariant (see \eqref{eq:Deltainv} below); moreover the map $\mu\mapsto E_{\mu}$ is injective. Thus, the functions $E_{\mu}$ form a subspace of dimension $[\overline W:\overline F]$ in the space of all $\Delta$-invariant 
bilinear forms $\overline W\times \overline W\to \overline F$. Any subspace of a larger dimension must contain a nonzero
degenerate form, which is impossible by the previous paragraph. Hence any $\Delta$-invariant bilinear form must equal $E_{\mu}$ for some $\mu$.
\end{proof}

\skv
We proceed with the proof of Proposition~\ref{G-inv2}.
Let $\mu$ be as in Claim~\ref{PR:Delta_invaraint}, but now think of it as an element of $\frac{1}{p}O/O\subset \grw=W_{ur}/O$.
Choose $\lam \in \grw $ such that $p^{k-1}\lam=\mu$
(if $\mu=0$, set $\lam=0$).
Define $E_1:O_i\times O_j \to \gra_k$ by
$E_1(\alpha,\beta)=\tr(\lam \alpha\sig^i(\beta)).$
We claim that $E_1$ is $\Delta$-invariant. Indeed, if $d\nmid (i+j)$, then
$E_1=0$ and there is nothing to prove. If $d\mid (i+j)$, then for
any $g\in\Delta$, $\alpha\in O_i$ and $\beta\in O_{j}$ we have
\begin{equation}
\label{eq:Deltainv}
\alpha^g \sig^i(\beta)^g=\alpha\frac{g}{\sig^i(g)}\sig^i(\beta)
\sig^i\lpar\frac{g}{\sig^j(g)}\rpar=\alpha\sig^i(\beta).
\end{equation}
Now $E-E_1$ is a bilinear $\Delta$-invariant map, and it follows from our construction
that the image of $E-E_1$ lies in $\gra_{k-1}$. By induction,
$(E-E_1)(\alpha,\beta)=\tr(\lam_1\alpha\sig^i(\beta))$
for some $\lam_1\in\grw $, whence $E$ has the desired form.
\QED
\skv
\begin{Definition}\rm
Let $C:\grg_N\times\grg_N\to\gra_{\infty}$ be a bilinear map.
A pair of integers $(i,j)$, with $i,j\geq N$, will be called \it{regular }\rm for $C$ if there exists
$\lam\in\grw $ such that
\begin{equation}
\label{lamij}
C_{i,j}(\alpha,\beta)=\tr(\lam\alpha\sig^i(\beta))
\mbox{ for all }\alpha\in O_i \mbox{ and }\beta\in O_j.
\end{equation}
The set of regular pairs will be denoted by $I_{reg}(C)$.
\end{Definition}\rm

One may ask if equation (\ref{lamij}) determines $\lam$ uniquely.
By Claim~\ref{torring}(e)(f), the answer is yes, unless $d=2$, $d\mid i$ and $d\mid j$. In the
latter case the set of all $\lam$ satisfying (\ref{lamij}) has the form $\lam_0+\sl(\grw)$ for some $\lam_0$; 
since $p\neq 2$ by our assumptions, the set $\lam_0+\sl(\grw)$ contains exactly one element of $\grf$.
This observation motivates our next definition.

\begin{Definition}\rm
Let $C$ be as above and $(i,j)\in I_{reg}(C)$. Define $\lam_{i,j}(C)\in \grw $
as follows:
\skv
If $d>2$ or $d\nmid i$ or $d\nmid j$, let $\lam_{i,j}(C)$ be the unique $\lam\in\grw$ such that
(\ref{lamij}) holds.
\skv
If $d=2$, $d\mid i$ and $d\mid j$, let $\lam_{i,j}(C)$ be the unique $\lam\in\grf$ such that
(\ref{lamij}) holds.
\end{Definition}

Proposition~\ref{G-inv2} can now be restated as follows:
if $C$ is $\Delta$-invariant, then
$I_{reg}(C)$ contains all pairs $(i,j)$ such that $d\nmid i$ or $d\nmid j$;
moreover $\lam_{i,j}(C)=0$ whenever $d\nmid (i+j)$.
If $C$ is also a cocycle, we can say much more:

\begin{Proposition}
Let $C$ be a $\Delta$-invariant cocycle of $\grg_{N}$. For $(i,j)\in I_{reg}(C)$ set $\lam_{i,j}=\lam_{i,j}(C)$.
\skv
(a) The following relations hold
provided all symbols occurring in them are defined:

\begin{itemize}
\item[{\rm (R1)}] $\lam_{i,j}=-\sig^i(\lam_{j,i})$

\item[{\rm (R2)}] $\lam_{i+j,k}=\lam_{i,j+k}+\sig^i(\lam_{j,i+k})=
\sig^j(\lam_{i,j+k})+\lam_{j,i+k}$
unless $i,j$ and $k$ are all divisible by $d$.

\end{itemize}

(b) Let $i\geq 2N$, $j\geq N$, with $d\mid i$ and $d\mid j$. Then $(i,j)\in I_{reg}(C)$
provided there exist $k,l\geq N $, with
$k+l=i$ and $k\equiv 1\mod d$, such that 
$$\lam_{k,l+j}+\sig^k(\lam_{l,k+j})\in \grf \eqno(***) $$
Moreover, such $k$ and $l$ exist if $N\equiv 1 \mod d$ and $p\nmid d$.

\label{G-inv3}
\end{Proposition}
\proofs

\bf{(a) }\rm Relation (R1) follows immediately from the equality $C(\alpha \pi^i,\beta \pi^j)=-C(\beta \pi^j,\alpha \pi^i)$,
so we will only prove (R2).

First note that if $d\nmid (i+j+k)$, then all expressions
in (R2) vanish by Proposition~\ref{G-inv2}(a). So, from now on
we assume that $d\mid (i+j+k)$.
Applying the equation
$C([u,v],w)+C([v,w],u)+C([w,u],v)=0$
with $u=\alpha\pi^i$, $v=\beta \pi^j$ and $w=\gamma \pi^k$
and simplifying we have
\begin{equation}
\label{abc}
\tr (\mu \alpha\sig^i(\beta)\sig^{-k}(\gamma)+\nu \beta\sig^j(\alpha)\sig^{-k}(\gamma))=0
\mbox{ whenever } \alpha \pi^i,\beta \pi^j,\gamma \pi^k\in \grg_N
\end{equation}
where $\mu=\lam_{i+j,k}+\sig^i(\lam_{j+k,i})+\sig^{-k}(\lam_{k+i,j})$
and $\nu=\lam_{i+j,k}+\sig^{-k}(\lam_{j+k,i})+\sig^{j}(\lam_{k+i,j})$.
It follows from (R1) that $\mu=\lam_{i+j,k}-\lam_{i,j+k}-\sig^i(\lam_{j,i+k})$
and $\nu=\lam_{i+j,k}-\sig^j(\lam_{i,j+k})-\lam_{j,i+k}$, and all we have to show is that
$\mu=\nu=0$. 

Recall that at least one of the numbers $i,j$ and $k$ is not divisible by $d$; in fact, since we assume that $d\mid(i+j+k)$,
at least two of these numbers are not divisible by $d$. It is easy to see that any permutation of $i$, $j$ and $k$
leaves the set $\{\mu,\sigma^{-i}(\mu),\sigma^k(\mu),\nu,\sigma^{-j}(\nu),\sigma^k(\nu)\}$ invariant. Thus, without loss
of generality we can assume that $d\nmid j$ and $d\nmid k$. 

Since $d\nmid k$, we can choose any $\gamma\in O$ in \eqref{abc}, and therefore by Claim~\ref{torring}(e), 
$\mu\alpha\sig^i(\beta)-\nu\beta\sig^j(\alpha)=0$ for all $\alpha\in O_i$ and $\beta\in O_j$.

Assume now that $\nu\neq 0$; it is clear that $\mu\neq 0$ as well. Arguing as in the proof of Claim~\ref{torring},
we deduce that there exist nonzero elements $\mu_{res},\nu_{res}\in \overline{W}$ (the residue field of $W$)
with the following property:
$$\mu_{res}\alpha\sig^i(\beta)-\nu_{res}\beta\sig^j(\alpha)=0$$ where $\beta\in \overline{W}$ is arbitrary  and
$\alpha\in \overline{W}$ is arbitrary if $d\nmid i$ and lies in $\sl(\overline{W})$ if $d\mid i$.

Since now we have an equation in a field, we can rewrite it as
$\frac{\mu_{res}}{\nu_{res}}=\frac{\beta\sig^j(\alpha)}{\alpha\sig^i(\beta)}$.
Thus, the quotient $\frac{\beta\sig^j(\alpha)}{\alpha\sig^i(\beta)}$ is constant for all nonzero $\alpha,\beta$
as above. Clearly, this is impossible if $d\nmid i$ (in which case $\alpha,\beta\in \overline{W}$
are arbitrary). Thus, we can assume that $d\mid i$, in which case the condition simplifies to
$\frac{\sig^j(\alpha)}{\alpha}$ being constant for all nonzero $\alpha\in \sl(\overline{W})$. 
A simple counting argument shows that the latter is possible
only if $d=2$, in which case $\alpha\in \sl(\overline{W})$
if and only if $\frac{\sig^j(\alpha)}{\alpha}=-1$ (since $j$ is odd).

\skv
Thus, we proved that $\mu=\nu=0$ unless $d=2$, $i$ is even and $j$ and $k$ are both odd,
so we proceed with the latter case. The original equation \eqref{abc} in this case immediately implies that
$\mu=-\nu$. But recall that in this case the definition of
$\lam$'s is different. Since $d\mid i$ and $d\mid (j+k)$,
we have $\lam_{j+k,i}\in \grf$ by definition, so $\sig^t(\lam_{j+k,i})=\lam_{j+k,i}$ for all $t$.
Therefore, $\mu-\nu=\sig^{-k}(\lam_{k+i,j})-\sig^{j}(\lam_{k+i,j})=0$
(since $-k\equiv j\mod d$). Thus, we finally showed that $\mu=\nu=0$.

\skv
\bf{(b) }\rm Let $k,l\geq N$ be such that $k+l=i$ and $k\equiv 1\mod d$.
Apply the equation $C([u,v],w)+C([v,w],u)+C([w,u],v)=0$ with $u=\alpha\pi^k$, $v=\pi^l$ and $w=\gamma \pi^j$,
where $\alpha\in O$ and $\gamma\in\sl(O)$. Since $(k,l+j), (l,k+j)\in I_{reg}(C)$
by Proposition~\ref{G-inv2}, after simplifications we get
\begin{equation}
C((\alpha-\sig^{-1}(\alpha))\pi^{i},\gamma \pi^j)=\tr(\gamma(\sig^{-1}(\nu\alpha)-\nu\alpha)),
\label{eqnu}
\end{equation}
where $\nu=\lam_{k,l+j}+\sig^k(\lam_{l,k+j})$.

If we know that $\nu\in \grf$, then $\sig^{-1}(\nu\alpha)-\nu\alpha=\nu(\sig^{-1}(\alpha)-\alpha)$.
So, (\ref{eqnu}) implies that
$C(\beta\pi^i,\gamma\pi^j)=\tr(\nu\beta\gamma)$ for all $\beta,\gamma\in \sl(O)$
since any $\beta\in \sl(O)$ is of the form $\alpha-\sig^{-1}(\alpha)$ for some $\alpha\in O$
by Claim~\ref{torring}(c).
Therefore, by definition $(i,j)\in I_{reg}(C)$.
\skv

Let us now prove the `moreover' part. First, if $N\equiv 1 \mod d$, we can always write
$i=k+l$ as in the first line of the proof of (b) by setting $k=N$ and $l=i-N$. It remains to show that $\nu\in\grf$ if $p\nmid d$.
Replace $\alpha$ by $\alpha+1$ in (\ref{eqnu}). The left-hand side does not change,
and the right-hand side changes by $\tr(\gamma(\sig^{-1}(\nu)-\nu))$.
Thus $\tr(\gamma(\sig^{-1}(\nu)-\nu))=0$ for any $\gamma\in\sl(O)$, whence
$\sig^{-1}(\nu)-\nu\in\grf$ by Claim~\ref{torring}(g). Now $\Tr(\sig^{-1}(\nu)-\nu)=0$; on the other hand,
$\Tr(\mu)=d\mu$ for any $\mu\in\grf$. Since $p\nmid d$, we conclude
that $\sig^{-1}(\nu)-\nu=0$, whence $\nu\in\grf$ by Claim~\ref{torring}(b).
\QED
\skv
\noindent
\bf{New notations. }\rm 

1. For the remainder of the section we fix $f\in\dbN$ and set $N=df+1$.

2. For $m,n\in\dbN$, with $m\leq n$, we set $[m,n]=\{k\in\dbN\colon m\leq k\leq n\}$.

3. For $n\geq 2\fe +1$, let $I_{n,\fe}=[d\fe +1,dn-(d\fe +1)]$.
\skv
4. (taken from \cite{PR2}). Given $\lam\in\grw$
and $i\geq 0$, let $$\lam(i)=\lam+\sig(\lam)+\ldots+\sig^{i-1}(\lam).$$
Note that $\lam(i)+\sig^i(\lam(j))=\lam(i+j)$.

\begin{Proposition}
Let $C$ be a $\Delta$-invariant cocycle of $\gref$ and set
$\lam_{i,j}=\lam_{i,j}(C)$ for $(i,j)\in I_{reg}(C)$.
Let $n\geq 4\fe +2$. The following hold:
\begin{itemize}
\item[(a)] $I_{reg}(C)$ contains  $(i,dn-i)$ for every $i\in I_{n,\fe}$.

\item[(b)] There exists $\kappa_n\in \grw$
such that $\lam_{i,dn-i}=\kappa_n(i)$ for all $i\in I_{n,\fe}$.
\end{itemize}
\label{kappa}
\end{Proposition}

First we make some preparations. Given $i$ such that $(i,dn-i)\in I_{reg}(C)$,
set $$\mu_i=\lam_{i, dn-i}(C).$$
By Proposition~\ref{G-inv2}, $\mu_i$ is defined
whenever $i\in I_{n,\fe}$ and $d\nmid i$.
Moreover, if $d=2$ (in which case $p\nmid d$ by our assumption),
$\mu_i$ is defined for all $i\in I_{n,\fe}$ by Proposition~\ref{G-inv3}(b),
so we already proved Proposition~\ref{kappa}(a) for $d=2$.

Relation (R2) of Proposition~\ref{G-inv3}(a) implies that
\begin{equation}
\mu_{i+j}=\mu_i+\sig^i(\mu_j)=\mu_j+\sig^j(\mu_i)
\mbox { unless } d\mid i \mbox{ and } d\mid j  
\label{mu1}
\end{equation}
(whenever $\mu_i,\mu_j,\mu_{i+j}$ are defined).

\begin{Claim}
\label{aux1}
Let $\mu=\mu_{d\fe +1}$. Then $\mu_{k}-\mu(k)\in \grf$ 
for all $k\in I_{n,\fe}$ with $k\equiv \pm 1\mod d$.
\end{Claim}
\proofs
Let $S=\{k\in I_{n,\fe} : \mu_{k}-\mu(k)\in \grf\}$.
We proceed in several steps.
\skv

\it{Step 1: }\rm $k\in S$ if $k\in[d\fe+1,n-(2d\fe +2)]$ and $k\equiv 1\mod d$.

\noindent
\it{Subproof: }\rm The restrictions on $k$ imply that $k+d\fe +1\in I_{n,\fe}$.
By (\ref{mu1}) we have $\mu_{d\fe +(k+1)}=\mu+\sig(\mu_{k})=\mu_{k}+\sig^k(\mu)$.
(Note that $\mu_{d\fe +(k+1)}$ is defined since either $d=2$ or $d\nmid (d\fe +(k+1))$).
Therefore,
$$\sig(\mu_{k})-\mu_k=\sig^k(\mu)-\mu=\sig(\mu(k))-\mu(k).$$
Hence $\sig(\mu_{k}-\mu(k))=\mu_k-\mu(k)$, and so
$\mu_{k}-\mu(k)\in \grf$ by Claim~\ref{torring}(b).
\skv

\it{Step 2: }\rm $k\in S$ if $k\in [2d\fe +2, n-(d\fe +1)]$ and $k\equiv -1\mod d$.

\noindent
\it{Subproof: }\rm By step 1, $dn-k\in S$, so $\mu_{dn-k}-\mu(dn-k)\in\grf$.
By Proposition~\ref{G-inv3}(a) we have $\mu_{k}=-\sig^k(\mu_{dn-k})$. Thus,
$\mu_k+\sig^k(\mu(dn-k))=-\sig^k(\mu_{dn-k}-\mu(dn-k))\in \grf$.
Since $\sig^k(\mu(dn-k))=\mu(dn)-\mu(k)=n\Tr(\mu)-\mu(k)$, it follows that
$\mu_k-\mu(k)=\mu_k+\sig^k(\mu(dn-k))-n\Tr(\mu)\in\grf$, and so $k\in S$.
\skv 
\it{Step 3: }\rm if $i,j\in I_{n,\fe}$ are such that $i\equiv j\equiv \pm 1 \mod d$ and $i+j\in I_{n,\fe}$,
then $i\in S$ if and only if $j\in S$.

\noindent
\it{Subproof: }\rm Since $i\equiv j\equiv \pm 1\mod d$, we have $\sig^i=\sig^j=\sig^{\pm 1}$, and
(\ref{mu1}) yields $\mu_i-\mu_j=\sig^{\pm 1}(\mu_i-\mu_j)$ whence $\mu_i-\mu_j\in\grf$.
Since $\mu(i)-\mu(j)=\frac{i-j}{d}\Tr(\mu)\in\grf$, the assertion of step 3 is clear.
\skv
\it{Step 4: }\rm $k\in S$ for any $k\in I_{n,\fe}$ with $k\equiv - 1 \mod d$.

\noindent
\it{Subproof: }\rm Let $i=d\fe +(d-1)$ and $j=2d\fe +(d-1)$. By step 2 we have $j\in S$.
Since $i+j=d(3\fe +2)-2\leq dn-d\fe -2$ (recall that $n\geq 4f+2$ by the hypotheses in Proposition~\ref{kappa}), 
step 3 implies that $i=d\fe +(d-1)\in S$.
Once again by step~3, we get that $k\in S$ for any $k\in [d\fe +d-1, 2d\fe +(d-1)]$ with
$k\equiv -1\mod d$. Combining this result with step 2, we conclude that $k\in S$
for any $k\in I_{n,\fe}$ with $k\equiv -1\mod d$.
\skv
\it{Step 5: }\rm $k\in S$ for any $k\in I_{n,\fe}$ with $k\equiv 1 \mod d$.

\noindent
\it{Subproof: }\rm This follows from Step 4 and equality $\mu_{k}=-\sig^k(\mu_{dn-k})$
by the same argument as in Step 2.
\QED
\skv
\skv
\noindent
\it{Proof of Proposition~\ref{kappa}(a). }\rm 
We already know that $(i,dn-i)\in I_{reg}(C)$ if $d\nmid i$.
Since $C$ is a cocycle, $(i,dn-i)\in I_{reg}(C)$
if and only if $(dn-i,i)\in I_{reg}(C)$. Thus it suffices to prove that
$(i,dn-i)\in I_{reg}(C)$ when $i\geq dn/2$ and $d\mid i$.
\skv
Fix such $i$. Note that $i\geq 2d\fe +d$ since $n\geq 4\fe +2$,
so we can choose $k,l>d\fe $ such that $k+l=i$, $l\equiv 1\mod d$
(hence $k\equiv -1\mod d$). We will show that $\mu_k+\sig^k(\mu_l)\in \grf$,
which would imply that $(i,dn-i)\in I_{reg}(C)$
by Proposition~\ref{G-inv3}(b).

By Claim~\ref{aux1}, we have $\mu_k-\mu(k)\in\grf$
and $\mu_l-\mu(l)\in\grf$ where $\mu=\mu_{d\fe +1}$ as before.
Therefore, $\mu_k+\sig^k(\mu_l)-(\mu(k)+\sig^k(\mu(l)))=(\mu_k-\mu(k))+\sig^k(\mu_l-\mu(l))\in\grf$.
Since $\mu(k)+\sig^k(\mu(l))=\mu(k+l)=\mu(i)=\frac{i}{d}\Tr(\mu)\in\grf$,
we deduce that $\mu_k+\sig^k(\mu_l)\in\grf$, as desired.
\skv

Once we proved part (a), we know that
$\mu_i$ is defined for all $i\in I_{n,\fe}$. 
We will now establish one more auxiliary result before proving part (b).

\begin{Claim}
Recall that $\mu=\mu_{d\fe +1}$.
For each $k\in I_{n,\fe}$ define $\nu_k=\mu_{k}-\mu(k)$.
Then
\begin{itemize}
\item[(i)] $\nu_k\in\grf$ for all $k\in I_{n,\fe}$
\item[(ii)] $\nu_{i+j}=\nu_i+\nu_j \mbox{ whenever } d\nmid i \mbox{ or }d\nmid j$.
\end{itemize}
\label{aux2}
\end{Claim}
\begin{proof} We will prove (i) and (ii) simultaneously. First if 
$d\nmid i$ or $d\nmid j$ and $i,j,i+j\in I_{n,\fe}$, \eqref{mu1}
yields 
\begin{equation}
\label{eq:nu}
\nu_{i+j}=\mu_{i+j}-\mu(i+j)=\mu_i+\sig^i(\mu_j)-\mu(i)- \sig^i(\mu(j))
=\nu_i+\sig^i(\nu_j).
\end{equation}

Thus, if we already know that $\nu_i,\nu_j\in \grf$ for some $i$ and $j$ with $d\nmid i$ or $d\nmid j$,
then $\nu_{i+j}\in\grf$ (resp. $\nu_{i-j}\in\grf)$ whenever $i+j\in I_{n,\fe}$ (resp. $i-j\in I_{n,\fe}$).
\skv
By Claim~\ref{aux1} $\nu_i\in \grf$ for all $i\in I_{n,\fe}$ with $i\equiv \pm 1\mod d$. Starting with this fact
and using the observation in the previous paragraph, it is easy to deduce that $\nu_k\in\grf$ for all $k\in I_{n,\fe}$.
This proves (i), and (ii) now follows directly from \eqref{eq:nu}.
\end{proof}

\skv
\noindent
\it{Proof of Proposition~\ref{kappa}(b). }\rm  We will use the notations introduced in Claim~\ref{aux2}.
\skv
\bf{Case 1: }\rm $d\neq 2$. For any $i\in [d\fe +1,dn-2d\fe -3]$ we have
$\nu_{d\fe +(i+2)}=\nu_{i+1}+\nu_{d\fe +1}=\nu_{i}+\nu_{d\fe +2}$, whence
$\nu_{i+1}-\nu_{i}=\nu_{d\fe +2}-\nu_{d\fe +1}$
(if $d=2$ and $i$ is even, both $i$ and $d\fe +2$
are divisible by $d$, so the above equalities may not hold).
Therefore, for each $i\in [d\fe +1,dn-2d\fe -2]$ we have $\nu_{i}=\nu_{d\fe +1}+(i-d\fe -1)\nu$
where $\nu=\nu_{d\fe +2}-\nu_{d\fe +1}$.
\skv
Since $\nu_{2d\fe +2}=2\nu_{d\fe +1}$ by Claim~\ref{aux2}(ii), we get
$\nu_{d\fe +1}+(d\fe +1)\nu=2\nu_{d\fe +1}$ and therefore
$\nu_{d\fe +1}=(d\fe +1)\nu$. So, for $i\in [d\fe +1, dn-2d\fe -2]$
we have $\nu_i=i\nu$, whence $\mu_i=\mu(i)+i\nu=\{\mu+\nu\}(i)$.
The formula $\mu_i=\{\mu+\nu\}(i)$ is easily seen to hold for $dn-2d\fe -2<i\leq dn-d\fe -1$
as well, e.g. by (\ref{mu1}).
\skv
\bf{Case 2: }\rm $d=2$. Let $i\in [2\fe +1,2n-4\fe -3]$.
If $i$ is odd,
$\nu_{i+1}-\nu_{i}=\nu_{2\fe +2}-\nu_{2\fe +1}$ as in case 1.
Similarly,
$\nu_{i+1}-\nu_{i}=\nu_{2\fe +3}-\nu_{2\fe +2}$ if $i$ is even.
Now let $\alpha=\nu_{2\fe +2}-\nu_{2\fe +1}$ and $\beta=\nu_{2\fe +3}-\nu_{2\fe +2}$.
Arguing as above, we conclude that
$$\nu_{2\fe +(2i+1)}=\nu_{2\fe +1}+i(\alpha+\beta)\,\mbox{ and }\,
\nu_{2\fe +2i}=\nu_{2\fe +1}-\beta+i(\alpha+\beta)\mbox{ for }i\in [0, n-3\fe -2].$$
The equation $\nu_{2\fe +1}+\nu_{2\fe +2}=\nu_{4\fe +3}$
yields $\nu_{2\fe +1}=\fe (\alpha+\beta)+\beta$, while
the equation $2\nu_{2\fe +1}=\nu_{4\fe +2}$
yields $\nu_{2\fe +1}=\fe (\alpha+\beta)+\alpha$.
It follows that $\alpha=\beta$ and $\nu_{2\fe +1}=(2\fe +1)\alpha$.
The rest of the proof is the same as in case 1.
\QED
\skv
The assertions of Propositions~\ref{G-inv2} and \ref{kappa} motivate the following definition.
\begin{Definition}\rm
A bilinear map $C:\gref\times\gref\to\gra_{\infty}$ will be called \it{regular }\rm
if there exists a sequence $\{\kappa_n\}_{n=2\fe +1}^{\infty}$ such that
for any $i,j\geq d\fe +1$ we have $(i,j)\in I_{reg}(C)$ and
$$
\lam_{i,j}(C)=
\left\{
\begin{array}{ll}
0&\mbox{ if } d\nmid (i+j)\\
\kappa_{(i+j)/d}(i)&\mbox{ if } d\mid (i+j)
\end{array}
\right.
$$

We will say that $\{\kappa_n\}$
is the \it{defining sequence }\rm of $C$
(each $\kappa_n$ is uniquely determined by $C$, e.g. since $\lam_{i+1,dn-(i+1)}(C)-\lam_{i,dn-i}(C)=\sig^i(\kappa_n)$
for any $i\in [df+1,2df]$).
\end{Definition}

Our next result asserts that $\Delta$-invariant cocycles are not far from
being regular.

\begin{Claim}
\label{claimres}
Let $C$ be a $\Delta$-invariant cocycle of $\gref$. The following hold:
\begin{itemize}
\item[(a)] Let $l=dm+1$ with $m\geq 2f+1$. Then 
the restriction of $C$ to $\grg_l\times \grg_l$ is a regular cocycle.

\item[(b)] Assume that $\fe \leq e$. Then $p^3 C$ is a regular cocycle of $\gref$.
\end{itemize}
\end{Claim}

\proofs
(a) is a direct consequence of Proposition~\ref{G-inv2} and Proposition~\ref{kappa}(b).
\vskip .12cm
(b) Let $D=p^3 C$. Clearly, $D$ is $\Delta$-invariant as well, so 
$\lam_{i,j}(D)=0$ if $d\nmid(i+j)$
by Proposition~\ref{G-inv2}. It remains to show that for
any $n\geq 2 \fe +1$ there exists $\kappa_n\in\grw$ such that
$\lam_{i,dn-i}(C)=\kappa_n(i)$ for all $i\in I_{n,\fe}$.

\skv
By Proposition~\ref{kappa}(b) we already know that such $\kappa_n$ exists for all $n\geq 4\fe+2$.

\skv
Recall that $\tau=\pi^d$ is a uniformizer of $F$. Since the extension $F/F_{ur}$
is totally ramified of degree $e$, $O_F$ is a free $O_{F_{ur}}$-module with basis $1,\tau,\ldots, \tau^{e-1}$.
Moreover, $\frac{p^3}{\tau^{3e}}\in O_F$, so there exists a unique sequence $\{d_k\in O_{F_{ur}}\}_{k=3e}^{4e-1}$ such that
$p^3=\sum_{k=3e}^{4e-1}d_k\tau^k$.
For any $n\geq 2 \fe +1$, $i\in I_{n,\fe}$,
$\alpha\in O_i$, and $\beta\in O_{dn-i}$ we have
\begin{multline*}
D(\alpha\pi^i,\beta\pi^{dn-i})= C(\alpha\pi^i, p^3 \beta\pi^{dn-i})=
\sum_{k=3e}^{4e-1} C(\alpha\pi^i, \beta d_k \pi^{dn-i+dk})=\\
\sum_{k=3e}^{4e-1} \tr(\kappa_{n+k}(i)\alpha\sig^i(\beta)d_k)=
\tr(\kappa'_n(i)\alpha\sig^i(\beta)),
\end{multline*}
where $\kappa'_n=\sum_{k=3e}^{4e-1}d_k \kappa_{n+k}$
(the right-hand side of the last equality is defined since
for $k\geq 3e$ we have $n+k\geq (2\fe +1)+3e\geq 5\fe +1 \geq 4\fe +2$).

Thus, $\lam_{i,dn-i}(C)=\kappa'_n(i)$ for all $n\geq 2\fe+1$ and
$i\in I_{n,\fe}$, as desired.
\QED
\skv

We are now ready to give a full characterization of regular $\Delta$-invariant cocycles.
This characterization involves the coefficients of the minimal polynomial of $\tau$
over $F_{ur}$. Since $\tau$ is a uniformizer of $F$, there exists a unique sequence
$\{c_k\in O_{F_{ur}}\}_{k=0}^{e-1}$ such that 
\begin{equation}
\label{eq:tau}
\tau^e=p\sum\limits_{k=0}^{e-1}c_k \tau^k.
\end{equation}
Moreover, $c_0$ must be a unit of $O_{F_{ur}}$.
\begin{Definition}\rm
A sequence $\{\kappa_n\in \grw \}_{n=2\fe +1}^{\infty}$
will be called \it{compatible }\rm if for any $n\geq 2\fe +1$ we have

(C1) $\kappa_{n+e}=p\sum\limits_{k=0}^{e-1}c_k\kappa_{n+k}$ (where $\{c_k\}$ are given by \eqref{eq:tau}) 

(C2) $n\,\Tr(\kappa_n)=0$.

\end{Definition}

\begin{Theorem}
A sequence $\{\kappa_n\}_{n=2\fe +1}^{\infty}$ is the defining sequence
of some regular $\Delta$-invariant cocycle of $\gref$ if and only
if $\{\kappa_n\}$ is compatible.
\label{G-inv4}
\end{Theorem}
\proofs
Let $C$ be a regular $\Delta$-invariant cocycle and
let $\{\kappa_n\}$ be the defining sequence of $C$.
By relation (R1) of Proposition~\ref{G-inv3} we have
$\kappa_n(i)+\sig^i(\kappa_n(dn-i))=0$, whence $n \Tr(\kappa_n)=0$,
so (C2) holds.

Condition (C1) is a consequence of the identity $C(u,pv)=pC(u,v)$.
Indeed, let $n\geq 2\fe +1$, $i\in I_{n,\fe}$. For any $\alpha\in O_i$ and $\beta\in O_{dn-i}$ we have
$C(\alpha\pi^i,\beta\pi^{d(n+e)-i})=C(\alpha\pi^i,\beta\pi^{dn-i}\tau^e)=
C(\alpha\pi^i,\beta\pi^{dn-i}\cdot p\,\sum\limits_{k=0}^{e-1}c_k \tau^k)=
p\, \sum\limits_{k=0}^{e-1}C(\alpha\pi^i,c_k\beta\pi^{d(n+k)-i}).$
Hence, $\tr(\kappa_{n+e}(i)\alpha\sig^i(\beta))=
p\, \sum\limits_{k=0}^{e-1}\tr(\kappa_{n+k}(i)c_k\alpha\sig^i(\beta))$,
and (C1) follows immediately.
\skv
Conversely, let $\{\kappa_n\}$ be compatible.
Given $i,j\geq d\fe +1$, let $\lam_{i,j}=\kappa_{(i+j)/d}(i)$ if $d\mid (i+j)$
and $\lam_{i,j}=0$ if $d\nmid (i+j)$.
By equation \eqref{eq:O} there exists a unique
bilinear map
$C:\gref\times\gref\to\gra_{\infty}$ such that
$C_{i,j}(\alpha,\beta)=\tr(\lam_{i,j}\alpha\sig^i(\beta))$
for all $i,j\in [d\fe +1,d\fe+de]$, $\alpha\in O_i$ and $\beta\in O_j$.
Condition (C1) ensures that the same formulas for $C_{i,j}$ are valid for all $i,j\geq d\fe+1$.
By direct computation, $C$ is $\Delta$-invariant. 

We also claim that $\lam_{i,j}(C)=\lam_{i,j}$.
This holds automatically except possibly when $d=2, d\mid i$, $d\mid j$. Recall that in the latter case
the definition of $\lam_{i,j}(C)$ is different and we need to check that $\lam_{i,j}\in\grf$.
This holds since by our definition $\lam_{i,dn-i}=\kappa_n(i)=\frac{i}{d}\Tr(\kappa_n)$ if $d\mid i$.
\skv

Condition (C2)
implies that $\{\lam_{i,j}\}$ satisfy relations (R1) and (R2)
of Proposition~\ref{G-inv3}; moreover (R2) holds even if $i,j$ and $k$ are all divisible by $d$.
It is now easy to deduce that $C$ is a cocycle. The equation $C([u,v],w)+C([v,w],u)+C([w,u],v)=0$
follows immediately from (R2) by the computation in \eqref{abc}.

To prove that $C$ is alternating it suffices to show that 
\begin{itemize}
\item[(i)]$C_{j,i}(\beta,\alpha)= -C_{i,j}(\alpha,\beta)$ for all $i,j\geq d\fe +1 $ and
\item[(ii)] $C_{i,i}(\alpha,\alpha)= 0$ for all $i\geq d\fe +1 $
\end{itemize}
Condition (i) follows immediately from (R1) of Proposition~\ref{G-inv3}. If $d\nmid 2i$, then (ii)
is automatic by construction, so assume that $d\mid 2i$. In this case we get
$$C_{i,i}(\alpha,\alpha)=\tr(\lambda_{i,i}\alpha \sigma^i(\alpha))=\tr(\sigma^i(\lambda_{i,i})\sigma^i(\alpha) \sigma^{2i}(\alpha))=
-C_{i,i}(\alpha,\alpha)$$
where the last equality holds since $d\mid 2i$ and $\sigma^i(\lambda_{i,i})=-\lambda_{i,i}$
by (R1). Thus, $2C_{i,i}(\alpha,\alpha)=0$, and since $p>2$, 
we deduce that $C_{i,i}(\alpha,\alpha)=0$, as desired.
\QED
\skv
\skv

Next we show that for any compatible sequence $\{\kappa_n\}$,
there is a better bound on the orders of $\Tr(\kappa_n)$ than
the one given by (C2) alone.
\begin{Lemma}
Let $p^w$ be the highest power of $p$ dividing $e$.
If $\{\kappa_n\}$ is a compatible sequence, then
$p^{w+1}\Tr(\kappa_n)=0$ for all $n$.
\end{Lemma}
\proofs
Let $\mu_n=\Tr(\kappa_n)$ for $n\geq 2\fe +1$. Note that the sequence
$\{\mu_n\}$ is compatible as well.

Let $l$ be the smallest integer such that $p^l \mu_n=0$ for all $n$,
and let $m$ be the largest integer such that $p^{l-1} \mu_m\neq 0$.
Such $l$ and $m$ indeed exist and moreover $m\leq 2\fe +e$ since
$\mu_{n}=p\sum\limits_{k=0}^{e-1}c_k\mu_{n-e+k}$ for $n\geq 2\fe +e+1$.
We know that $m\mu_m=0$, so $p^l$ divides $m$.

Now consider the equality
$\mu_{m+e}=p\sum\limits_{k=0}^{e-1}c_k\mu_{m+k}$.
The element $\sum\limits_{k=0}^{e-1}c_k\mu_{m+k}$
has order $p^{l}$ because
$c_0\mu_m$ has order $p^{l}$ (as $c_0$ is a unit in $O_{F_{ur}}$)
and $c_k \mu_{m+k}$ has order at most $p^{l-1}$ for $k>0$ (by the choice of $m$).
So, $\mu_{m+e}$ has order $p^{l-1}$.

On the other hand, $(m+e)\mu_{m+e}=0$. Since
$p^{l-1}$ divides $m$,
$p^{l-1}$ must divide $e$ as well. Therefore, $l\leq w+1$.
\QED
\skv
\skv

\begin{Proposition}
Let $\grh=\gref$ for some $\fe$, let $C$ be a regular 
$\Delta$-invariant cocycle of $\grh$ and let
$\{\kappa_n\}_{n\geq 2\fe +1}$ be the defining sequence of $C$.
Let $v$ be any integer such that
$C(\grh,\grh)\subseteq \gra_{v}$ or, equivalently,
any integer such that $p^v\kappa_n=0$ for all $n\geq 2\fe +1$.
Then there exists a regular cocycle $C_1$ of $\grh$ such that
\begin{itemize}
\item[(a)] $C$ and $C_1$ represent the same class in $H^2(\grh,\gra_v)$

\item[(b)] If $\{\theta_n\}$ is the defining sequence of $C_1$, then
$p^{w+1}\theta_n=0$ for all $n\geq 2\fe +1$.
\end{itemize}
\label{tracecorr}
\end{Proposition}
\proofs
If $v\leq w+1$, we can simply set $C_1=C$, so we will assume
that $v>w+1$.
\skv
We already know that $p^{w+1}\Tr(\kappa_n)=0$ for all $n\geq 2 \fe +1$.
Hence $\Tr(\kappa_n)\in \frac{1}{p^{w+1}} O_{F_{ur}}/O_{F_{ur}}$,
and therefore by Claim~\ref{torring}(d) there exists $\theta_n\in \grw $
such that $\Tr(\theta_n)=\Tr(\kappa_n)$ and $p^{w+1} \theta_n=0$.

We claim that the sequence $\{\theta_n\}$ can be chosen compatible.
We start by choosing any $\theta_n$ such that $\Tr(\theta_n)=\Tr(\kappa_n)$ and $p^{w+1} \theta_n=0$
for $n\in [2\fe +1, 2\fe +e]$.
Then there exists a unique way
to choose the remaining $\theta_n$ so that (C1) holds, that is, $\theta_{n+e}=p\sum\limits_{k=0}^{e-1}c_k\theta_{n+k}$.

Since $\{\kappa_n\}$ satisfies (C1) as well, it follows that
$\Tr(\theta_n)=\Tr(\kappa_n)$ for all $n$, whence
$\{\theta_n\}$ satisfies (C2). It remains to show that
$p^{w+1}\theta_n=0$ for all $n\geq 2\fe +1$.
The latter is true for $n\in [2\fe +1, 2\fe + e]$ by construction,
and follows from (C1) for $n\geq 2\fe + e+1$.
\skv

We proceed with the proof. Since $\{\theta_n\}$ is compatible, there exists
a regular cocycle $C_1$ whose defining sequence is $\{\theta_n\}$.
Since $p^{w+1} \theta_n=0$ for all $n$, we have
$C_1(\grh,\grh)\subseteq \gra_{w+1}\subseteq \gra_v$.
It remains to prove the following claim:
\begin{Claim}
The cocycles $C$ and $C_1$ represent the same class in $H^2(\grh,\gra_v)$.
\label{cc1}
\end{Claim}
\proofs
Let $B=C-C_1$ and let $\{\kappa'_n\}$ be the defining sequence of $B$.
Then $\kappa'_n=\kappa_n-\theta_n$,
whence $\Tr(\kappa'_n)=0$. Hence by Claim~\ref{torring}(c) there exists
$\{\mu_n\}_{n=2\fe +1}^{\infty}$ such that $\kappa'_n=\mu_n-\sig(\mu_n)$.
Since $p^v \kappa'_n=0$, we can assume that $p^v \mu_n=0$.
Similarly, since ${\kappa'_n}$ satisfies (C1),
we can assume that $\{\mu_n\}$ satisfies (C1).

Let $h:\grh\to \gra_v$ be the unique $\dbZ_p$-linear function such that
$$h(\alpha \pi^n)=
\left\{
\begin{array}{ll}
\tr(\alpha \mu_{n/d})&\mbox{ if } d\mid n\\
0&\mbox{ if } d\nmid n\\
\end{array}
\right.
\mbox{ for } d\fe +1\leq n\leq d\fe+de \mbox{ and } \alpha\in O_n.
$$
Since $\{\mu_n\}$ satisfies (C1), the above formula for $h(\alpha \pi^n)$
holds for all $n\geq d\fe +1$.
\vskip .1cm

We claim that $B(u,v)=h([u,v])$ for any $u,v\in\grh$.
This would imply that $B$ is a coboundary and thus finish
the proof of the claim and Proposition~\ref{tracecorr}.

Let $u=\sum \alpha_i \pi^i$ and $v=\sum \beta_i \pi ^i$
(where $\alpha_i,\beta_i\in O_i$ for all $i$). Then
\begin{multline*}
B(u,v)=\sum_{i,j} B(\alpha_i \pi^i,\beta_j \pi^ j)=
\sum_{n=2\fe +1}^{\infty}\sum_{i+j=dn} \tr(\kappa_n'(i)\alpha_i\sig^i(\beta_j))=\\
\sum_{n=2\fe +1}^{\infty}\sum_{i+j=dn} \tr((\mu_n-\sig^i(\mu_n))\alpha_i\sig^i(\beta_j))=
\sum_{n=2\fe +1}^{\infty}\sum_{i+j=dn} \tr\lpar\mu_n(\alpha_i\sig^i(\beta_j)-\sig^{-i}(\alpha)\beta_j)\rpar\\
=\sum_{n=2\fe +1}^{\infty}\sum_{i+j=dn} h([ \alpha_i \pi^i,\beta_j \pi^ j])=
\sum_{i,j} h([\alpha_i \pi^i,\beta_j \pi^ j])=h([u,v]).
\quad
\qedhere
\end{multline*}

We are now ready to prove Theorems~\ref{Thm62}~and~\ref{Thm61} whose statements are recalled below. 
Recall that $p^w$ is the largest power of $p$ dividing $e$.

\begin{Theorem6.2}
Let $p\geq 3$ with $(p,d,|\overline F|)\neq (3,2,3)$, and assume that
$n\equiv l\equiv 1\mod d$ and $l>2n$. Let $c\in H^2(\grg_n,\gra)^G$ and let $c_1$ be the image of $c$ in 
$H^2(\grg_{l},\gra)$. Then $ord(c_1)\leq p^{w+1}$. 
Furthermore, if $m\geq l+(w+1)de$, there exists 
$c_2\in H^2(\grg_{l}/\grg_m,\gra)$ such that $ord(c_2)\leq p^{w+1}$ and $c_2$
maps to $c_1$ under the inflation map $H^2(\grg_{l}/\grg_m,\gra)\to H^2(\grg_{l},\gra)$. 
\end{Theorem6.2}

\proofs 
Since $c\in H^2(\grg_n,\gra)^G$, by Proposition~\ref{G-inv} it is represented by some $\Delta$-invariant cocycle $C$ of $\grg_n$.
Let $C_1$ be the restriction of $C$ to $\grg_l\times\grg_l$.
By Claim~\ref{claimres}(a), $C_1$ is a regular cocycle of $\grg_l$. Therefore,
by Proposition~\ref{tracecorr}, there exists a cocycle $C_1'$ of $\grg_l$, 
such that $C_1$ and $C_1'$ represent the same cohomology class in $H^2(\grg_{l},\gra)$ and $p^{w+1}C_1'=0$. 
This implies that $c_1=[C_1]=[C_1']$ has order at most $p^{w+1}$. 

Now let $m\geq l+(w+1)de$. Define the $\gra$-valued $2$-cocycle $C_2$ of $\grg_l/\grg_m$
by setting $C_2(u+\grg_m,v+\grg_m)=C_1'(u,v)$. Then $C_2$ is well-defined since
$p^{w+1}C_1'=0$ and $\grg_m\subseteq \grg_{l+(w+1)de}=p^{w+1}\grg_l$.
Let $c_2=[C_2]\in H^2(\grg_l/\grg_m,\gra)$. By construction, $ord(c_2)\leq ord(C_1')\leq p^{w+1}$,
and the inflation image of $c_2$ in $H^2(\grg_l,\gra)$ is equal to $[C_1']=c_1$.
\QED

\begin{Theorem6.1}
Let $p\geq 3$ with $(p,d,|\overline F|)\neq (3,2,3)$, and
let $n=de+1$. 
Then the group $H^2(\grg_n,\gra)^{G}$ has exponent $\leq p^{w+4}$.
\end{Theorem6.1}

\proofs
Let $c\in H^2(\grg_n,\gra)^G$, and let
$C$ be a $\Delta$-invariant cocycle of $\grg_n$ representing $c$.
By Claim~\ref{claimres}(b), $p^3 C$ is a regular cocycle, so arguing as in the proof of Theorem~\ref{Thm62}, 
we conclude that $[p^3 C]\in H^2(\grg_n,\gra)$ has order at most $p^{w+1}$.
Therefore, $[C]$ has order at most $p^{w+4}$.
\QED

\section{Reduction to the small field case}

The purpose of this section is to prove part (c) of
Theorem~\ref{Thmain} whose statement is recalled below.
The author is grateful to Gopal Prasad for suggesting
several ideas used in the proof.

\begin{Theorem}
Assume that $p\geq 19$. Let $F$ be a $p$-adic field containing a primitive
$p^2{\rm th}$ root of unity and such that the extension 
$F/\dbQ_p$ is Galois. Let $D$ be a central division algebra over $F$ whose degree is 
not a power of $p$, and let $G=SL_1(D)$. Then $|H^2(G,\dbR/\dbZ)|\leq p^{w+1}$
where $p^w$ is the largest power of $p$ dividing the ramification index of $F$.
\label{maint2}
\end{Theorem}

\noindent
\bf{Notation:} \rm Throughout this section we set $H^2(G)=H^2(G,\dbR/\dbZ)$ for any group $G$.
\skv

We start with a simple fact about division algebras over local fields.
\begin{Proposition}
\label{divext}
Let $K'/K$ be an extension of local fields, let $n=[K':K]$,
and let $d\in\dbN$ be coprime to $n$.
The following hold:
\begin{itemize}
\item[(a)] Let $D$ be a central 
division algebra over $K$ of degree $d$. 
Then $D\otimes_{K} K'$ is a central division algebra over $K'$ (also of degree $d$).
\item[(b)] Conversely, if $D'$ is a central division algebra over $K'$ of degree $d$,
then $D'\cong D\otimes_{K} K'$ for some central division algebra $D$ over $K$.
\end{itemize}
\end{Proposition}
\proofs
If $F$ is a local field, the Brauer group $Br(F)$ is canonically isomorphic to $\dbQ/\dbZ$.
Under this isomorphism, division algebras of degree $d$ over $F$ correspond to generators of the
subgroup $\frac{1}{d}\dbZ/\dbZ$ of $\dbQ/\dbZ$. The map $E_{K,K'}:Br(K)\to Br(K')$ given
by $D\mapsto D\otimes_{K} K'$ corresponds to multiplication by $n=[K':K]$
under the above identification. Since $n$ is coprime to $d$, $E_{K,K'}$
maps $\frac{1}{d}\dbZ/\dbZ$ onto itself (and sends generators to generators). This yields
both assertions of the proposition.
\QED

\skv
\bf{Cohomology of $SL_d$ over $p$-adic fields. }\rm
\skv

We start by recalling the definition of norm-residue symbols. 
Let $F$ be a $p$-adic field, and let $\mu_F$ be the group of roots
of unity in $F$. Let $F^{ab}$ be the maximal abelian extension of $F$, and let
$\rho_F: F^{*}\to Gal(F^{ab}/F)$ be the (local) Artin map.
\footnote{The map $\rho_F$ itself is frequently called the norm-residue symbol.} 

Now let $n\in\dbN$ be any divisor of $|\mu_F|$. The {\it norm-residue symbol on $F^*$ of order $n$}
is the map $(\cdot,\cdot)_{n,F}:F^*\times F^*\to \mu_F$ given by
\begin{equation}
\label{eq:nrdef}
(a,b)_{n,F}=\frac{\rho_F(a)(\sqrt[n]{b})}{\sqrt[n]{b}}.
\end{equation}
Note that $F$ contains primitive $n^{\rm th}$ root of unity by assumption, so
$\sqrt[n]{b}\in F^{ab}$ and the right-hand side of \eqref{eq:nrdef} is independent of the choice of $\sqrt[n]{b}$.

We will need the following 3 properties of norm-residue symbols. In all statements below $n$ is a divisor
of $|\mu_F|$.
\begin{itemize}
\item[(NR1)] $(a,b)_{n,F}$ is bi-multiplicative, that is, $(ab,c)=(a,c)(b,c)$ and $(a,bc)=(a,b)(a,c)$
for all $a,b,c\in F^*$.
\item[(NR2)] For any $m$ dividing $n$ we have $(a,b)_{m,F}=(a,b)_{n,F}^{n/m}$.
\item[(NR3)] If $K/F$ is a finite extension, then $(a,b)_{n,K}=(N_{K/F}(a),b)_{n,F}$ for all 
$a\in K^{*}$ and $b\in F^{*}$.
\end{itemize}
(NR1) and (NR2) follow immediately from the definition, and (NR3) follows from a basic result in local
class field theory which asserts that $\rho_K(a)_{| F^{ab}}=\rho_F(N_{K/F}(a))$ (see, e.g., \cite[Theorem~6.9]{Iw}).
\skv

To simplify our notations below we will denote $(\cdot,\cdot)_{|\mu_F|,F}$, the norm-residue symbol on $F^{*}$
of maximal possible order, simply by $(\cdot,\cdot)_F$.
\skv

We now turn to the discussion of the second cohomology groups.
Moore~\cite{Mo1} showed that $H^2(SL_d(F))$
is isomorphic to $\mu_F$. Elements of $H^2(SL_d(F))$
can be explicitly described as follows \cite[Theorem~B]{Rp}.
Let $T$ be the diagonal subgroup of $SL_d(F)$. Then
there is a canonical cocycle $c_F:SL_d(F)\times SL_d(F)\to \mu_F$
such that
\begin{itemize}
\item[(i)] the cohomology class $[c_F]$ generates $H^2(SL_d(F))$
\item[(ii)] the restriction of $c_F$ to $T\times T$ is given by
\begin{equation}
c_F(\diag(\lam_1,\ldots,\lam_d),\diag(\mu_1,\ldots,\mu_d))=
\prod_{i\geq j} (\lam_i,\mu_j)_{F}.
\label{MooreRap}
\end{equation}
\end{itemize}
\skv
Now let $\mu_{F, wild}$ be the $p$-primary component of $\mu_F$
and $H^2(SL_d(F))_{wild}$ the corresponding subgroup of $H^2(SL_d(F))$.
If $n=|\mu_F|$ and $q=|\mu_{F,wild}|$, then clearly
$H^2(SL_d(F))_{wild}$ is generated by $[c_F]^{n/q}=[c_F^{n/q}]$. 

Thus, by property (NR2) above,
the restriction of $c_F^{n/q}$ to $T\times T$ is given by
\begin{equation}
c_F^{n/q}(\diag(\lam_1,\ldots,\lam_d),\diag(\mu_1,\ldots,\mu_d))=
\prod_{i\geq j} (\lam_i,\mu_j)_{q,F}.
\label{MooreRap2}
\end{equation}
\begin{Remark}\rm
By \cite[Lemma~3]{Rp}, the restriction map 
$H^2(SL_d(F))\to H^2(T)$ is injective if $d\geq 3$ and has kernel
of order $2$ if $d=2$. Thus, \eqref{MooreRap2} determines the cohomology
class $[c_F^{n/q}]$ uniquely unless $p=d=2$.
\end{Remark}
\skv

The following result is established in \cite[8.2]{PR2}: 
\begin{Lemma}
\label{lem:PR2}
Let $F$ be a $p$-adic field, $D$ a central division algebra over $F$
whose degree is not divisible by $p$. Let $W$ be a maximal unramified
extension of $F$ in $D$, and let $r_{W,D}: H^2(SL_d(W))\to H^2(SL_1(D))$
be the natural restriction map.
Then $|\Im r_{W,D}|=|\mu_{F,wild}|$, and therefore $r_{W,D}$
is injective on $H^2(SL_d(W))_{wild}$. 
\QED
\end{Lemma}

Using this lemma and the above description of cohomology of $SL_n$, 
we can relate the cohomology groups $H^2(SL_1(D))$ and $H^2(SL_1(D'))$
when $D'$ is obtained from $D$ by a field extension.

\begin{Proposition}
Let $K'/K$ be an extension of $p$-adic fields and $l=[K':K]$.
Let $D$ be a central division algebra over $F$ whose degree $d$ is coprime to both $l$ and $p$,
let $D'=D\otimes_{K} K'$, and let $r_{D',D}:H^2(SL_1(D'))\to H^2(SL_1(D))$ be the restriction map.
Let $p^s$ be the largest power of $p$ dividing $l$, and
assume that $|\mu_{K,wild}|\geq p^{s+1}$. Then $|\Ker r_{D',D}|=p^s$.
\label{PrasIAS}
\end{Proposition}
\proofs
First, since $d$ is coprime to $l$, by Proposition~\ref{divext}(a) $D'=D\otimes_{K} K'$
is a central division algebra of degree $d$ over $K'$.
Let $W$ be a maximal unramified extension of $K$ contained in $D$.
Then $[W:K]=d$ and $W'=W\otimes_{K} K'$ is an unramified extension of $K'$ of degree $d$,
so $W'$ is a maximal unramified extension of $K'$ in $D'$. Moreover we have the following
commutative diagram:
$$
\begin{CD}
H^2(SL_d(W'))_{wild}@>r_{W',D'}>>H^2(SL_1(D'))\\
@VV r_{W',W} V @VV r_{D',D} V\\
H^2(SL_d(W))_{wild}@>r_{W,D} >>H^2(SL_1(D))
\end{CD}
$$
We claim that it is sufficient to show that
\begin{equation}
\label{twokernels}
|\Ker r_{W',W}|=p^s.
\end{equation}
Indeed, assume that \eqref{twokernels} holds.
Both maps $r_{W',D'}$ and $r_{W,D}$ are injective by
Lemma~\ref{lem:PR2}. Thus, if $P'=\Im r_{W',D'}$ and $P=\Im r_{W,D}$,
then $$|P'\cap \Ker r_{D',D}|=|\Ker r_{W', W}|=p^s.$$ On the other hand,
$P'$ and $\Ker r_{D',D}$ both lie in $H^2(SL_1(D'))$ which
is cyclic of $p$-power order, so either $\Ker r_{D',D}\subseteq P'$
or $P'\subseteq \Ker r_{D',D}$. 
In the former case, we get $|\Ker r_{D',D}|=p^s$, as desired.
If $P'\subseteq \Ker r_{D',D}$, then $r_{W',W}$ must be the zero map
which contradicts \eqref{twokernels} since $|H^2(SL_d(W'))_{wild}|=|\mu_{W',wild}|\geq |\mu_{K,wild}|\geq p^{s+1}$
by our hypotheses.
\skv

We proceed with proving \eqref{twokernels}. 
Let $n'=|\mu_{W'}|$, $n=|\mu_W|$,
$q'=|\mu_{W',wild}|$ and $q=|\mu_{W,wild}|$.
Let $c'=(c_{W'})^{n'/q'}$ and $c=(c_W)^{n/q}$
where $c_{W'}$ and $c_W$ are as in \eqref{MooreRap}. 
Then $H^2(SL_d(W'))_{wild}$ is generated by $[c']$ 
and $H^2(SL_d(W))_{wild}$ is generated by $[c]$. 
Given $\alpha,\beta\in W^*$,
by properties (NR1)-(NR3) of the norm-residue
symbols we have
$$((\alpha,\beta)_{q',W'})^{q'/q}=(\alpha,\beta)_{q,W'}=(\alpha^l,\beta)_{q,W}=
((\alpha,\beta)_{q,W})^l$$
which by \eqref{MooreRap2} and the remark after it yields $r_{W',W}([c']^{q'/q})=[c]^l$.
The element $[c]^l$ has order $q/p^s>1$.  
Therefore, $\Ker r_{W',W}$ is generated by $([c']^{q'/q})^{q/p^s}=[c']^{q'/p^s}$. 
Since $ord([c'])=q'$, we conclude that $|\Ker r_{W',W}|=p^s$.
\QED
\proofs[Proof of Theorem~\ref{maint2}] Let $d=\deg(D)$. Write $d=d_1d_2$ where $d_1$
is coprime to $p$ and $d_2$ is a power of $p$. By our assumption, $d_1> 1$. 

As before, let $W$ be a maximal unramified extension of $F$ in $D$, and let
$K$ be the unique extension of $F$ of degree $d_2$ inside $W$. Note that
$K$ and $F$ have the same ramification index. Let $D'$ be the centralizer
of $K$ in $D$. By \cite[4.7]{PR2}, $D'$ is a central division algebra
of degree $d_1$ over $K$, and since $d_1>1$, by \cite[4.8]{PR2}
the restriction map $H^2(SL_1(D))\to H^2(SL_1(D'))$
is injective. Thus it is sufficient to show that $|H^2(SL_1(D'))|\leq p^{w+1}$.

Since $F/\dbQ_p$ is Galois and $K/F$ is unramified, the extension $K/\dbQ_p$ is Galois as well
(this holds since any local field has a unique unramified extension of any given degree).
Therefore, there exists an intermediate field $\dbQ_p\subset L\subset K$ such that
$L/\dbQ_p$ is tamely ramified and $K/L$ is wildly ramified (specifically we let
$L$ be the fixed field of a Sylow $p$-subgroup of $Gal(K/K_{ur})$). 
Thus, $|K:L|=p^w$. Since $K$ contains primitive $p {\rm th}$ root of unity, so does $L$.

Since $K/L$ is Galois and $\Gal(K/L)$ is a $p$-group, there is a tower of fields
$L=L_0\subset L_1\subset\ldots \subset L_w=K$ such that $[L_{i+1}:L_i]=p$ for each $i$.
Furthermore, since $|\mu_{K,wild}|\geq p^2$ by hypotheses of the theorem, 
we can assume that $|\mu_{L_1,wild}|= p^2$.
By Proposition~\ref{divext}(b), there exists a central division algebra $D_0$ of degree $d_1$ over $L$ such that 
$D_0\otimes_L K\cong D'$. Let $D_i=D_0\otimes_L L_i$ for $1\leq i\leq w$. We shall prove that
$|H^2(SL_1(D_i)|\leq p^{i+1}$ for $1\leq i\leq w$ by induction on $i$.

The base case $i=1$ follows from Theorem~\ref{Thmain}(b) since $w(L_1)=1$.
Now suppose that $|H^2(SL_1(D_i))|\leq p^{i+1}$ for some $i$. Since $D_{i+1}\cong D_i\otimes_{L_i}L_{i+1}$,
$[L_{i+1}:L_i]=p$ and $|\mu_{L_i,wild}|\geq p^2$, Proposition~\ref{PrasIAS} yields 
$|\Ker \{H^2(SL_1(D_{i+1}))\to H^2(SL_1(D_{i}))\}|\leq p$, whence 
$|H^2(SL_1(D_{i+1}))|\leq p\cdot |H^2(SL_1(D_i))|\leq p\cdot p^{i+1}=p^{i+2}$.
\QED

\appendix

\section{$\exp-\log$ correspondence for powerful $p$-central $\dbZ_p$-Lie algebras and pro-$p$ groups}
\skv
\centerline{\sc by Mikhail Ershov and Thomas Weigel}
\skv

In this appendix we provide additional details on the $\exp$-$\log$
correspondence between the category $\grG_{ppc}$ of finitely generated powerful 
$p$-central pro-$p$ groups and the category $\grL_{ppc}$ of powerful $p$-central 
$\dbZ_p$-Lie algebras of finite rank. Our main goal is to give a proof of Theorem~\ref{Weigel0}
for the pair $(\grL_{ppc},\grG_{ppc})$. The proof given here is based almost entirely
on the results and ideas from \cite{We}, but the order of exposition and some
of our terminology is quite different.
\skv

For the rest of this section we fix a prime $p\geq 5$ and let $\grG=\grG_{ppc}$
and $\grL=\grL_{ppc}$.

\skv
Recall that according to our convention in \S~3, the functors $\exp:\grL\to\grG$
and $\log:\grG\to\grL$ act as identity on the underlying sets, that is,
$\exp(L)=L$ as sets for any $L\in Ob(\grL)$ and $\log(G)=G$ as sets for any $G\in Ob(\grG)$. Such a definition may lead to confusion between group and Lie algebra operations defined on the same sets, most importantly, the group commutator and
the Lie bracket both of which are denoted by $[\cdot,\cdot]$. To avoid such
problems we will slightly change the definitions here and assume
that as a set $\exp(L)$ is the set of formal symbols $\{\exp(x): x\in L\}$,
and similarly for the $\log$ functor. Of course, this only changes the functors up to equivalence. 
\skv
{\bf Notation:} For a set $S$ and $n\in\dbN$ we will denote by $S^{\times n}$ the direct product of $n$ copies of $S$.
We do not use the standard notation $S^n$ here to avoid confusion with the subgroup generated by $n^{\rm th}$ powers.

\subsection{Evaluating Lie power series in powerful $p$-central $\dbZ_p$-Lie algebras}

Let $X=\{x_1,\ldots,x_d\}$ be a finite set. Let $\dbQ_p\la X\ra$
(resp. $\dbQ_p\lla X\rra$) be the associative $\dbQ_p$-algebra of polynomials 
(resp. power series) in non-commuting variables $x_1,\ldots, x_d$. Let $A(X)$ (resp. $\calA(X)$)
be the subalgebras of $\dbQ_p\la x_1,\ldots, x_d\ra$ (resp. $\dbQ_p\lla x_1,\ldots, x_d\rra$)
consisting of polynomials (resp. power series) with zero constant term. 
For each $n\in\dbN$ let $\calA^{(n)}(X)$ be the set of all power series in $\calA(x)$ which have
no terms of degree $<n$. By the {\it degree topology} on $\calA(X)$ we will mean the unique translation-invariant topology
such that $\{\calA^{(n)}(X)\}$ is a base of neighborhoods of $0$.
\vskip .1cm

Let $L(X)$ be the $\dbQ_p$-Lie subalgebra of $A(X)$ generated by $X$, and let $\calL(X)$
be the $\dbQ_p$-Lie subalgebra of $\calA(X)$ consisting of power series all of whose homogeneous
components lie in $L(X)$. Thus, $\calA(X)\cap \dbQ_p\la X\ra=A(X)$, 
$\calL(X)\cap \dbQ_p\la X\ra=L(X)$, and $A(X)$ (resp. $L(X)$) is dense
in $\calA(X)$ (resp. $\calL(X)$) in the degree topology. It is well known that $L(X)$
is a free $\dbQ_p$-Lie algebra on $X$.
\skv

The following terminology is non-standard, but very convenient for our purposes:

\begin{Definition}\rm $\empty$
\begin{itemize}
\item A {\it monic $X$-commutator of weight $1$} is just an element of $X$.
\item For an integer $k\geq 2$, we define {\it monic $X$-commutators of weight $k$} inductively
as formal expressions $[c,d]$ where $c$ and $d$ are monic $X$-commutators of weight $<k$ with $wt(c)+wt(d)=k$.
\item An {\it $X$-commutator of weight $k$} is a formal expression $\lam c$ where $c$ is a monic $X$-commutator of weight $k$
and $\lam\in\dbQ_p$. We define the {\it $p$-adic valuation} $\nu_p(\lam c)$ to be the $p$-adic valuation of $\lam$, that is,
the unique $k\in\dbZ$ such that $\lam\in p^k\dbZ_p\setminus p^{k+1}\dbZ_p$.
\end{itemize}
\end{Definition}

We will usually think of $X$-commutators as elements of $\calL(X)$; however there are two advantages of defining them as formal
expressions. First, we will sometimes need to distinguish between distinct commutators which represent the same element of $\calL(X)$
(e.g. $[[x_1,x_2],x_3]$ and $[x_3,[x_2,x_1]]$); second we will occasionally need to consider $X$-commutators as elements of the algebras of polynomials or power series over other coefficient rings. 

\skv
We will now introduce an important class of $\dbZ_p$-Lie subalgebras of $\calL=\calL(X)$.
Fix a function $\gamma:\dbN\to \dbZ_{\geq 0}$. 

\begin{Definition}\rm $\empty$
\begin{itemize}
\item[(a)] A commutator $w$ of weight $k$ will be called {\it $\gamma$-integral} 
if $\nu_p(w)\geq \gamma(k)$. In other words, $\gamma$-integral commutators of weight $k$ are precisely
elements of the form $\frac{\lam}{p^{\gamma(k)}}c$ where $c$ is a monic $X$-commutator of weight $k$ and $\lam\in\dbZ_p$.
\item[(b)] We define $\calL^{full}_{\gamma}=\calL^{full}_{\gamma}(X)$ to be the $\dbZ_p$-subspace of $\calL(X)$ consisting of all
infinite sums of $\gamma$-integral commutators which converge in the degree topology. In other words, 
$w\in \calL^{full}_{\gamma}$ if and only if $w$ can be written as $\sum\limits_{i=1}^{\infty}w_i$ where each $w_i$ is a $\gamma$-integral
$X$-commutator and $wt(w_i)\to\infty$ as $i\to\infty$.
\item[(c)] An element of $\calL^{full}_{\gamma}$ will be called {\it $\gamma$-integral} if it can be written as 
$\sum\limits_{i=1}^{\infty}\lam_i w_i$ where each $w_i$ is a $\gamma$-integral
$X$-commutator, $wt(w_i)\to\infty$ as $i\to\infty$, each $\lam_i\in\dbZ_p$ and $\lam_i\to 0$ in $\dbZ_p$.
The set of all $\gamma$-integral elements (which is a clearly a $\dbZ_p$-submodule) will be denoted by 
$\calL_{\gamma}=\calL_{\gamma}(X)$.
\end{itemize}
\end{Definition}

In general neither $\calL^{full}_{\gamma}$ nor $\calL_{\gamma}$ is a subalgebra of $\calL(X)$. It is easy to see that
$\calL^{full}_{\gamma}$ is a subalgebra if and only if $\calL_{\gamma}$ is a subalgebra if and only if $\gamma$ is super-additive,
that is, $\gamma(i+j)\geq \gamma(i)+\gamma(j)$ for all $i,j\in\dbN$. We will be interested in a slightly more restrictive condition
on $\gamma$ which guarantees that $\calL^{full}_{\gamma}$ and $\calL_{\gamma}$ are also {\it substitution-closed}.

\begin{Definition}\rm A subset $S$ of $\calL$ will be called {\it substitution-closed} if for any 
$f,f_1,\ldots, f_d\in S$
(recall that $d=|X|$) the element $f(f_1,\ldots, f_d)$ also lies in $S$. Here $f(f_1,\ldots, f_d)$ is defined to be the image
of $f$ under the unique continuous $\dbQ_p$-Lie algebra endomorphism of $\calL$ which sends $x_i$ to $f_i$ for all $i$.
\end{Definition}

The following result can be proved by direct verification.

\begin{Lemma}
\label{subclosed}
The following are equivalent:
\begin{itemize}
\item[(i)] $\gamma(i_1+\ldots+i_k)\geq \gamma(i_1)+\ldots+\gamma(i_k)+\gamma(k)$ for any finite sequence $i_1,\ldots, i_k\in\dbN$
\item[(ii)] $\calL^{full}_{\gamma}$ is a substitution-closed subalgebra of $\calL$
\item[(iii)] $\calL_{\gamma}$ is a substitution-closed subalgebra of $\calL$
\end{itemize}
\end{Lemma}
\skv

Now let $m\in\dbN$ and define the function $\gamma_m:\dbN\to\dbZ_{\geq 0}$ by 
$\gamma_m(k)=\max\{0,k-m\}$. We define
$$\calL_m^{full}(X)=\calL_{\gamma_m}^{full}(X)\mbox{ and }\calL_m(X)=\calL_{\gamma_m}(X)$$
It is straightforward to check that $\gamma_m$ satisfies condition (i) of Lemma~\ref{subclosed}, so
$\calL_m^{full}(X)$ and $\calL_m(X)$ are both substitution-closed subalgebras.

We will refer to elements of $\calL_m (X)$ as {\it $m$-integral}. In particular, if $c$
is a monic $X$-commutator of weight $k$ and $\lam\in\dbZ_p$, then $\lam c$ is $m$-integral if and only
if $k\leq m$ and $\lam\in\dbZ_p$ or $k>m$ and $p^{k-m}\lam\in\dbZ_p$.
\skv
We also let $\calL_{\infty}(X)=\bigcap\limits_{m=1}^{\infty}\calL_m(X)$. It is easy to
see that $\calL_{\infty}(X)=\calL(X)\cap\dbZ_p\lla X\rra$, that is, $\calL_{\infty}(X)$ consists
of all elements $w\in\calL(X)$ such that all coefficients of $w$ (considered as a power series in $X$)
lie in $\dbZ_p$. We will refer to elements of $\calL_{\infty}(X)$ as {\it integral}.
\skv
\noindent
{\bf Topologies on $L_{\gamma}$.} Note that $\calL_{\gamma}$ is never complete with respect to the degree topology. However, $\calL_{\gamma}$ is complete with respect to the {\it $p$-adic topology} where the sets
$\{p^n \calL_{\gamma}\}_{n=1}^{\infty}$ form a base of neighborhoods of $0$. All continuity statements
involving $\calL_{\gamma}$ will be made with respect to the $p$-adic topology unless mentioned otherwise.
It is clear that for any $m\leq n$, the topology on $\calL_n(X)$ induced from the $p$-adic topology
on $\calL_m(X)$ coincides with the $p$-adic topology on $\calL_n(X)$.

\skv
We now turn to the central problem of this subsection. Given an element $w\in \calL(X)$
and a Lie algebra $L\in Ob(\grL)$, we want to define a function $w_L: L^{\times d}\to L$
such that $w_L(u_1,\ldots, u_d)$ can be reasonably interpreted as the value of $w$
under the substitution $x_i\mapsto u_i$ for $1\leq i\leq d$.
\skv

\noindent
{\bf Case 1:} As a warm-up, let us start with the case where $w\in \calL_{\infty}$, that is, $w$ is integral.
In this case we can proceed in an obvious way.
Indeed, since $L$ is powerful, it is clear that for any $d$-tuple $\bfu=(u_1,\ldots, u_d)$,
there exists a unique continuous homomorphism
$\phi_{\bfu}:\calL_{\infty}\to L$ which sends $x_i$ to $u_i$ for each $i$. We define
$w_L(\bfu)=\phi_{\bfu}(w)$.
\skv

\noindent
{\bf Case 2:} Next consider the case where $w$ is a $1$-integral $X$-commutator. In this case we define $w_L$
using certain choices in $L$, but such that the value of $w_L$ modulo $\Omega$ does not depend
on these choices where $\Omega=\{z\in L: pz=0\}$ (recall that $\Omega\subseteq Z(L)$ since $L$ is $p$-central).
Here and throughout the appendix, given an element $v\in pL$, by $\frac{1}{p}v$ we will denote any element
$v'\in L$ satisfying $pv'=v$.

We proceed by induction on the weight of $w$. If $wt(w)=1$, then $w$ is integral, and $w_L$
is already defined (no choices needed here). Now fix $k\geq 2$ and a $1$-integral $X$-commutator $w$ of weight $k$, and assume
we already defined $v_L$ for any $1$-integral $X$-commutator $v$ of weight $k$ such that $v_L$ is choice-independent modulo $\Omega$.

Since $w$ is $1$-integral, by definition $w=\frac{\lam}{p^{k-1}} [y,z]$ where $y$ and $z$ are monic $X$-commutators
of weights $i$ and $j$ with $i+j=k$ and $\lam\in\dbZ_p$. Note that $\frac{1}{p^{i-1}}y$ and
$\frac{1}{p^{j-1}}z$ are $1$-integral, so we can define $w_L(\bfu)=\lam\cdot \frac{1}{p}[\frac{1}{p^{i-1}}y (\bfu),\frac{1}{p^{j-1}}z(\bfu)]$.
Since $L$ is $p$-central and powerful and the values $\frac{1}{p^{i-1}}y (\bfu)$ and $\frac{1}{p^{j-1}}z(\bfu)$ are well-defined modulo $\Omega$,
the element $v=[\frac{1}{p^{i-1}}y (\bfu),\frac{1}{p^{j-1}}z(\bfu)]$ lies in $pL$ and does not depend on any choices. Thus,
$w_L(\bfu)=\lam(\frac{1}{p}v)$ is well defined modulo $\Omega$, as desired.
\skv

\noindent
{\bf Case 3:} Next assume that $w$ is a $2$-integral $X$-commutator. In this case the scalar $\lam$ in the previous paragraph must lie
in $p\dbZ_p$, and therefore the function $w_L$ we just defined does not involve any choices in $L$.

We can now define $w_L: L^{\times d}\to L$ for any $w\in \calL_2(X)$ as follows: write $w=\sum\limits_{i=1}^{\infty} \lam_i w_i$
where $w_i$ is a $2$-integral $X$-commutator and $\lam_i \to 0$ in $\dbZ_p$ and define 
\begin{equation}
\label{comm:function}
w_L(\bfu)=\sum\limits_{i=1}^{\infty} {(\lam_i w_i)}_L(\bfu). 
\end{equation}
The series on the right of \eqref{comm:function} converges since $L$ has finite rank and $\lam_i\to 0$ in $\dbZ_p$. 

The obtained function $w_L$ does not involve any choices in $L$, but it may depend on the expansion of $w$ as a sum $\sum\limits_{i=1}^{\infty} \lam_i w_i$. As we will see shortly (see Corollary~\ref{cor:L3} below), this issue will not arise if we assume that $w\in \calL_3(X)$.
The following is the key result in this subsection:

\begin{Theorem}
\label{LA:pseudofree}
Let $X=\{x_1,\ldots, x_d\}$ for some $d\in\dbN$. Let $L\in Ob(\grL)$ or $L=\calL_3(X)$ and $\bfu=(u_1,\ldots, u_d)$ any $d$-tuple in $L$. 
Then there exists a unique continuous Lie algebra homomorphism $f=f_{\bfu}:\calL_{3}(X)\to L$ such that
$f(c)=c_L(\bfu)$ for every 3-integral commutator $c$. In particular, $f(x_i)=u_i$ for all $1\leq i\leq d$.
\end{Theorem}

The uniqueness of $f$ is clear since by definition $3$-integral commutators form a dense subset of $\calL_3(X)$.
If $L=\calL_3(X)$, the existence is also immediate -- we just start with the substitution homomorphism 
$\widetilde f:\calL(X)\to \calL(X)$ given by $\widetilde f(v(x_1,\ldots, x_d))=v(u_1,\ldots, u_d)$ (recall that $\calL(X)$ is a 
$\dbQ_p$-Lie algebra, so $\widetilde f$ is clearly well-defined) and let $f$ be the restriction of $\widetilde f$ to $\calL_3(X)$. 
Since $\calL_3(X)$ is substitution closed, we have $f(\calL_3(X))\subseteq \calL_3(X)$ .

The proof of the existence part of Theorem~\ref{LA:pseudofree} for $L\in Ob(\grL)$ is fairly long and will be postponed till the last subsection of this Appendix. 

\begin{Corollary}
\label{cor:L3}Let $w\in\calL_3(X)$ and write $w=\sum\limits_{i=1}^{\infty} \lam_i w_i$ where each $w_i$
is a 3-integral commutator and $\lam_i\to 0$ as $i\to \infty$. Then for any $L\in Ob(\grL)$, the function $w_L$ defined by
\eqref{comm:function} (relative to the chosen expansion)
depends only on $w$ and not on the expansion $\sum \lam_i w_i$. 
\end{Corollary}
\begin{proof}
Let $\bfu=(u_1,\ldots, u_d)$ be any $d$-tuple in $L$, and let $f=f_{\bfu}:\calL_{3}\to L$ be the homomorphism
from Theorem~\ref{LA:pseudofree}. Then
$w_L(\bfu)=\sum\limits_{i=1}^{\infty} (\lam_i w_i)_{L}(\bfu)=\sum\limits_{i=1}^{\infty} f(\lam_i w_i)=f(\sum\limits_{i=1}^{\infty}\lam_i w_i)=
f(w)$, so $w_L$ is determined entirely by $w$.
\end{proof}

Note that a continuous homomorphism $f:\calL_3(X)\to L$  need not be uniquely determined by its values on $X$
since the $\dbZ_p$-Lie subalgebra generated by $X$ is not dense in $\calL_3(X)$. Nevertheless we will prove the following: 

\begin{Proposition}
\label{LA:pseudofree2}
Let $L\in Ob(\grL)$ or $L=\calL_3(X)$, let $\bfu=(u_1,\ldots, u_d)$ be a $d$-tuple in $L$, and let $f:\calL_3(X)\to L$ be a continuous Lie algebra
homomorphism such that $f(x_i)=u_i$ for all $i$. Then $f(w)=w_L(\bfu)$ for any $w\in \calL_4(X)$. In particular,
the restriction of $f$ to $\calL_4(X)$ is completely determined by $\bfu$.
\end{Proposition}

\begin{proof}
Let $C$ be the set of elements of the form $\frac{1}{p^{\max\{k-4,0\}}}c$ where $c$ is a monic {\it left-normed} $X$-commutator of length $k$.
Since every monic $X$-commutator is a $\dbZ$-linear combination of monic left-normed commutators of the same weight, the $\dbZ_p$-submodule generated by $C$ is dense in $\calL_4(X)$. Thus, it suffices to prove that $f(z)=f_{\bfu}(z)$ for every $z\in C$ 
(where $f_{\bfu}$ is the homomorphism from Theorem~\ref{LA:pseudofree}).

We argue by induction on $k=wt(c)$. If $k\leq 4$, the result follows directly from the assumption that $f$ is a Lie algebra homomorphism.

Let us assume that $k\geq 5$ and take any $z\in C$ with $wt(z)=k$. Then $z=\frac{1}{p^{k-4}}[y_1,\ldots, y_k]$ where each $y_i=x_{n_i}$ for some $n_i\in \{1,\ldots, d\}$. 

Let
$w=\frac{1}{p^{k-4}}[y_1,\ldots, y_{k-1}]$, so that $z=[w,y_k]$. Clearly, $w\in\calL_3$, so
$$\textstyle
f(z)=[f(w),f(y_k)]=[\frac{1}{p}f(pw),f(y_k)]=[\frac{1}{p}f(pw),f_{\bfu}(y_k)]$$ 
Note that $pw\in C$ and $wt(pw)=k-1$,
so by the induction hypothesis $f(pw)=f_{\bfu}(pw)$. Hence
$f(z)=[\frac{1}{p}f_{\bfu}(pw),f_{\bfu}(y_k)]=[f_{\bfu}(w),f_{\bfu}(y_k)]=f_{\bfu}([w,y_k])=f_{\bfu}(z)$.
\end{proof}

\begin{Corollary}
\label{lem:L2} Let $X=\{x_1,\ldots, x_d\}$ for some $d$. Assume that $L\in Ob(\grL)$ or $L=\calL_3(X)$,
let $M\in Ob(\grL)$, $\phi:L\to M$ be a continuous $\dbZ_p$-Lie algebra homomorphism and $w\in \calL_4(X)$.
Then $\phi$ commutes with $w$, that is, $\phi(w_L(u_1,\ldots,u_d))=w_M(\phi(u_1),\ldots, \phi(u_d))$
for any $d$-tuple $\bfu=(u_1,\ldots, u_d)$ in $L$. 
\end{Corollary}
\begin{proof}
Choose any continuous Lie algebra homomorphism $f:\calL_3(X)\to L$ such that $f(x_i)=u_i$ for $1\leq i\leq d$
(such $f$ exists by Theorem~\ref{LA:pseudofree}). 
Applying Proposition~\ref{LA:pseudofree2} to the maps
$f:\calL_3(X)\to L$ and $\phi\circ f:\calL_3(X)\to M$, we get that 
$\phi(w_L(u_1,\ldots,u_d))=\phi(f(w))=w_M(\phi\circ f(x_1),\ldots, \phi\circ f(x_d))=w_M(\phi(u_1),\ldots, \phi(u_d))$,
as desired.
\end{proof}

\subsection{The $exp$ functor}

Let $X=\{x_1,x_2\}$. Recall from Theorem~\ref{CHF} that the Baker-Campbell-Hausdorff (BCH) series $\Phi\in \calL(X)$
can be written as $\Phi=\sum\limits_c \lam_c c$ where $c$ ranges over all left-normed monic $X$-commutators and
for any such $c$ 
we have $p^{\delta(wt(c))}\lam_c \in\dbZ_p$ where $\delta(k)=\lfloor \frac{k-1}{p-1}\rfloor$ for $k\in\dbN$.
Note that $\delta(k)\leq \gamma_{p-1}(k)=\max\{0,k-(p-1)\}$ and $\lim\limits_{k\to\infty}\gamma(k)-\delta(k)=\infty$.
Hence by definition $\Phi\in \calL_{p-1}(X)$ and in particular $\Phi\in \calL_{4}(X)$ for $p\geq 5$.

Let us now take any $L\in Ob(\grL)$. By Corollary~\ref{cor:L3}, $\Phi$ defines a function from 
$L^{\times 2}$ to $L$ which does not involve any choices in $L$ and does not depend on representation
of $\Phi$ as a convergent sum of $3$-integral commutators. For simplicity of notation this function will also be 
denoted by $\Phi$ rather than $\Phi_L$. 

We already defined the group
$\exp(L)$ in \S~3. Recall (with our new convention) that $\exp(L)=\{\exp(u): u\in L\}$, and 
the group operation
on $\exp(L)$ is defined by
$$\exp(u)\cdot \exp(v)=\exp \Phi(u,v).$$
Define the topology on $\exp(L)$ simply by transferring the topology from $L$ via the map $u\mapsto \exp(u)$. It is straightforward to check that $\exp(L)$ is a topological group.

Also recall that for any morphism $\phi:L\to M$ in $\grL$, the corresponding morphism $\exp(\phi):\exp(L)\to \exp(M)$
is given by $(\exp(\phi))(\exp(u))=\exp(\phi(u))$.
\skv

We need to prove the following:
\begin{Proposition} $\empty$
\label{exp:properties}
\begin{itemize}
\item[(1)] The operation $\cdot$ satisfies the group axioms.
\item[(2)] If $L,M\in Ob(\grL)$ and $\phi:L\to M$ is a homomorphism
of $\dbZ_p$-Lie algebras, then $\exp(\phi):\exp(L)\to\exp(M)$ is a group homomorphism.
\item[(3)] For any $L\in Ob(\grL)$ the group $\exp(L)$ is an object of $\grG$, that is,
$\exp(L)$ is a finitely generated powerful $p$-central pro-$p$ group.
\end{itemize}
\end{Proposition}

\begin{proof}
We start by observing that $\exp(0)$ is clearly the identity element with respect to $\cdot$ and $\exp(u+v)=\exp(u)\exp(v)$ 
whenever $u$ and $v$ commute; in particular, $\exp(nu)=\exp(u)^n$ for all $n\in\dbZ$.
Also since the BCH series $\Phi$ is $4$-integral, (2) follows directly from Corollary~\ref{lem:L2} (applied with $w=\Phi$).
\skv

(1) Let us now prove that $\cdot$ is associative.
First recall that for any $a,b\in \calL (X)$ we have
$$\Exp(a)\cdot \Exp(b) =\Exp( \Phi(a,b))$$ where $\Exp(a)=1+a+\frac{a^2}{2}+\ldots$
and the product on the left-hand side is taken in 
$\calA(X)$. Since multiplication in $\calA(X)$ is associative, for any 
$a,b,c\in \calL(X)$ we have the equality $\Phi(a,\Phi(b,c))=\Phi(\Phi(a,b),c)$. 

Now take any $u_1,u_2,u_3\in L$ and set $X=\{x_1,x_2,x_3\}$. 
By Theorem~\ref{LA:pseudofree} there exists a continuous Lie algebra homomorphism
$f:\calL_3(X)\to L$ such that $f(x_i)=u_i$ for $1\leq i\leq 3$. 
By Corollary~\ref{lem:L2} applied to $w=\Phi$ and $\phi=f$ we have $f(\Phi(a,b))=\Phi(f(a),f(b))$ for all $a,b\in\calL_{3}(X)$. 

Hence 
\begin{multline*}
(u_1\cdot u_2)\cdot u_3=\Phi(\Phi(u_1,u_2),u_3)=\Phi(\Phi(f(x_1),f(x_2)),f(x_3))\\
=\Phi(f(\Phi(x_1,x_2)),f(x_3))=f(\Phi(\Phi(x_1,x_2),x_3)).
\end{multline*}
Similarly $u_1\cdot(u_2\cdot u_3)=f(\Phi(x_1,\Phi(x_2,x_3)))$. Since $\Phi(\Phi(x_1,x_2),x_3)=\Phi(x_1,\Phi(x_2,x_3))$,
we proved that $(u_1\cdot u_2)\cdot u_3=u_1\cdot(u_2\cdot u_3)$, as desired.
\skv

(3) Finally, let us prove that $\exp(L)$ is an object of $\grG$.
First, the sets $\{p^n L\}_{n=1}^{\infty}$ form a base of neighborhoods of $0$ for the topology on
$L$. It is straightforward to check that each $\exp(p^n L)$ is a normal subgroup of
$p$-power index in $G=\exp(L)$. Thus, $G$ (considered with the same topology as $L$)
is a topological group which has a base of neighborhoods consisting of normal subgroups
of $p$-power index, so by definition $G$ is a pro-$p$ group.

Next observe that for any $u,v\in L$ we have
$[\exp(u),\exp(v)]=\exp \Psi(u,v)$ where $\Psi\in\calL(\{x_1,x_2\})$ is the series
defined by $$\Psi=\Phi(-x_1, \Phi(-x_2, \Phi(x_1,x_2))).$$
Since $\Phi$ is $(p-1)$-integral and $\calL_{p-1}$ is substitution-closed, $\Psi$ is also $(p-1)$-integral.
Also, by a simple direct computation $\Psi$ has no linear terms (the leading term of $\Psi$ is $[x_1,x_2]$).
Hence $\Psi(u,v)$ is equal to a convergent sum $\sum\limits_{i=1}^{\infty} h_i$
where each $h_i=[a_i,b_i]$ for some $a_i,b_i\in L$. Since $L$ is powerful,
$[a_i,b_i]=pc_i$ for some $c_i\in L$; moreover, we can assume that the sequence
$\{c_i\}$ is also convergent. Thus, if we let $c=\sum_{i=1}^{\infty}c_i$, then
$\Psi(u,v)=pc$ whence $[\exp(u),\exp(v)]=\exp pc =(\exp c)^p\in G^p$. Thus,
we proved that $G$ is powerful.
\skv
Next we prove that $G$ is a finitely generated. It is well known that for a pro-$p$ group $G$
the minimal number of generators is equal to $\dim_{\dbF_p}G/[G,G]G^p$. We already proved
that $[G,G]G^p=G^p$ and $G^p=\exp(pL)$. Moreover, it is easy to see
that $\exp$ maps additive cosets of $pL$ to multiplicative cosets of $G^p$. Thus,
$|G/[G,G]G^p|=|L/pL|$, and $|L/pL|$ is finite since $L$ is finitely generated as 
$\dbZ_p$-module (in fact, since $L$ is powerful, its minimal number of generators
is equal to $\dim_{\dbF_p}L/pL$, so we actually proved that $G$ and $L$ have the same 
minimal number of generators).
\skv

It remains to prove that $G$ is $p$-central. Suppose that $\exp(u)\in G$ is an element of order $p$.
Since $\exp(u)^p=\exp(pu)$, it follows that $pu=0$. Since $L$ is $p$-central,
$u\in Z(L)$. It is clear from definitions that $\Phi(u,v)=u+v=\Phi(v,u)$ for any $v\in L$,
so $\exp(u)\in Z(G)$. Thus, we proved that $G$ is $p$-central.
\end{proof}

\subsection{The $log$ functor}

The functor $\log:\grG\to\grL$ will be defined quite similarly to $\exp$, although
both its definition and justification of its basic properties involve extra technicalities. 

First we claim that for any integer $1\leq m\leq p-1$,
the set $\calG_{m}(X)=\Exp(\calL_{m}(X))$ is a group (with respect to the usual multiplication on $\calA(X)$). 
Indeed, for any $y,z\in \calL_{m}(X)$ we have
$\Exp(y)\cdot \Exp(z)=\Exp(\Phi(y,z))$ (by Theorem~\ref{CHF})
and $\Phi(y,z)\in \calL_{m}(X)$ since $\Phi$ is 
$(p-1)$-integral (and hence $m$-integral) and $\calL_m(X)$ is substitution-closed.

We define the topology on $\calG_m(X)$ by transferring the topology
from $\calL_m(X)$ via the map $u\mapsto \Exp(u)$. It is easy to check that
the subgroups $\{\calG_m(X)^{p^n}\}_{n\in\dbN}$ form a base of neighborhoods of $1$ for
this topology.
\skv

It is also easy to show that the group $\calG_1(X)$ is powerful -- this can be verified directly
or deduced from the fact that the Lie algebra $\calL_1(X)$ is powerful (the latter is immediate from the definition).
Similarly to the Lie algebra case, the group $\calG_{3}(X)$ will play the role
of a free-like object for the category $\grG$.
\skv

Recall that in the Lie algebra case we showed that every element $w\in\calL_3(X)$ 
naturally defines a function $w_L:L^{\times d}\to L$ (where $d=|X|$) for every $L\in Ob(\grL)$.
In the group case we will take a slightly different approach. Instead of turning elements
of $\calG_3(X)$ into functions, we introduce certain formal expressions $\{H_c\}$ and $\{R_c\}$
(see the definitions below) which can be interpreted as functions $G^{\times d}\to G$ for all
$G\in Ob(\grG)$ as well as for $G=\calG_m(X)$ for certain $m$. This approach will lead
to a more ad hoc, but also more intuitive, definition of the $log$ functor.
\skv

Consider a new set of formal variables $\{t_k\}_{k=1}^{\infty}$.
For each monic $X$-commutator $c$ we define the expression $H_c$ as follows:
\begin{itemize}
\item If $wt(c)=1$, so $c=x_i$ for some $i$, we let $H_c=t_i$.
\item If $wt(c)=k>1$, so $c=[d,e]$ for some monic $X$-commutators $d$ and $e$
with $wt(d)+wt(e)=k$, we define $H_c$ inductively as $H_c=\sqrt[p]{[H_d,H_e]}$.
\end{itemize}
For instance, if $c=[x_1,[x_2,x_3]]$, then $H_c=\sqrt[p]{t_1,\sqrt[p]{[t_2,t_3]}}$.
One should think of $H_c$ as a multiplicative analogue of the $1$-integral commutator $\frac{1}{p^{wt(c)-1}}c$.

\skv
Now assume that $c$ is a monic {\it left-normed} $X$-commutator. Define the expression $\{R_c\}$ as follows:
\begin{itemize}
\item If $wt(c)\leq 4$, so $c=[x_{i_1},\ldots, x_{i_k}]$ for some $1\leq k\leq 4$, let $R_c=[t_{i_1},\ldots, t_{i_k}]$.
\item If $wt(c)=k\geq 5$ and $c=[d,x_j]$ (where $d$ is a monic left-normed $X$-commutator of weight $k-1$), we let 
$R_c=[\sqrt[p]{R_d},t_j]$.
\end{itemize}
For instance, if $c=[x_1,x_2,x_3,x_4,x_5]$, then $R_c=[\sqrt[p]{[t_1,t_2,t_3,t_4]},t_5]$.
One should think of $R_c$ as a multiplicative analogue of the $4$-integral commutator $\frac{1}{p^{\max\{wt(c)-4,0\}}}c$.

The expressions $\{R_c\}$ will be called {\it root-commutators} (this name is equally applicable to the elements $\{H_c\}$,
and our choice is dictated by the fact that we will work with $R_c$ more frequently).
\skv

Similarly to the Lie algebra case, the expressions $H_c$ and $R_c$ naturally define functions $G^{\times d}\to G$
for certain groups $G$. We claim the following:
\begin{itemize}
\item[(i)] If $G\in Ob(\grG)$, each $H_c$ defines a function whose values are well-defined modulo $\Omega=\{g\in G: g^p=1\}$
(recall that $\Omega\subseteq Z(G)$ since $G$ is $p$-central), and each $R_c$ defines a function not involving any choices.
\item[(ii)] If $G=\calG_1(X)$, then $G$ is powerful and torsion-free, so 
each $H_c$ and $R_c$ becomes a well-defined function on $G$. 
Moreover for each $1\leq m\leq p-1$, the subgroup $\calG_m(X)$ of $G$ is invariant under these functions.  
\end{itemize}

To justify (i) we need the following results which will be discussed below in \S~A.5: 
\begin{itemize}
\item[(a)] (see Lemma~\ref{lem:powerful}(b)(i)) Let $G$ be a finitely generated powerful pro-$p$ group. Then every element of
$G^p$ is equal to $g^p$ for some $g\in G$.
\item[(b)] (see Lemma~\ref{lem:pcentral2}) Let $G\in Ob(\grG)$. If $x^p=y^p$ for some $x,y\in G$, then $(x^{-1}y)^p=1$.
\end{itemize}

Property (ii) can be established by a straightforward computation.
We will need a slightly stronger version of (ii) in the special case $m=3$:

\begin{Lemma} 
\label{lem:F4}
Let $c$ be a monic $X$-commutator of weight $\geq 2$ and
$g_1,\ldots, g_m\in \calG_3(X)$. Then $R_c(g_1,\ldots, g_m)\in\calG_4(X)$.
\end{Lemma}
\begin{Remark} Note that $p^{\min\{3,wt(c)-1\}}\cdot \frac{1}{p^{wt(c)-1}}c=\frac{1}{p^{\max\{wt(c)-4,0\}}}c$,
so one should expect the elements 
$R_c (\Exp (x_1),\ldots, \Exp (x_d))$ and $H_c^{p^{\min\{3,wt(c)-1\}}}(\Exp (x_1),\ldots, \Exp (x_d))$ of $\calG_3(X)$
to be related. If $wt(c)\geq 3$, these elements are distinct but have the same leading term.
\end{Remark}

\begin{proof} We argue by induction on $k=wt(c)$.
As before, define the series $\Psi\in\calL(\{x_1,x_2\})$ by $\Psi=\Phi(-x_1, \Phi(-x_2, \Phi(x_1,x_2))).$ 
In the base case $k=2$,
Lemma~\ref{lem:F4} is equivalent to the assertion
\begin{equation}
\label{eq:psi}
\Psi(u,v)\in \calL_4(X) \mbox{ for all }u,v\in \calL_3(X). 
\end{equation}
which can be verified by direct computation.

Let us proceed with the induction step. Suppose that $k=wt(c)>2$, and write 
$c=[d,x_j]$ where $wt(d)=k-1$. Below we will give an argument for $k\geq 5$
(the cases $k=3,4$ are similar and easier). By definition of $R_c$ we have
\begin{equation}
\label{eq:psi2}
R_c(g_1,\ldots, g_m)=[\sqrt[p]{R_d(g_1,\ldots, g_m)},g_j] 
\end{equation}
By assumption
$g_j=\Exp(v_j)$ for some $v_j\in \calL_3(X)$, and by the induction hypothesis
$R_d(g_1,\ldots, g_m)=\Exp(u)$ for some $u\in \calL_4(X)$. Also it is clear that
$u$ has no terms of degree $\leq wt(d)=k-1$. Since $k-1\geq 4$ and $u\in \calL_4(X)$, we have
$\frac{1}{p}u\in \calL_3(X)$ and therefore $\Psi(\frac{1}{p}u,v_j)\in\calL_4(X)$ by \eqref{eq:psi}.

On the other hand, from \eqref{eq:psi2} we get $R_c(g_1,\ldots, g_m)=\Exp (\Psi(\frac{1}{p}u,v_j))$,
and thus $R_c(g_1,\ldots, g_m)\in\calG_4(X)$, as desired.
\end{proof}
\skv

Our next result (Proposition~\ref{BCH:inverse})
asserts that the addition and the Lie bracket on $\calL_{3}(X)$ can be expressed as convergent products of integer powers of the root-commutators
$\{R_c\}$ in $\calG_{3}(X)$. This result, which should be thought of as a formal inversion of the Baker-Campbell-Hausdoff formula, is well known, but we are not aware of a reference where it is stated exactly in this form.
\skv

Let $X=\{x_1,x_2\}$. For each $n\in\dbN$ let $ C_n$ be the set of all left-normed monic $X$-commutators of weight $n$.

\begin{Proposition}
\label{BCH:inverse}
There exist formal expressions $\{a_n\}_{n=2}^{\infty}$ and $\{b_n\}_{n=3}^{\infty}$
in the variables $\{t_1,t_2\}$ and an integer sequence $\{\delta_n\}$ with $\delta_n\to\infty$ 
satisfying the following conditions:
\begin{itemize}
\item[(a)] Each $a_n=\prod\limits_{c\in  C_n} R_c^{\lam_c p^{\delta_n}}$ for some 
$\lam_c\in\dbZ_p$.
\item[(b)] $\Exp(x_1+x_2)=\Exp(x_1)\Exp(x_2) S_A(x_1,x_2)$ where 
$$S_A(x_1,x_2)=\prod\limits_{n=2}^{\infty}a_n (\Exp(x_1),\Exp(x_2))$$
\item[(c)] Each $b_n=\prod\limits_{c\in  C_n} R_c^{\mu_c p^{\delta_n}}$ for some 
$\mu_c\in\dbZ_p$.
\item[(d)] $\Exp([x_1,x_2])=[\Exp(x_1),\Exp(x_2)] S_B(x_1,x_2)$ where 
$$S_B(x_1,x_2)=\prod\limits_{n=3}^{\infty}b_n(\Exp(x_1),\Exp(x_2))$$
\end{itemize}
\end{Proposition}
\begin{Remark} Recall that for each $n\in\dbN$ and $c\in C_n$, the root-commutator $R_c$
defines a function $G^{\times 2}\to G$ for any $G\in Ob(\grG)$ and also for $G=\calG_1(X)$. Hence
each $a_n$ and $b_n$ also naturally defines a function on any of these groups. Also note that since $\delta_n\to\infty$, the products in (b) and (d) converge in the $p$-adic topology on $\calG_1(X)$.
\end{Remark}
\begin{proof} 

We will prove parts (a) and (b) of Proposition~\ref{BCH:inverse} dealing with the addition operation. The proof of (c) and (d) which deal with the Lie bracket is analogous. It will be convenient to prove a slightly stronger and more precise version of Proposition~\ref{BCH:inverse} (which in particular will include a precise formula for $\delta_n$). To state it we need some
additional notations.

Define the function $\gamma:\dbN\to\dbZ_{\geq 0}$ by 
$\gamma(k)=\lfloor \frac{k-1}{p-1}\rfloor$. It is straightforward to check
that $\gamma$ satisfies condition (i) in Lemma~\ref{subclosed}, so $\calL_{\gamma}$ is substitution-closed. Also note that
$\Phi\in\calL_{\gamma}$ by Theorem~\ref{CHF} and therefore $\calG_{\gamma}=\Exp\calL_{\gamma}$ is a group.
For each $n\in\dbN$ define
$$\textstyle \delta_n=\gamma_4(n)-\gamma(n)=\max\{0,n-4\}-\lfloor \frac{n-1}{p-1}\rfloor.$$
Clearly, $\delta_n\geq 0$ and $\delta_n\to \infty$ as $n\to\infty$.

We will now prove that there exists a sequence $\{a_n\}_{n=1}^{\infty}$ with  $a_1=t_1 t_2$ such that
\begin{itemize}
\item[(i)] each $a_n$ satisfies (a) in Proposition~\ref{BCH:inverse} and $a_n(\Exp(x_1),\Exp(x_2))$
lies in $\calG_{\gamma}(X)$
\item[(ii)] for each $n$ the series $\Exp(x_1+x_2)$ and $\prod\limits_{k=1}^n a_k(\Exp(x_1),\Exp(x_2))$
coincide up to (and including) degree $n$ (and hence (b) in Proposition~\ref{BCH:inverse}
holds as well)
\end{itemize}
We will prove (i) and (ii) by simultaneous induction on $n$. In the base case $n=1$
(i) and (ii) obviously hold by definition of $a_1$. 

Now fix $n\geq 2$ and assume that (i) and (ii) hold for all $m<n$.
Then the element $g={\left(\prod\limits_{k=1}^{n-1} a_k(\Exp(x_1),\Exp(x_2))\right)}^{-1}\cdot\Exp(x_1+x_2)$ has no terms of positive degree $<n$ and lies in $\calG_{\gamma}(X)$. Thus, the degree $n$ term of $g$
has the form $\sum\limits_{c\in C_n}\frac{\lam_c}{p^{\gamma(n)}}c$ for some $\lam_c\in\dbZ_p$.

Now define $a_n=\prod\limits_{c\in C_n} R_c^{\lam_c p^{\delta_n}}$.
Then (a) in Proposition~\ref{BCH:inverse} holds (we will verify the second part of condition (i) later). Let $r_c=R_c(\Exp(x_1),\Exp(x_2))$, so
$$r_c=1+\frac{1}{p^{\gamma_4(n)}}c+\mbox{ higher order terms}.$$ Hence the series
${\left(\prod\limits_{k=1}^{n-1} a_k(\Exp(x_1),\Exp(x_2))\right)}^{-1}\cdot\Exp(x_1+x_2)$ and $a_n(\Exp(x_1),\Exp(x_2))$
both have no terms of positive degree $<n$ and the same term
of degree $n$ (since $\delta_n-\gamma_4(n)=-\gamma(n)$)
and hence $\Exp(x_1+x_2)$ and $\prod\limits_{k=1}^n a_k(\Exp(x_1),\Exp(x_2))$
coincide up to (and including) degree $n$, so (ii) holds for $n$.

It remains to show that $a_n(\Exp(x_1),\Exp(x_2))$ lies in $\calG_{\gamma}(X)$. Note that
$r_c=R_c(\Exp(x_1),\Exp(x_2))\in \calG_{4}(X)$ by Lemma~\ref{lem:F4}, so $r_c=\Exp(v_c)$ for some $v_c\in \calL_{4}(X)$.
Then $r_c^{\lam_c p^{\delta_n}}=\Exp(\lam_c p^{\delta_n} v_c)$.
For each $m\in\dbN$ the degree $m$ coefficients of $v_c$ lie in $\frac{1}{p^{\gamma_{4}(m)}}\dbZ_p$
and hence the degree $m$ coefficients of $\lam_c p^{\delta_n} v_c$ lie in 
$p^{\delta_n-\gamma_{4}(m)}\dbZ_p$.

Recall that $v_c$ has no terms of positive degree $<n$. Also 
it is clear from the definitions of $\gamma_4$ and $\gamma$
that for $n\geq m$ we have
$\gamma_4(n)-\gamma_4(m)\geq \gamma(n)-\gamma(m)$ and hence
$\delta_n-\gamma_4(m)=\gamma_{4}(n)-\gamma(n)-\gamma_4(m)\geq -\gamma(m)$.
It follows that $\lam_c p^{\delta_n} v_c\in \calL_{\gamma}(X)$ whence
$r_c^{\lam_c\delta_n}\in \calG_{\gamma}(X)$ for each $c\in C_n$ and
hence $a_n(\Exp(x_1),\Exp(x_2))\in \calG_{\gamma}(X)$ as well.
\
\end{proof}

We can now use Proposition~\ref{BCH:inverse} to define the $\log$ functor.
Let $G\in Ob(\grG)$, that is, $G$ is a finitely generated powerful $p$-central pro-$p$ group.
We define $L=\log(G)$ to be the set of formal symbols $\{\log(g):g\in G\}$.
Given any $g,h\in G$, we set $S_A(g,h)=\prod\limits_{n=2}^{\infty} a_n(g,h)$ and $S_B(g,h)=\prod\limits_{n=3}^{\infty} b_n(g,h)$
where $a_n$ and $b_n$ come from Proposition~\ref{BCH:inverse}, and  
define 
\begin{equation}
\label{def:log}
\log(g)+\log(h)=\log(gh S_A(g,h))
\mbox{ and }
[\log(g),\log(h)]=\log([g,h] S_B(g,h))
\end{equation}

Also, similarly to $\exp$, for any morphism $\phi:G\to H$ in $\grG$, the corresponding morphism $\log(\phi):\log(G)\to\log(H)$
is defined by $(\log(\phi))(\log(g))=\log(\phi(g))$.

\begin{Remark}
The expressions $S_A$ and $S_B$ in Proposition~\ref{BCH:inverse}
are not uniquely determined as $a_n$ and $b_n$ are not uniquely determined; however, our definition of $\log$
does not depend on these choices. The latter is not clear at this point and will be established in the next subsection 
once we prove that the compositions $\exp \circ \log$ and $\log \circ \exp$ are equivalent to
the identity functor.
\end{Remark}
\skv

The following theorem is a (weaker) group-theoretic counterpart of Theorem~\ref{LA:pseudofree}. Its proof is analogous (albeit slightly more technical) and will also be given in the last subsection.

\begin{Theorem}
\label{thm:ppc}
Let $X=\{x_1,\ldots, x_d\}$ and $G\in Ob(\grG)$ or $G=\calG_3(X)$ .
For any finite subset $\{g_1,\ldots,g_d\}$ of $G$ there exists a continuous group homomorphism $f:\calG_{3}(X)\to G$ such that $f(\Exp(x_i))=g_i$ for all $i$. 
\end{Theorem}

\skv

We proceed with the rest of the proof of the Lie algebra axioms for $L=\log(G)$ (where
$G\in Ob(\grG)$ is arbitrary).
Let $\{g_1,\ldots,g_d\}$ be any finite generating set for $G$, let $X=\{x_1,\ldots, x_d\}$,
and let $f:\calG_{3}(X)\to G$ 
be the homomorphism from Theorem~\ref{thm:ppc}.
We claim that $f$ preserves all root-commutators:
\begin{Proposition}
\label{prop:root}
Let $c$ be a left-normed monic $X$-commutator. Then for any $u_1,\ldots,u_d\in  \calG_{3}(X)$
we have
\begin{equation}
\label{eq:hpreserved}
f(R_c(u_1,\ldots, u_d))=R_c(f(u_1),\ldots, f(u_d))
\end{equation}
\end{Proposition}
\begin{proof}
This is proved by induction on $wt(c)$ similarly to Proposition~\ref{LA:pseudofree2}. The full power of Lemma~\ref{lem:F4} (which we have not used so far) is needed to complete the induction step.
\end{proof}
Since the functions $S_A$ and $S_B$ are convergent products of root-commutators,
we conclude from Proposition~\ref{prop:root} that 
$$f(S_A(g,h))=S_A(f(g),f(h))\mbox{ and  }f(S_B(g,h))=S_B(f(g),f(h))
\mbox{ for all }g,h\in \calG_3(X).$$
As in the construction of the $\exp$ functor
(see the last paragraph in the proof of Proposition~\ref{exp:properties}(1)),
we deduce that $L$ is a Lie ring. Moreover, we can turn $L$
into a $\dbZ_p$-Lie algebra by setting $\lam\, \log(g)=\log(g^{\lam})$ for $\lam\in\dbZ_p$
(the expression $g^{\lam}$ is well defined since $G$ is a pro-$p$ group).

To finish proving that $\log:\grL\to\grG$ is a functor, we need to check
the following: 
\begin{itemize}
\item[(i)] Let $G_1,G_2\in Ob(\grG)$ and a let $\phi:G_1\to G_2$ be a group
homomorphism. Then $\log(\phi):\log(G_1)\to \log(G_2)$ is a $\dbZ_p$-Lie algebra homomorphism.
\item[(ii)] For any $G\in Ob(\grL)$, the Lie algebra $L=\log(G)$ is
finitely generated, powerful and $p$-central.
\end{itemize}
These are direct analogues of Proposition~\ref{exp:properties}(2)(3)
and are proved by nearly identical arguments.

\subsection{Why $\exp$ and $\log$ are mutually inverse.}

We have now constructed both functors $\exp:\grL\to\grG$ and $\log:\grG\to\grL$, and 
parts (a)-(f) of Proposition~\ref{Weigel} (which are stated for $\log$) and their analogues for $\exp$
clearly hold by construction. Proposition~\ref{Weigel} also asserts that $\exp$ and $\log$ are mutually inverse,
but of course this is based on the convention that $\exp(L)=L$ and $\log(G)=G$ as sets for any $L\in Ob(\grL)$ and $G\in Ob(\grG)$.
With the different convention adopted in this appendix, the latter statement translates as follows:
\begin{Proposition}$\empty$
\label{prop:final}
\begin{itemize}
\item[(1)] For any $G\in Ob(\grG)$ the map $G\to \exp(\log(G))$ given by $g\mapsto \exp(\log(g))$ is a group homomorphism.
\item[(2)] For any $L\in Ob(\grL)$ the map $L\to \log(\exp(L))$ given by $u\mapsto \log(\exp(u))$ is a  Lie algebra homomorphism.
\end{itemize}
\end{Proposition}

Below we will prove (1). The proof of (2) is similar. 

\begin{proof}[Proof of Proposition~\ref{prop:final}(1)]
Take any $G\in Ob(\grG)$, let $L=\log(G)$ and $G'=\exp(\log(G))$. 
By Theorem~\ref{thm:ppc} there exists a finite set $X$ and a continuous
surjective homomorphism $f:\calG_{3}(X)\to G$. Define $f_*:\calL_{3}(X)\to L$
by 
\begin{align*}
&f_*(u)=\log f(\Exp(u))\mbox{ for } u\in\calL_3(X) \mbox{ or, equivalently, }& \\
&f_*(\Log(y))=\log(f(y))\mbox{ for } y\in\calG_3(X)&.
\end{align*} 
The two expressions for $f_*$ are indeed equivalent since $\Exp$ and $\Log$ are mutually inverse.

Clearly $f_*$ is continuous. We claim that $f_*$ is a Lie algebra homomorphism.
Indeed, take any $u,v\in \calL_{3}(X)$, and let $g=f(\Exp(u))$ and $h=f(\Exp(v))$ (both $g$ and $h$ lie in $G$).
By Proposition~\ref{BCH:inverse}(b),
\begin{multline*}
f_*(u+v)=\log f(\Exp(u+v))=\log f(\Exp(u)\cdot \Exp(v)\cdot S_A(\Exp(u),\Exp(v))) \\
=\log \left( f(\Exp(u))\cdot f(\Exp(v))\cdot f(S_A(\Exp(u),\Exp(v))\right).
\end{multline*}
By Proposition~\ref{prop:root}, $f$ preserves all root commutators and hence commutes with 
$S_A$, so the last expression is equal to $$\log(f(\Exp(u))\cdot f(\Exp(v))\cdot S_A(f(\Exp(u)),f(\Exp(v))))=
\log(g\cdot h\cdot S_A(g,h)).$$ By the definition of addition in $L=\log(G)$ (see \eqref{def:log}), the obtained expression
is equal to $\log(g)+\log(h)=f_*(u)+f_*(v)$, as desired. Similarly one checks that 
$f_*$ preserves the Lie bracket.
\vskip .1cm

Now define $f': \calG_{3}(X)\to G'$ by
\begin{equation}
\label{eq:fprime}
f'(y)=\exp f_*(\Log(y))=\exp\log f(y). 
\end{equation}
We will prove that $f'$ is a group homomorphism (which is obviously continuous) 
arguing similarly to the previous paragraph.
Take any $y_1,y_2\in \calG_{3}(X)$, and let $u_i=\Log(y_i)\in \calL_3(X)$ for $i=1,2$.
Note that $y_1 y_2=\Exp(u_1)\Exp(u_2)=\Exp\, \Phi(u_1,u_2)$, so $\Log(y_1 y_2)=\Phi(u_1,u_2)$.

 Since $f_*:\calL_{3}(X)\to L$ is a continuous Lie algebra homomorphism, it preserves all the
$4$-integral commutators by Corollary~\ref{lem:L2}. In particular, it commutes with $\Phi$,
so $f_*(\Phi(u_1,u_2))=\Phi(f_*(u_1),f_*(u_2))$. Putting everything together we get
\begin{multline*}
f'(y_1)f'(y_2)=\exp f_*(\Log(y_1))\exp f_*(\Log(y_2))=\exp f_*(u_1) \exp f_*(u_2)\\
=\exp \Phi(f_*(u_1),f_*(u_2))=\exp f_*(\Phi(u_1,u_2))=\exp f_*(\Log(y_1 y_2))=
f'(y_1 y_2).
\end{multline*}

Denote the map $g\mapsto \exp(\log(g))$ from $G$ to $G'$ by $\iota$. Then $f'=\iota\circ f$.
We just proved that $f'$ is a homomorphism. Since $f$ is a surjective homomorphism, it follows that
$\iota$ is a homomorphism, as desired.
\end{proof}

\subsection{$(\lam,p)$-groups and $p$-complete Lie algebras}

In this subsection we introduce the notions of $(\lam,p)$-groups
(called $\lam$-groups in \cite{We}) and their Lie-theoretic counterparts
called {\it $p$-complete Lie algebras} and state some results about them.
These results will be used to prove Theorem~\ref{LA:pseudofree}~and~\ref{thm:ppc}.

\begin{Definition}\rm Let $p$ be a fixed prime and let $G$ be a group.
\begin{itemize}
\item[(a)] The {\it lower $p$-series $\{\lam_n(G)\}_{n\geq 1}$} of $G$ is defined by $\lam_1(G)=G$ and 
$\lam_n(G)=\lam_{n-1}(G)^p [\lam_{n-1}(G),G]$ for $n\geq 2$ (note that the quotients $\lam_{n-1}(G)/\lam_n(G)$ are $\dbF_p$-vector spaces).
\item[(b)] The {\it $(\lam,p)$-topology} on $G$ is the unique
translation-invariant topology on $G$ in which $\{\lam_n(G)\}$ is a base of
neighborhoods of identity (since $\lam_n(G)$ are normal, this turns $G$ into
a topological group). The completion of $G$ with respect to its $(\lam,p)$-topology
will be called the {\it $(\lam,p)$-completion}.
\item[(c)] $G$ is called a {\it $(\lam,p)$-group} if $G$ is a complete Hausdorff group with respect to its $(\lam,p)$-topology.
\end{itemize}
\end{Definition}

\begin{Definition}\rm Let $L$ be a topological $\dbZ_p$-Lie algebra. We will say that $L$ is
{\it $p$-complete} if its additive group $(L,+)$ is a $(\lam,p)$-group. 
\end{Definition}

It is easy to see that any $(\lam,p)$-group $G$ is a $\dbZ_p$-powered group (that is, for any $\mu\in\dbZ_p$ there exists a group homomorphism $g\mapsto g^{\mu}$ from $G$ to $G$ such that
$g^{\mu}\cdot g^{\nu}=g^{\mu+\nu}$, $(g^{\mu})^{\nu}=g^{\mu \nu}$ and
$g^{\mu}=\lim g^{\mu_n}$ for any integer sequence $\{\mu_n\}$ such that
$\mu_n\to\mu$ in $\dbZ_p$).

It is not hard to show that finitely generated $(\lam,p)$-groups are exactly finitely generated pro-$p$ groups, but in the 
infinitely generated case neither of the two classes contains the other. For instance, the group $\calG_{1}(X)$ considered 
in this appendix and the additive group of $\calL_{1}(X)$ are $(\lam,p)$-groups but not pro-$p$ groups while 
$\dbZ_p^{\infty}$ with product topology is a pro-$p$ group but not a $(\lam,p)$-group.

\skv

In \S~3 we defined the property of being {\it powerful} only for pro-$p$ groups. The definition actually makes
sense for an arbitrary topological (possibly discrete) group, but we should keep in mind that the definition depends on $p$.

\begin{Lemma}[see Proposition~VI.1.3 in \cite{We}]
\label{lem:powerful} Let $G$ be a powerful $(\lam,p)$-group (with respect to a prime $p$). The following hold for all $n\in\dbN$:
\begin{itemize}
\item[(a)] $\lam_n(G)=G^{p^{n-1}}$.
\item[(b)] 
\begin{itemize}
\item[(i)] Every element of $G^{p^{n}}$ is equal to $g^{p^{n}}$ 
for some $g\in G$. 
\item[(ii)] Moreover, for any $g_1,\ldots, g_k\in G$
there exists $g\in G$ such that $\prod\limits_{i=1}^k g_i^{p^n}=g^{p^n}$
and $g\equiv \prod\limits_{i=1}^k g_i\mod \lam_2 (\la g_1,\ldots, g_k\ra)$.
\end{itemize}
\end{itemize}
\end{Lemma}

Parts (a) and (b)(i) are well known for finitely generated pro-$p$ groups (for example, see Theorem~1.3 and Proposition~1.7 in \cite{LM1})
and can be proved quite similarly for arbitrary $(\lam,p)$-groups. Part (b)(ii) follows easily from the proof of 
(b)(i) given in \cite{LM1}.

As an immediate consequence of Lemma~\ref{lem:powerful}(b)(i) 
and its analogue for $p$-complete Lie algebras (which it trivially true), we obtain the following corollary.

\begin{Corollary}
\label{cor:cont}
Any homomorphism $\phi$ between powerful $(\lam,p)$-groups or $p$-complete Lie algebras is continuous. If $\phi$ is surjective, it is also an open map.
\end{Corollary}

We will also need several results dealing with $p$-central groups and Lie algebras.

\begin{Lemma} 
\label{lem:pc1}
Let $A$ be a group (resp. Lie ring). Then $A$ admits the unique largest $p$-central quotient, that is, there exists a normal subgroup (resp. ideal)
$N$ of $A$ such that $A/N$ is $p$-central and any homomorphism from $A$
to a $p$-central group (resp. Lie ring) factors through $A/N$.
\end{Lemma}
\begin{proof} This follows directly from the obvious fact that the class of
$p$-central groups (resp. Lie rings) is closed under subgroups (resp. subrings)
and direct products.
\end{proof}

Before stating the next result, we need to define the notion of $\Omega$-subgroups/ideals.

\begin{Definition}\rm Again let $p$ be a fixed prime.
\begin{itemize}
\item[(a)] Let $L$ be a topological Lie ring. For each $k\in\dbN$ define
$\Omega_k(L)=\{z\in L: p^k z=0\}$ (obviously $\Omega_k(L)$ is an ideal of $L$).
\item[(b)] Let $G$ be a topological group. For each $k\in\dbN$ define
$\Omega_k(G)$ to be the closed subgroup generated by
$\{g\in G: g^{p^k}=1\}$.
\end{itemize}
\end{Definition}
\skv
Note that by definition a topological group $G$ (resp. a topological Lie algebra $L$)
is $p$-central if and only if $\Omega_1(G)\subseteq Z(G)$ 
(resp. $\Omega_1(L)\subseteq Z(L)$). 
The next result provides some information about the sets $\Omega_k$ 
in $p$-central $(\lam,p)$-groups and Lie algebras.

\begin{Lemma}$\empty$
\label{lem:pcentral}
\begin{itemize}
\item[(a)] Let $L$ be a $p$-central topological Lie algebra.
Then $[\Omega_k(L),L]\subseteq \Omega_{k-1}(L)$ for all $k\geq 2$.

\item[(b)] Let $G$ be a powerful $p$-central $(\lam,p)$-group. Then
$[\Omega_2(G),G]\subseteq \Omega_{1}(G)G^{p^2}$.
\end{itemize}
\end{Lemma}
\begin{proof} (a) Take any $x\in \Omega_k(L)$. Then $p^{k-1}x\in \Omega_1(L)$. Since
$L$ is $p$-central, for any $y\in L$ we have $p^{k-1}[x,y]=[p^{k-1}x,y]=0$ and hence $[x,y]\in \Omega_{k-1}(L)$. 

(b) Using the commutator identity $[h_1 h_2,g]=[h_1,g]^{h_2}[h_2,g]$, we are reduced to showing that
$[h,g]\in \Omega_{1}(G)G^{p^2}$ whenever $h^{p^2}=1$.

So let us assume that $h^{p^2}=1$ and take any $g\in G$. Since $G$ is $p$-central, we have $[h^p,g]=1$.
On the other hand, by \eqref{eq:HP} (which is a direct consequence of the Hall-Petrescu formula) we have
$[h^p,g]\equiv [h,g]^p\mod \lam_4 G$ (since $p\geq 3$). Since $G$ is powerful, parts (a) and (b)(i)
of Lemma~\ref{lem:powerful} imply that $[h,g]^p= z^{p^3}$ for some $z\in G$. On the other hand, by 
Lemma~\ref{lem:powerful}(b)(ii) we have $1=[h,g]^p (z^{-p^2})^p=([h,g]z^{-p^2} r)^p$ for some  
$r\in \lambda_2(\lambda_2(G))=\lambda_3(G)= G^{p^2}$. Thus, $[h,g]z^{-p^2} r\in \Omega_1(G)$
and so $[h,g]\in \Omega_1(G) G^{p^2}$, as desired.
\end{proof}

\begin{Lemma}
\label{lem:pcentral2}
Let $G$ be a powerful $p$-central $(\lam,p)$-group, let $x,y\in G$, and suppose that $x^p=y^p$. Then
$x^{-1}y\in \Omega_1(G)$, that is, $(x^{-1}y)^p=1$.
\end{Lemma}
\begin{proof} Since $\Omega_1(G)$ is closed with respect to the $(\lam,p)$-topology on $G$, it suffices
to prove that $x^{-1}y\in\lam_k(G)\Omega_1(G)$ for each $k\in\dbN$. We will argue by induction on $k$, with the base case
$k=1$ being obvious.

Suppose now that $x^{-1}y\in\lam_{k}(G)\Omega_1(G)$ for some $k$. Thus, $y=xwz$ where  $w\in\lam_k(G)$ and $z\in \Omega_1(G)$.
Since $z^p=1$ and $z$ is central in $G$, we have $y^p=(xw)^p$, and by the Hall-Petrescu formula
$(xw)^p\equiv x^p w^p\mod \lam_{k+2}(G)$ (since $w\in\lam_k G$).

Since $x^p=y^p$ by assumption and $G$ is powerful, Lemma~\ref{lem:powerful}(b)(i) implies that $w^p=h^{p^{k+1}}$ for some $h$ in $G$.
By Lemma~\ref{lem:powerful}(b)(ii), $(wh^{-p^k}r)^p=1$ for some $r\in \lam_2(\la w,h^{p^k}\ra)\subseteq \lam_{k+1}(G)$.
Thus, $w\in r^{-1} h^{p^k}\Omega_1(G)\subseteq \lam_{k+1}(G)\Omega_1(G)$, which completes the induction step.
\end{proof}

\subsection{Constructing homomorphisms from $\calL_3(X)$ and $\calG_3(X)$}
In our last subsection we will prove  
Theorem~\ref{LA:pseudofree} and its group-theoretic counterpart
Theorem~\ref{thm:ppc}. The proofs of these results will be quite similar to each other, but the group case will involve some additional technicalities. Thus we will give a complete proof of 
Theorem~\ref{LA:pseudofree} and then explain which steps require non-trivial
modifications in the group case.

We start by recalling the statement of Theorem~\ref{LA:pseudofree}.
\begin{TheoremA.2}
Let $X=\{x_1,\ldots, x_d\}$ for some $d\in\dbN$. Let $L\in Ob(\grL)$ or $L=\calL_3(X)$ and $\bfu=(u_1,\ldots, u_d)$ any $d$-tuple in $L$. 
Then there exists a unique continuous Lie algebra homomorphism $f=f_{\bfu}:\calL_{3}(X)\to L$ such that
$f(c)=c_L(\bfu)$ for every 3-integral commutator $c$. In particular, $f(x_i)=u_i$ for all $1\leq i\leq d$.
\end{TheoremA.2}

\begin{proof} 
Recall that we only need to prove the existence of $f$ in the case $L\in Ob(\grL)$ (the other parts of the statement are easy and have already been established). The proof will be divided into several steps. 

{\it Step 0: a brief outline}.
For a monic $X$-commutator $c$ we set $\overline{c}=\frac{1}{p^{wt(c)-1}}c$,
and let $S\subseteq \calL_1(X)$ be the set of all elements of the
form $\overline{c}$. Thus, $$\textstyle
S=X\cup\{\frac{1}{p}[x_i,x_j]\}\cup
\{\frac{1}{p^2}[[x_i,x_j],x_k],\frac{1}{p^2}[x_i,[x_j,x_k]]\}\cup\ldots$$
Since each $\overline{c}$ is $1$-integral, it yields the corresponding function $\overline{c}_L: L^{\times d}\to L$
which is well defined up to $\Omega_1(L)=\{z\in L: pz=0\}$. From now on we fix such a function for each $\overline{c}\in S$.
If $wt(c)=1$, so $\overline c=c=x_i$ for some $i$, we naturally require that $\overline{c}((v_1,\ldots,v_d))=v_i$.
Define the function $f_0:S\to L$ by $$f_0(s)=s_L(\bfu) \mbox{ for all } s\in S.$$

Note that the $\dbZ_p$-span of $S$ is dense in $\calL_1(X)$. If it were
possible to extend $f_0$ to a Lie algebra homomorphism 
$f:\calL_1(X)\to L$, then the restriction of $f$ to $\calL_3(X)$
would be a map with desired properties; however, one can show that 
in general it is not even possible to extend $f_0$ to a $\dbZ_p$-linear map on $\calL_1(X)$. We will overcome this problem as follows.

We will construct a $\dbZ_p$-Lie algebra $\calF(X)$ generated by a set $\Shat$
 and continuous Lie algebra homomorphisms
$\widetilde f:\calF(X)\to L$ and $\pi:\calF(X)\to \calL_1(X)$
such that $\pi(\Shat\cup -\Shat)=S$ and $f_0(\pi(s))=\widetilde f (s)$ for 
all $s\in \Shat$. The Lie algebra $\calF(X)$ will be obtained from a certain Lie algebra
defined by generators and relations by first taking the $(\lam,p)$-completion and
then considering its largest $p$-central quotient.

We will then show that there exists an open subalgebra
$\calF_3(X)$ of $\calF(X)$ such that $\pi$ is injective
on $\calF_3(X)$ and maps $\calF_3(X)$ onto $\calL_3(X)$. Thus,
there exists a Lie algebra isomorphism 
$\pi^{-1}:\calL_3(X)\to \calF_3(X)$, and $\pi^{-1}$ is  automatically continuous.
since $\pi$ is open by Corollary~\ref{cor:cont}.
It is now straightforward to check that 
$f=\widetilde f \circ \pi^{-1}: \calL_3(X)\to L$ is a homomorphism with
desired properties.
\skv
{\it Step 1: Construction of $\calF(X)$.}
Choose any total ordering on the set of monic $X$-commutators
such that a longer commutator is always larger than a shorter one.
Let $T$ be the unique set of monic $X$-commutators such that
$X\subseteq T$ and an $X$-commutator $c$ of weight $\geq 2$ lies in $T$
if and only if $c=[t,t']$ with $t,t'\in T$ and $t>t'$. 

Let $\Shat$ be the set of formal symbols $\{s_t: t\in T\}$,
and let $F(X)$ be the $\dbZ_p$-Lie algebra generated by
$\Shat$ subject to relations $ps_{[c,d]}=[s_c,s_d]$ whenever
$c,d,[c,d]\in T$.

Recall that if $e$ is a monic $X$-commutator, then $\overline{e}=\frac{1}{p^{wt(e)-1}}e$.
Note that $p\overline{[c,d]}=[\overline c,\overline d]$ and $p\overline{[c,d]}_L=[\overline c_L, \overline d_L]$ (as functions on $L$)
for all monic $X$-commutators $c,d$. Hence there exist (unique) Lie algebra homomorphisms
 $\pi:F(X)\to \calL_1(X)$ and $\widetilde f: F(X)\to L$
such that $\pi(s_t)=\overline {t}$ and
$\widetilde f (s_t)=\overline t_L(\mathbf u)=f_0(\overline t)$ for all $t\in T$.
\skv

We define $\calF(X)$ as the Lie algebra obtained from 
$F(X)$ by taking first the $(\lam,p)$-completion and then
the largest $p$-central quotient (see Lemma~\ref{lem:pc1}). Since 
$L$ and $\calL_1(X)$ are both $p$-complete Lie algebras,
the homomorphisms  $\pi:F(X)\to \calL_1(X)$ and $\widetilde f: F(X)\to L$ induce the corresponding homomorphisms from
$\calF(X)$. For simplicity of notations we will denote
these induced homomorphisms by the same symbols $\pi$ and $\widetilde f$, and the image of $\Shat$ in $\calF(X)$ will also be denoted by $\Shat$.
\skv
Note that $\calF(X)$ is $p$-central (by construction) and powerful since
$F(X)$ is powerful (also by construction).
\skv

{\it Step 2:} Let $S_X=\{s_x\in \Shat: x\in X\}$ be the subset
of $\Shat$ corresponding to $X$. 
Let $\la S_X\ra$ be the abstract $\dbZ_p$-Lie subalgebra of $\calF(X)$ generated by $S_X$. 
Let
$\calF_3(X)=p^2\calF(X)+\la S_X\ra$. 
We claim that
\begin{itemize}
\item[(i)] $\pi$ is injective on $p^2\calF(X)$
\item[(ii)] $\pi$ is injective on $\la S_X\ra$
\item[(iii)] $\pi(p^2\calF(X)\cap\la S_X\ra)=
\pi(p^2\calF(X))\cap \pi(\la S_X\ra)$
\item[(iv)] $\pi(\calF_3(X))=\calL_3(X)$.
\end{itemize}

We will need the following elementary result:

\begin{Lemma}
\label{injective} Let $\phi:C\to D$ be a homomorphism of groups,
and let $A$ and $B$ be subgroups of $C$ such that $\phi$ is injective
on both $A$ and $B$ and $\phi(A\cap B)=\phi(A)\cap \phi(B)$. Then
$\phi$ is injective on $AB$.
\end{Lemma}
\begin{Remark} Lemma~\ref{injective} automatically applies to homomorphisms
of Lie rings as well.
\end{Remark}
\begin{proof} Suppose $\phi(ab)=1$ for some $a\in A$, $b\in B$. Then
$\phi(a)=\phi(b^{-1})\in \phi(A)\cap \phi(B)=\phi(A\cap B)$. Since
$\phi$ is injective on $A$ and $B$, it follows that $a,b\in A\cap B$
whence $ab\in A\cap B$, so $\phi(ab)=1$ forces $ab=1$.
\end{proof}

If we prove (i)-(iv), it would follow from Lemma~\ref{injective}
that $\calF_3(X)=p^2\calF(X)+\la S_X\ra$ satisfies all the requirements described in Step~0 and
thus finish the proof.

It is clear that $\pi(\la S_X\ra)=\la X\ra$ (the abstract 
$\dbZ_p$-Lie subalgebra of $\calL_1(X)$ generated by $X$).
This immediately implies (ii) since $\la X\ra$ is free on $X$ and $|S_X|=|X|$.
Conditions (iii) and (iv) can be checked by an easy direct computation.
Both sides of (iii) consist of all elements of $\la X\ra$ in which the coefficients of all degree $1$ monomials are divisible by $p^2$ and the
coefficients of all degree $2$ monomials are divisible by $p$.
Thus, it remains to prove (i).

\skv 
{\it Step 3:} In this step we reduce condition (i) from Step~2 to showing that
$\pi$ induces an injective (and hence bijective) map 
$\pi_2: p^2\calF(X)/p^3 \calF(X)\to p^2\calL_1(X)/p^3 \calL_1(X)$.

For each $k\in\dbN$ we have a commutative diagram
\begin{equation}
\label{CD:ppc}
\begin{CD}
p^k\calF(X)/p^{k+1} \calF(X)@ >\pi_k >> p^k\calL_1(X)/p^{k+1} \calL_1(X)\\
@ VV\theta_k V @ VV\theta'_k V \\
p^{k+1}\calF(X)/p^{k+2} \calF(X) @ >\pi_{k+1}>> 
p^{k+1}\calL_1(X)/p^{k+2} \calL_1(X)\\
\end{CD}
\end{equation}
where the horizontal maps $\pi_k$ and $\pi_{k+1}$ are induced by
$\pi$ and the vertical maps $\theta_k$ and $\theta_k'$ are induced by
multiplication by $p$.

The maps $\theta_k$ and $\theta_k'$ are automatically surjective. Moreover,
$\theta_k'$ is also injective since $\calL_1(X)$ is torsion-free. Therefore, 
\eqref{CD:ppc} shows that injectivity of $\pi_k$ for some $k$ implies
injectivity of $\pi_{k+1}$. Thus if we know that $\pi_3$
is injective, then $\pi_k$ is injective for all $k\geq 3$ and since
$\bigcap\limits_k p^k\calF(X)=0$ (as $\calF(X)$ is a $p$-complete Lie algebra), it follows that $\pi$ is injective on $p^2\calF(X)$, as desired.
\skv
{\it Step 4:} We are now reduced to proving that the map 
$$\pi_2:p^2\calF(X)/p^3 \calF(X)\to p^2\calL_1(X)/p^3 \calL_1(X)$$ from Step~3 is injective. 

Let $B$ be Hall's set of basic $X$-commutators. We will recall the definition
of $B$ in Step~5 below; at this point all we need to know is that $B\subseteq T$
and $B$ forms a basis for the free $\dbZ$-Lie algebra on $X$ (and hence also a basis for the free $R$-Lie algebra on $X$ for any commutative ring $R$ with $1$).
\skv

Let $M$ be the abstract $\dbZ_p$-submodule of $\calF(X)$ generated by $\{s_b: b\in B\}$. The image of $p^2M$ in $p^2\calF(X)/p^3 \calF(X)$ is an $\dbF_p$-vector space
spanned by $\{p^2 s_b+p^3 \calF(X): b\in B\}$ and  
$\pi_2 (p^2 s_b)= \frac{1}{p^{wt(b)-3}}b+p^3 \calL_1(X)$ for all $b\in B$. 
\skv
The elements
$\{\frac{1}{p^{wt(b)-3}}b+p^3 \calL_1(X): b\in B\}$ are clearly linearly independent
(in fact, if $\mathrm{Lie}_{\dbF_p}(X)$ is the free $\dbF_p$-Lie algebra on $X$,
then $\mathrm{Lie}_{\dbF_p}(X)\cong p^2\calL_1(X)/p^3 \calL_1(X)$ as 
$\dbF_p$-vector spaces via the map $b\mapsto \frac{1}{p^{wt(b)-3}}b+p^3 \calL_1(X)$
for all $b\in B$), so $\pi_2$ is injective on 
$(p^2M+p^3\calF(X))/p^3 \calF(X)$. Thus, to finish the proof it suffices
to show that $p^2\calF(X)=p^2M+p^3\calF(X)$ (where only the inclusion
``$\subseteq$'' needs to be justified).
\skv

{\it Step 5:} For $k\in\dbN$ let $\Omega_k=\Omega_k(\calF(X))=\{z\in \calF(X): p^k z=0\}$. In this step we will prove the inclusion 
\begin{equation}
\label{eq:omega}
\calF(X)\subseteq M+\Omega_2 + p\calF(X). 
\end{equation}
Multiplying both sides of \eqref{eq:omega} by $p^2$, we get
$p^2\calF(X)\subseteq p^2M+p^3\calF(X)$, which would finish the proof
by Step~4.
\skv

Since $\calF(X)$ is powerful, to prove \eqref{eq:omega} it suffices
to show that $s_t\in M+\Omega_2+p\calF(X)$ for any $t\in T$. We will split the
proof into two parts.
\skv
{\it Substep 1:} First we will show that $s_t\in M+\Omega_2$ for all
$t\in T$ such that $t=[b,c]$ for some $b,c\in B$
(our assumptions on $T$ imply that $b>c$). 
\skv
Recall that $B$ is Hall's set of basic $X$-commutators (relative to the chosen ordering on the set of all monic $X$-commutators) which are defined as follows:
\begin{itemize}
\item[(i)] $B$ contains $X$ (all monic $X$-commutators of weight $1$).
\item[(ii)] Suppose $k>1$ and we already defined which monic $X$-commutators
of weight $<k$ lie in $B$. Then a monic $X$-commutator $c$ of weight $k$ lies in $B$
if $c=[y,z]$ where $y,z\in B$, $y>z$ and either $wt(y)=1$ or
$y=[v,w]$ with $w\leq z$.
\end{itemize}

By construction $B\subseteq T$ and all elements of $T$ of weight $\leq 2$ lie in $B$, so we only need to prove the assertion of Substep~1 for $t$ with $wt(t)\geq 3$.

So assume that $t=[b,c]$ with $b,c\in B$, $b>c$ and $wt(t)=wt(b)+wt(c)\geq 3$.
Then $wt(b)>1$, so $b=[d,e]$ for some $d,e\in B$ with $d>e$. We consider 3 cases.
In Case~3 below we will argue by downward induction on $wt(c)$.
It is clear that Case~3 can only occur if  $wt(c)\leq \frac{k}{3}$,
so all $c$ with $wt(c)>\frac{k}{3}$ are covered by Cases~1, 2 and 2' where
we give a direct argument.
\skv

{\it Case 1: $e\leq c$}. In this case 
$t=[b,c]=[d,e,c]\in B$, so $s_t\in M$. 
\skv
{\it Case 2: $e>c$ and $[e,c]> d$}. In this case $[e,c,d]\in B$
and $[d,c,e]\in B$, so $s_{[e,c,d]}, s_{[d,c,e]}\in M$. Moreover,
since $[d,e,c]$, $[e,c,d]$ and $[d,c,e]$ all lie in $T$, we have
\begin{multline*}
0=[s_d,s_e,s_c]+ [s_e,s_c,s_d]+[s_c,s_d,s_e]\\
=
[s_d,s_e,s_c]+ [s_e,s_c,s_d]-[s_d,s_c,s_e]=
p^2(s_{[d,e,c]}+s_{[e,c,d]}-s_{[d,c,e]}).
\end{multline*}
Hence, $s_{[d,e,c]}+s_{[e,c,d]}-s_{[d,c,e]}\in \Omega_2$,
so $$s_t=s_{[d,e,c]}=s_{[d,c,e]}-s_{[e,c,d]}+
(s_{[d,e,c]}+s_{[e,c,d]}-s_{[d,c,e]})\in M+\Omega_2.$$
\skv
{\it Case 2':} $e>c$ and $[e,c]=d$. In this case $[e,c,d]=0$ and
$[s_e,s_c,s_d]=p[s_{[e,c]},s_d]=0$. We can essentially repeat
the argument from Case 2 (just delete $s_{[e,c,d]}$ from the above
calculations).  
\skv
{\it Case 3:} $e>c$ and $[e,c]<d$. In this case
$[d,c,e]\in B$ and $[d,[e,c]]\in T$, so
\begin{multline*}
0=[s_d,s_e,s_c]+ [s_e,s_c,s_d]+[s_c,s_d,s_e]\\
=
[s_d,s_e,s_c]- [s_d,[s_e,s_c]]-[s_d,s_c,s_e]=
p^2(s_{[d,e,c]}-s_{[d,[e,c]]}-s_{[d,c,e]})
\end{multline*}
and thus $s_{[d,e,c]}-s_{[d,[e,c]]}-s_{[d,c,e]}\in \Omega_2$.
We also have $s_{[d,c,e]}\in M$. While
$s_{[d,[e,c]]}$ may not lie in $M$, since 
$wt([e,c])>wt(c)$, we know that $s_{[d,[e,c]]}\in M+\Omega_2$
by the induction hypothesis, and we are done as in case~2.
\skv

Thus, we have now proved that
$s_t\in M+\Omega_2$ for all $t\in T$ such that $t=[b,c]$ for some $b,c\in B$.
\skv
{\it Substep 2:} Let us now prove that $s_t\in M+\Omega_2+p\calF(X)$ for all 
$t\in T$ by induction on $k=wt(t)$.
\skv
As in substep 1, there is nothing to prove for $k\leq 2$.
So fix $k\geq 3$, $t\in T$ with $wt(t)=k$ and suppose that
$s_{t'}\in M+\Omega_2 + p\calF(X)$ for all $t'\in T$ with $wt(t')<k$.

Write $t=[t_1,t_2]$ with $t_1,t_2\in T$ and $t_1>t_2$.
By the induction hypothesis we can write 
$s_{t_1}=\sum\limits_{b\in B} n_b s_b+py_1+r_1$ and 
$s_{t_2}=\sum\limits_{b\in B} m_b s_b+py_2+r_2$
where $n_b,m_b\in\dbZ_p$, $y_1,y_2\in\calF(X)$ and $r_1,r_2\in \Omega_2$.
Therefore (using again that $\calF(X)$ is powerful) we have
\begin{multline}
\label{eq:final} 
ps_t=[s_{t_1},s_{t_2}]= \sum\limits _{b,b'\in B}n_b m_{b'}[s_b,s_{b'}]+
p^2 y+pz\\
=
\sum\limits_{b>b'\in B}(n_b m_{b'}-m_b n_{b'})p s_{[b,b']}+
p^2 y+pz.
\end{multline}
where $y\in \calF(X)$ and $pz\in [\Omega_2,\calF(X)]$.

Equation \eqref{eq:final} implies that
$s_t-(\sum\limits_{b>b'\in B}(n_b m_{b'}-m_b n_{b'}) s_{[b,b']}+
py+z)\in \Omega_1$.
By Substep~1, $s_{[b,b']}\in M+\Omega_2$ for all $b,b'\in B$
with $b>b'$. Finally, by Lemma~\ref{lem:pcentral}(a)
$pz\in \Omega_1$, so $z\in \Omega_2$. Thus, $s_{t}\in M+\Omega_2+p\calF(X)$, as desired.
\skv
The proof of Theorem~\ref{LA:pseudofree} is now complete.
\end{proof}

Let us finally prove Theorem~\ref{thm:ppc}. As with Theorem~\ref{LA:pseudofree}, we recall the statement.

\begin{TheoremA.9}
Let $X=\{x_1,\ldots, x_d\}$ and $G\in Ob(\grG)$ or $G=\calG_3(X)$ .
For any finite subset $\{g_1,\ldots,g_d\}$ of $G$ there exists a continuous group homomorphism $f:\calG_{3}(X)\to G$ such that $f(\Exp(x_i))=g_i$ for all $i$. 
\end{TheoremA.9}

\begin{proof}
As we already mentioned, this proof is very similar to that of
Theorem~\ref{LA:pseudofree}. We will briefly go over the argument concentrating on the non-obvious changes. 

The role of $\calL_1(X)$ will be played by the group
$\calG_1(X)$. Recall that $\calG_1(X)$ is powerful and torsion-free, and moreover
any element of $\calG_1(X)^p$ has a unique $p^{\rm th}$ root.
\skv

For each monic $X$-commutator $c$ let $h_c=H_c(\Exp(x_1),\ldots, \Exp(x_d))\in \calG_1(X)$
where $H_c$ is defined as at the beginning of \S~A.3.
 
\skv
{\it Step 1:} We define $T$ exactly as in the Lie algebra case, and again we let 
$\Shat=\{s_t: t\in T\}$. We define $F(X)$ to be the abstract group
generated by $\Shat$ subject to relations $s_{[c,d]}^p=[s_c,s_d]$ whenever
$c,d,[c,d]\in T$.

As in the Lie algebra case, we have (unique) homomorphisms with dense images
$\pi:F(X)\to \calG_1(X)$ and $\widetilde f: F(X)\to G$ such that
$\pi(s_t)=h_t$ and $\widetilde f(s_t)=H_t(g_1,\ldots, g_d)$ for all $t\in T$
(recall from \S~A.3 that
the values of $H_t: G^{\times d}\to G$ are well defined modulo $\Omega_1(G)$).

If $\calF(X)$ is the group obtained from $F(X)$ by taking first the $(\lam,p)$-completion and then taking the largest $p$-central quotient, we obtain
induced continuous surjective homomorphisms 
$\pi:\calF(X)\to \calG_1(X)$ and $\widetilde f: \calF(X)\to G$.
\skv
{\it Step 2:} Let $\calF_3(X)=\calF(X)^{p^2}\la S_X\ra$ where
$S_X$ is the abstract subgroup generated by $S_X=\{s_x: x\in X\}$.
As in the Lie algebra case we reduce Theorem~\ref{thm:ppc} 
first to showing that $\pi$ is injective on $\calF_3(X)$
and then to showing that $\pi$ is injective on $\calF(X)^{p^2}$.

\skv
{\it Step 3:} Now we reduce to proving that the induced map
$\pi_2:\calF(X)^{p^2}/\calF(X)^{p^3}\to \calG_1(X)^{p^2}/\calG_1(X)^{p^3}$ is injective.
Similarly to the Lie algebra case this is done via the commutative diagram

\begin{equation}
\label{CD:ppc2}
\begin{CD}
\calF(X)^{p^k}/\calF(X)^{p^{k+1}} @ >\pi_k >> \calG_1(X)^{p^k}/\calG_1(X)^{p^{k+1}}\\
@ VV\theta_k V @ VV\theta'_k V \\
\calF(X)^{p^{k+1}}/\calF(X)^{p^{k+2}} @ >\pi_{k+1}>> 
\calG_1(X)^{p^{k+1}}/\calG_1(X)^{p^{k+2}} \\
\end{CD}
\end{equation}
where the horizontal maps $\pi_k$ and $\pi_{k+1}$ are induced by
$\pi$ and the vertical maps $\theta_k$ and $\theta_k'$ are induced by
raising to power $p$. Since $\calF(X)$ and $\calG_1(X)$ are powerful 
$(\lam,p)$-groups, the maps $\theta_k$ and $\theta_k'$ are surjective
by Lemma~\ref{lem:powerful}. Since every element of $\calG_1(X)^p$
has a unique $p^{\rm th}$-root, Lemma~\ref{lem:powerful} also implies that
$\theta'_k$ is injective.
\skv
{\it Step 4:} Define $M$ to be the abstract
subgroup generated by $\{s_b: b\in B\}$ (where $B$ is still Hall's set of basic $X$-commutators). We claim that
$\pi_2$ is injective on $M^{p^2}\calF(X)^{p^3}/\calF(X)^{p^{3}}$. This reduces
to the corresponding result in the Lie algebra case using the simple observation
that the exponential map $\Exp:\calL_1(X)\to \calG_1(X)$ induces an isomorphism of
$\dbF_p$-vector spaces 
$p^k \calL_1(X)/p^{k+1} \calL_1(X)\to \calG_1(X)^{p^k}/\calG_1(X)^{p^{k+1}}$.
\skv
{\it Step 5:} As in the Lie algebra case (this time using the fact that
$\calF(X)$ and $\calG_1(X)$ are powerful), we reduce Theorem~\ref{thm:ppc}(a)
to showing that
\begin{equation}
\label{eq:omega2}
\calF(X) \subseteq M\cdot\Omega_2\cdot \calF(X)^p. 
\end{equation}
where $\Omega_2=\Omega_2(\calF(X))=\{g\in \calF(X): g^{p^2}=1\}$ (the order of the product on the right-hand
side of \eqref{eq:omega2} is irrelevant since $\Omega_2$ and $\calF(X)^p$
are both normal). 

The argument from the Lie algebra case carries over without any changes until
Case~2 of Substep~1 where in the group case we have the equality
$$
[s_d,s_e,s_c]\cdot [s_e,s_c,s_d]\cdot {[s_d,s_c,s_e]}^{-1}
={s_{[d,e,c]}}^{p^2} {s_{[e,c,d]}}^{p^2} {({s_{[d,c,e]}}^{-1})}^{p^2}.
$$
It is still true that $[e,c,d],[d,c,e]\in B$, so $s_{[e,c,d]}, s_{[d,c,e]}\in M$.
Unlike the Lie algebra counterpart of this equality, the left-hand side
need not be trivial; however the element $y=[s_d,s_e,s_c]\cdot [s_e,s_c,s_d]\cdot {[s_d,s_c,s_e]}^{-1}$ still lies in $\gamma_4 \calF(X)$ by basic commutator identities, and hence by Lemma~\ref{lem:powerful}(b)(i)
$y=w^{p^3}$ for some $w\in \calF(X)$. Thus,
${s_{[d,e,c]}}^{p^2} {s_{[e,c,d]}}^{p^2} ({s_{[d,c,e]}}^{-1})^{p^2}(w^{-p})^{p^2}=1$,
and applying the Lemma~\ref{lem:powerful}(b)(ii) we get that
there exists $r\in \lambda_2 \calF(X)=\calF(X)^p$ such that
$(s_{[d,e,c]}s_{[e,c,d]}{s_{[d,c,e]}}^{-1}w^{-p}r)^{p^2}=1$, whence
$s_{[d,e,c]}\in s_{[d,c,e]}{s_{[e,c,d]}}^{-1}\calF(X)^p \Omega_2\subseteq
M\cdot  \calF(X)^p\cdot\Omega_2$, as desired.
\skv
Making similar modifications in Case~3 of Substep~1 and in Substep~2, we complete the proof of Theorem~\ref{thm:ppc}.
\end{proof}


\begin{thebibliography}{AAAA}
\bibitem[De]{De} V. V. Deodhar,
\it On central extensions of rational points of algebraic groups.
\rm Amer. J. Math. 100 (1978), no. 2, 303--386.

\bibitem[DDMS]{DDMS} J. D. Dixon, M. P. F. du Sautoy, A. Mann and D. Segal,
\it Analytic pro-$p$ groups.
\rm Second edition. Cambridge Studies in Advanced Mathematics, 61. Cambridge University
Press, Cambridge, 1999.

\bibitem[HZ]{HZ} M. Horn, S. Zandi,
\it Second cohomology of Lie rings and the Schur multiplier. 
\rm Int. J. Group Theory 3 (2014), no. 2, 9--20.

\bibitem[Iw]{Iw} K. Iwasawa,
\it Local class field theory. 
\rm 
Oxford Science Publications. Oxford Mathematical Monographs. 
The Clarendon Press, Oxford University Press, New York, 1986. viii+155 


\bibitem[K]{K} B. Klopsch,
\it On the Lie theory of $p$-adic analytic groups. 
\rm Math. Z. 249 (2005), no. 4, 713--730.

\bibitem[La1]{La1} M. Lazard,
\it Sur les groupes nilpotents et les anneaux de Lie. (French) 
\rm Ann. Sci. \'Ecole Norm. Sup. (3) 71 (1954), 101--190.

\bibitem[La2]{La} M. Lazard,
\it Groupes analytiques $p$-adiques. (French)
\rm Inst. Hautes \`Etudes Sci. Publ. Math. No. 26 (1965), 389--603.

\bibitem[LM1]{LM1} A. Lubotzky, A. Mann,
\it Powerful $p$-groups. I. Finite groups.
\rm J. Algebra 105 (1987), no. 2, 484--505.


\bibitem[LM2]{LM} A. Lubotzky, A. Mann,
\it Powerful $p$-groups. II. $p$-adic analytic groups.
\rm J. Algebra 105 (1987), no. 2, 506--515.

\bibitem[Ma]{Ma} H. Matsumoto,
\it Sur les sous-groupes arithm\'etiques des groupes semi-simples d\'eploy\'es. (French)
\rm Ann. Sci. \'Ecole Norm. Sup. (4) no. 2 (1969), 1--62.

\bibitem[Mo1]{Mo1} C. C. Moore,
\it Group extensions of $p$-adic and adelic linear groups.
\rm Inst. Hautes Etudes Sci. Publ. Math. No. 35 (1968), 157--222.

\bibitem[Mo2]{Mo2} C. C. Moore,
\it Group extensions and cohomology for locally compact groups. IV.
\rm Trans. Amer. Math. Soc. 221 (1976), no. 1, 35--58.

\bibitem[NS1]{NS1} N. Nikolov, D. Segal,
\it On finitely generated profinite groups. I. Strong completeness and uniform bounds. 
\rm Ann. of Math. (2) 165 (2007), no. 1, 171--238. 

\bibitem[NS2]{NS2} N. Nikolov, D. Segal,
\it On finitely generated profinite groups. II. Products in quasisimple groups. 
\rm Ann. of Math. (2) 165 (2007), no. 1, 239--273.

\bibitem[Pr1]{Pr1} G. Prasad,
\it On some work of Raghunathan. 
\rm Algebraic groups and arithmetic, 25--40, Tata Inst. Fund. Res., Mumbai, 2004. 

\bibitem[Pr2]{Pr2} G. Prasad,
\it Deligne's topological central extension is universal. 
\rm Adv. Math. 181 (2004), no. 1, 160--164.

\bibitem[PR1]{PR1} G. Prasad, M. S. Raghunathan,
\it Topological central extensions of
semisimple groups over local fields I, II.
\rm Ann. of Math. (2) 119 (1984),
no. 1, 143--201 and no. 2, 203--268.

\bibitem[PR2]{PR2} G. Prasad, M. S. Raghunathan,
\it Topological central extensions of $SL_1(D)$.
\rm Invent. Math. 92 (1988), 645--689.


\bibitem[PRp]{PRp} G. Prasad, A. Rapinchuk,
\it Computation of the metaplectic kernel. 
\rm Inst. Hautes \'Etudes Sci. Publ. Math. No. 84 (1996), 91--187 (1997).


\bibitem[Ra]{Ra} M. S. Raghunathan,
\it On the congruence subgroup problem.
\rm Inst. Hautes \'Etudes Sci. Publ. Math. No. 46 (1976), 107--161.

\bibitem[Rp]{Rp} A. Rapinchuk,
\it The multiplicative arithmetic of division algebras over number fields and the metaplectic problem. 
\rm Math. USSR-Izv. 31 (1988), no. 2, 349--379. 

\bibitem[Ri]{Ri} C. Riehm,
\it The norm 1 group of $\mathfrak{p}$-adic division algebra.
\rm Amer. J. of Math. 92 (1970), no. 2, 499-523.


\bibitem[Se]{Se} J.-P. Serre,
\it Local fields.
\rm Translated from the French by Marvin Jay Greenberg. 
Graduate Texts in Mathematics, 67. Springer-Verlag, New York-Berlin, 1979.

\bibitem[Weib]{Wb} C. Weibel,
\it An introduction to homological algebra. 
\rm Cambridge Studies in Advanced Mathematics, 38. 
Cambridge University Press, Cambridge, 1994. xiv+450 pp.



\bibitem[Weig]{We} T. Weigel,
\it Exp and log functors for the categories of powerful p-central
groups and Lie algebras,
\rm Habilitationsschrift, Freiburg, 1994.

\bibitem[Wil]{Wil} J. Wilson,
\it Profinite groups,
\rm London Mathematical Society Monographs. New Series, 19.
The Clarendon Press, Oxford University Press, New York, 1998.
\end{thebibliography}
\end{document}